\documentclass[final,3p,times,fleqn]{elsarticle}
\usepackage{mathtools,booktabs}
\usepackage{bm}
\usepackage{bbm}
\usepackage{makecell}
\usepackage{amsmath}
\usepackage{amsthm}
\usepackage{gensymb}
\usepackage{subfigure}
\usepackage{epstopdf}
\usepackage[T1]{fontenc}
\usepackage{color}
\usepackage{hyperref}

\usepackage{threeparttable}
\newtheorem{remark}{Remark}
\biboptions{sort&compress}
\begin{document}
\begin{frontmatter}
	\title{A general fourth-order mesoscopic multiple-relaxation-time lattice Boltzmann model and equivalent macroscopic finite-difference scheme for two-dimensional diffusion equations}
	\author[a]{Ying Chen}
	\author[a,b,c]{Zhenhua Chai \corref{cor1}}
	\ead{hustczh@hust.edu.cn}
	\author[a,b,c]{Baochang Shi}
	\address[a]{School of Mathematics and Statistics, Huazhong University of Science and Technology, Wuhan, 430074, China}
	\address[b]{Institute of Interdisciplinary Research for Mathematics and Applied Science, Huazhong University of Science and Technology, Wuhan 430074, China}
	 \address[c]{Hubei Key Laboratory of Engineering Modeling and
	 Scientific Computing, Huazhong University of Science and Technology,
	 Wuhan 430074, China}
	\cortext[cor1]{Corresponding author.}
\begin{abstract}
	In this work, we first develop a general mesoscopic multiple-relaxation-time lattice Boltzmann (MRT-LB) model for the two-dimensional diffusion equation with the constant diffusion coefficient and source term, where the D2Q5 (five discrete velocities in two-dimensional space) lattice structure is considered. Then we exactly derive the equivalent macroscopic finite-difference scheme of the MRT-LB model. Additionally, we also propose a proper MRT-LB model for the diffusion equation with a linear source term, and obtain an equivalent macroscopic six-level finite-difference scheme. After that, we conduct the accuracy and stability analysis of the finite-difference scheme and the mesoscopic MRT-LB model. It is found that at the diffusive scaling, both of them can achieve a fourth-order accuracy in space based on the Taylor expansion. The stability analysis also shows that they are both unconditionally stable. Finally, some numerical experiments are conducted, and the numerical results are also consistent with our theoretical analysis. 
\end{abstract}
\end{frontmatter}
\section{Introduction}
The mesoscopic lattice Boltzmann (LB) method, as one of the kinetic-theory based approaches, has received increasing attention in the past three decades. Compared to some traditional numerical approaches, the LB method is a simple and versatile numerical approach, which  not only has gained a great success in the study of complex flows governed by the Navier-Stokes equations \cite{Chen1998, Succi2001, Aidun2010, Guo2013, Kruger2017, Wang2019}, but also can be considered as a numerical solver to the partial differential equations (PDEs), including the diffusion equations \cite{Huber2010, Ancona1994, Suga2010, Lin2022, Silva2023}, the convection-diffusion equations \cite{Van2000, Gin2005, Rasin2005, Shi2009, Chopard2009, Yoshida2010, Gin2013, Chai2010, Chai2014, Jettestuen2016, Li2017,  Michelert2022, Dellacherie2014, Cui2016, Chen2023}, Poisson equation \cite{Hirabayashi2001, Chai2008, Chai2019} and Burgers equation \cite{Li2012}. Based on the collision term in the LB method, the LB models can be classified into three basic kinds, the lattice Bhatnagar-Gross-Krook (BGK) or single-relaxation-time LB (SRT-LB) model \cite{Qian1992}, two-relaxation-time LB (TRT-LB) model \cite{Gin2005} and multiple-relaxation-time LB (MRT-LB) model \cite{Chai2020}. Although the SRT-LB model is the most efficient, it suffers from the stability and/or accuracy problem when the relaxation parameter is close to its limit value \cite{Chai-1-2014, Yoshida2010}, and what is more, it is usually limited to study the isotropic diffusion problems \cite{Guo2013, Lallemand2000} since it does not have sufficient relaxation parameters to describe the anisotropic diffusion process. To overcome these inherent deficiencies in the SRT-LB model, the TRT-LB model with a free relaxation parameter that can be used to improve the stability and/or accuracy by selecting the so-called magic parameter $\Lambda^{eo}$ properly \cite{Kruger2017}. Compared to the SRT-LB and the TRT-LB models, the MRT-LB model is more general, and could be more stable through adjusting additional free relaxation parameters \cite{Lallemand2000, Pan2006, Cui2016, Luo2011}. For this reasons, here we only consider the MRT-LB model.

Up to now, although some asymptotic analysis approaches, including the Chapman-Enskog analysis \cite{Chapman1990}, Maxwell iteration \cite{Ikenberry1956, Yong2016}, direct Taylor expansion \cite{Holdych2004, Wagner2006}, recurrence equations \cite{Gin2012,GinzburgD2009} and equivalent equations \cite{Dubois2008, Dubois2009, Dubois2019}, can be used to illustrate that the LB method is suitable for the specified macroscopic PDEs, the relation between the LB model and the macroscopic PDEs based numerical schemes (hereafter named macroscopic numerical schemes) is still unclear, which means that it is inconvenient or difficult to construct a rigorous notion of consistency and perform an accuracy analysis of the LB model for the specified PDE.

In the past years, some researchers have made more contributions on this aspect. For instance, Junk \cite{Junk} and Inamuro \cite{Inamuro} found that the SRT-LB model is equivalent to a macroscopic two-level finite-difference scheme if the relaxation parameter is equal to unity, and at the diffusive scaling, a second-order convergence rate in space can be achieved for the incompressible Naiver-Stokes equations \cite{Junk}. Li et al. \cite{Li2012} demonstrated that for one-dimensional Burgers equation, the SRT-LB model with the D1Q2 lattice structure can be written as a macroscopic three-level second-order finite-difference scheme. In addition to above works, more researchers pay attention to the macroscopic numerical schemes of the LB models for the diffusion and convection-diffusion equations. On the one hand, for one-dimensional diffusion equations, Ancona \cite{Ancona1994} first pointed out the SRT-LB model with the D1Q2 lattice structure is consistent with a macroscopic three-level second-order DuFort-Frankel scheme \cite{Fort1953}. Suga \cite{Suga2010} found that the SRT-LB model with the D1Q3 lattice structure could be corresponding to a macroscopic four-level fourth-order finite-difference scheme, and then Lin et al. \cite{Lin2022} considered a more general MRT-LB model, and derived an equivalent macroscopic four-level sixth-order finite-difference scheme. Silva \cite{Silva2023} focused on the TRT-LB model for the diffusion equation with a linear source term, and obtained a macroscopic fourth-order finite-difference scheme. On the other hand, for one-dimensional convection-diffusion equations, Dellacherie \cite{Dellacherie2014} carried out an analysis on the SRT-LB model with D1Q2 lattice structure, and demonstrated that this LB model is equivalent to a three-level finite-difference scheme named LFCCDF schemes \cite{Kwok1993} (LFCCDFS means that the Leap-Frog difference is used for the temporal derivative, the cenetral difference is adopted for the convection term and the Du Fort-Frankel approximation is applied for the diffusion term). Cui et al. \cite{Cui2016} illustrated that for one-dimensional steady convection-diffusion equation, the MRT-LB model can be written as a macroscopic second-order centeral-difference scheme. Recently, we developed a fourth-order MRT-LB model with the D1Q3 lattice structure, and also obtained an equivalent four-level fourth-order finite-difference scheme from the MRT-LB model \cite{Chen2023}. However, it should be noted that all the works mentioned above are limited to the one-dimensional case. Recently, Fu\v{c}\'{i}k et al. \cite{Fucik2021} focused on a decomposition of the LB model by using a hollow matrix, and deveoped an algorithm to obtain the equivalent macroscopic finite-difference scheme.	They also presented a computational tool for the derivation of equivalent PDEs from the LB model with a general collision operator \cite{Fucik2023}. However, the algorithm is relatively complicated. Bellotti et al. \cite{Bellotti2022, Bellotti-1-2022} presented a precise algebraic characterization of LB models rather than a purely algorithmic approach developed by Fu\v{c}\'{i}k et al. \cite{Fucik2021, Fucik2023}, and performed some detailed studies on the relation between the MRT-LB model and the macroscopic numerical scheme. They also demonstrated that any LB models with DdQq (q discrete velocites in d-dimensional space) lattice structure can be rewritten as the macroscopic second-order finite-difference schemes, which is consistent with the results based on the asymptotic analysis methods \cite{Chai2020}. Besides, d$'$Humi$\grave{{\rm e}}$re and Ginzburg \cite{GinzburgD2009} also conducted a theoretical analysis on the TRT-LB model with recurrence equations, and illustrated that when the magic parameter  satisfies the condition $\Lambda^{eo} = 1/4$, the TRT-LB model would reduce to a macroscopic three-level finite-difference scheme with a second-order accuracy in space. Subsequently, Ginzburg \cite{Gin2012} further found that the TRT-LB model could achieve a fourth-order accuracy in space for diffusion equations with $\Lambda^{eo} = 1/6$ and $\Lambda^{BGK} = 1/12$, and a third-order accuracy in space for convection-diffusion equations with $\Lambda^{eo} = \Lambda^{BGK} = 1/12$. However, some details on the TRT-LB model (e.g., the equilibrium distribution function) have not been provided in this work \cite{Gin2012}. Now the question is whether we can obtain an equivalent high-order finite-difference scheme for the two-dimensional diffusion equation from the more general MRT-LB model with the simple D2Q5 lattice structure? To this end, we first conduct a theoretical analysis on the MRT-LB model with the D2Q5 lattice structure, and obtain an equivalent macroscopic finite-difference scheme of the MRT-LB model. Then we also show that the TRT-LB model with $\Lambda^{eo} = 1/6$ and $\Lambda^{BGK} = 1/12$ \cite{Gin2012} is just a special case of our work.

The rest of this paper is organized as follows. In Sec. \ref{sec2}, a general MRT-LB model with the D2Q5 lattice structure is developed for the two-dimensional diffusion equation, followed by the equivalent macroscopic six-level finite-difference (SLFD) scheme of the MRT-LB model and a simplified five-level finite-difference (FLFD) scheme in Sec. \ref{sec3}. In Sec. \ref{sec4}, we perform the accuracy and stability analysis of the MRT-LB model and equivalent finite-difference scheme with the Taylor expansion and von Neumann stability analysis method, and demonstrate that both of them are not only fourth-order accurate in space, but also unconditionally stable. In Sec. \ref{secNumer}, some numerical experiments are carried out, and the results show that both the MRT-LB model and finite-difference scheme can achieve a fourth-order convergence rate in space, which is consistent with our theoretical analysis. Finally, some conclusions are given in Sec. \ref{sec6}.

	\section{The multiple-relaxation-time lattice Boltzmann model for two-dimensional diffusion equations}\label{sec2}
In this section, we first present the MRT-LB model for the two-dimensional diffusion equation with the constant diffusion coefficient and source term, and then give the condition which the four-order MRT-LB model should satisfy.
\subsection{The two-dimensional diffusion equation}
From the mathematical point of view, the two-dimensional governing equation for the diffusion process of mass and heat or any other scalar quantity can be written as
\begin{equation}\label{equde}
	\partial_t \phi= \nabla \cdot(\kappa\nabla \phi)+R,
\end{equation}	
where $\phi(\bm{r},t)$ is a scalar variable dependent on the space $\bm{r}=(x,y)^T$ and time $t$. $\kappa$ is the diffusion coefficient, $R$ is the source term, and they are assumed to be constant in this work.

	\subsection{Spatial and temporial discretization}
To simplify the following analysis, the two-dimensional space in the LB method
is discretized by the lattice $\mathcal{L}_x:=\Delta x \mathbbm{Z}, \mathcal{L}_y:=\Delta x \mathbbm{Z}$ with the lattice spacing $\Delta x=\Delta y>0$, the time is uniformly discretized with $t_{n}:=n\Delta t$, $\Delta t$ is the time step. Then the so-called lattice velocity is defined by $c=\Delta x/\Delta t$.
\subsection{The multiple-relaxation-time lattice Boltzmann model}
Here we only consider the more general MRT-LB model for its good accuracy and stability in the study of complex problems \cite{Chai-1-2016, Chai2020, Chai2014, Lallemand2000, Chai2010}. For the two-dimensional diffusion equation (\ref{equde}), the evolution equation of the MRT-LB model can be written as \cite{Chai2014, Cui2016}
\begin{equation}\label{distributionf}
	\begin{aligned}
		&f_i(\bm{r}+\bm{c_i}\Delta t, t+\Delta t)=f_i(\bm{r}, t)-\big(\bm{M}^{-1}\bm{S}\bm{M}\big)_{i,k}
		[f_k(\bm{r},t)-f^{eq}_k(\bm{r},t)]+\Delta t\big[\bm{M}^{-1}\big(\bm{I}-
		\frac{\bm{S}}{2}\big)\bm{M}\big]_{i,k}R_k,(i=0,1,2,3,4),
	\end{aligned}
\end{equation}
where $f_i(\bm{r},t)$, $f^{eq}_{i}(\bm{r},t)$ and $R_i(\bm{r},t)$ are the distribution function, the equilibrium distribution function and the discrete source term at position $\bm{r}$ and time $t$, respectively. In the D2Q5 lattice structure, the discrete velocity $\bm{c_{i}}$ reads
\begin{equation}\label{velocity}
	\bm{c_0}=c(0,0)^T,\bm{c_1}=c(1,0)^T,\bm{c_2}=c(0,1)^T,\bm{c_3}=c(-1,0)^T,\bm{c_4}=c(0,-1)^T,
\end{equation}
the transform matrix $\bm{M}$ and the diagonal relaxation matrix $\bm{S}$ can be given by
\begin{equation}\label{ms}
	\begin{aligned}
		&\bm{M}=\left (
		\begin{matrix}
			1& 1 &1&1&1\\
			0 & c & 0&-c &0\\
			0 & 0&c & 0&-c\\
			-4c^2 & c^2 & c ^2&c^2&c^2\\
			0 & c^2 & -c^2&c^2 &-c^2\\
		\end{matrix}
		\right ),
		\bm{S}=\left (
		\begin{matrix}
			s_1& 0 &0&0&0\\
			0 & s_2 & 0 &0&0\\
			0 & 0 & s_2&0&0\\
			0 & 0 & 0&s_4&0\\
			0 & 0 & 0&0&s_5\\
		\end{matrix}
		\right ).
	\end{aligned}
\end{equation}
It should be noted that the choice of the transfrom matrix $\bm{M}$ is not unique (see another one in \ref{AM}), and the diagonal element $s_i$ of the relaxation matrix $\bm{S}$ is the relaxation parameter corresponding to the orthogonal moment of the distribution functions. To ensure the  diffusion coefficient $\kappa$ to be positive, $s_i$ should be located in the range $(0,2)$. 

In the LB method, to derive correct macroscopic diffusion
equation (\ref{distributionf}), the equilibrium distribution function can be simply defined as
\begin{equation}\label{feq}
	f_i^{eq}=w_i\phi,
\end{equation}
and the zeroth, first and second-order velocity moments of the equilibrium distribution function $f_i^{eq}$ satisfy the following conditions \cite{Chai2014},
\begin{subequations}\label{jieju}
	\begin{align}
		&\sum_if_i^{eq}=\phi,\label{0jieju}\\
		&\sum_if_i^{eq}\bm{c_i}=0,\label{1jieju}\\
		&\sum_if_i^{eq}\bm{c_ic_i}=\frac{1-w_0}{2}c^2\phi\bm{I},\label{2jieju}
	\end{align}
\end{subequations}
where $\bm{I}$ is the unity matrix. The discrete source term $R_i$ is defined by
\begin{equation}
	R_i=w_iR,(i=0,1,2,3,4),
\end{equation}
where 
$w_0\in (0, 1)$ and $w_{1-4}\in (0,1/4)$ are the weight coefficients.

From Eq. (\ref{jieju}), one can obtain
\begin{subequations}\label{w}
	\begin{align}
		&w_1=w_2=w_3=w_4=\frac{1-w_0}{4},\\
		&c^2_s=\frac{1-w_0}{2}c^2,
	\end{align}
\end{subequations}
where the weight coefficient $w_0$ can be considered as a free parameter. In addition, the macroscopic variable $\phi(\bm{r},t)$ (or the conservative moment) is calculated by
\begin{equation}\label{phi_fi}
	\phi(x,t)=\sum_i f_i(x,t)+\frac{\Delta t}{2}R.
\end{equation}
Through the direct Taylor expansion analysis \cite{Holdych2004, Wagner2006, Dubois2008}, one can correctly recover the diffusion equation (\ref{equde}) from the present MRT-LB model (\ref{distributionf}) with the following relation between the diffusion coefficient $\kappa$ and relaxation parameter $s_2$,
\begin{equation}\label{kappa}
	\kappa=c_s^2\big(\frac{1}{s_2}-\frac{1}{2}\big)\Delta t=\varepsilon\frac{\Delta ^2x}{\Delta t},\varepsilon :=\frac{1-w_0}{2}\big(\frac{1}{s_2}-\frac{1}{2}\big).
\end{equation}

Considering the relaxation parameter $s_5\in(0,2)$ and $\varepsilon$ as two free parameters, one can obtain a fourth-order MRT-LB model with the following relaxation matrix and weight coefficients,
\begin{equation}\label{SW}
	\begin{aligned}
		\bm{S}=\left (
		\begin{matrix}
			s_1& 0 &0&0&0\\
			0 & \frac{6s_5-12}{s_5-6} & 0 &0&0\\
			0 & 0 & \frac{6s_5-12}{s_5-6}&0&0\\
			0 & 0 & 0&\frac{12\varepsilon s_5^2-24\varepsilon s_5+2s_5^2}{2s_5-36\varepsilon+24\varepsilon s_5+\varepsilon s_5^2}&0\\
			0 & 0 & 0&0&s_5\\
		\end{matrix}
		\right ),
		&w_0=\frac{s_5+(6s_5-12)\varepsilon}{s_5},w_{1-4}=\frac{1-w_0}{4},
	\end{aligned}
\end{equation}
the details can be found in Sec. \ref{sec4}. 

	\section{The macroscopic six and five-level finite-difference schemes}\label{sec3}	
This section will show some details on how to derive the equivalent finite-difference scheme of the general MRT-LB model for the diffusion equation (\ref{equde}). Based on the fact that any LB models can be divided into two processes: a local collision step performed on each site of the lattice and a propagation step conducted among neighboring sites of the lattice, we multiply $\bm{M}$ on both sides of Eq. (\ref{distributionf}), and express the collision step in a vector form,
\begin{equation}\label{collision}
	\bm{m}^{\star}(\bm{r},t) = \big(\bm{I}-\bm{S}\big)\bm{m}(\bm{r},t)+\bm{S}\bm{m}^{eq}(\bm{r},t)+\Delta t(\bm{I}-\frac{\bm{S}}{2})\bm{M}\bm{R},\:\bm{r}\in(\mathcal{L}_x,\mathcal{L}_y),
\end{equation}
where the post-collision state is denoted by $\star$. The relation between the post-collision distribution function vector $\bm{f} ^{\star}(\bm{r},t)=\big(f^{\star}_{k}\big)_{k=0}^4$ and $\bm{m}^{\star}$ is given by $\bm{m}^{\star}=\bm{M}\bm{f}^{\star}$, and the source term vector $\bm{R} =(R_k)_{k=0}^4$.

According to the shift operator $T^{{\bm{c_k}}/c}_{\Delta x}$ $\big(T^{{\bm{c_k}}/c}_{\Delta x}[f_i(\bm{r},t)]:=f_i(\bm{r}-\bm{c_k}\Delta t,t)\big)$ defined in Ref. \cite{Bellotti2022}, the propagation step can be rewritten as
\begin{equation}\label{streaming}
	\bm{f} (\bm{r}+\bm{c}\Delta t,t+\Delta t) =\bm{ {\rm diag}}\big(T^{{\bm{c_0}}/c}_{\Delta x}, T^{\bm{c_1}/c}_{\Delta x}, T^{\bm{c_2}/c}_{\Delta x},T^{\bm{c_3}/c}_{\Delta x},T^{\bm{c_4}/c}_{\Delta x}\big)\bm{f} ^{\star}(\bm{r}+\bm{c}\Delta t,t),\:\bm{r}\in(\mathcal{L}_x,\mathcal{L}_y),
\end{equation}
where $\bm{f}(\bm{r},t)=\big(f_{k}\big)_{k=0}^4$, $\bm{c}=\big(\bm{c_k}\big)_{k=0}^4$. Multiplying $\bm{M}$ on both sides of Eq. (\ref{streaming}) and with the aid of Eq. (\ref{collision}), one can obtain
\begin{equation}\label{lb0}
	\bm{m}^{n+1}(\bm{r}) = \bm{P}\bm{m}^{n}(\bm{r})+\bm{Q}\bm{m}^{eq|n}(\bm{r})+\Delta t\tilde{\bm{R}}, \:\bm{r}\in(\mathcal{L}_x,\mathcal{L}_y),
\end{equation}
where $\bm{m}^{n}(\bm{r}):=\bm{m}(\bm{r},t_n)$, $\bm{T}:=\bm{M}\bm{{\rm diag}}\big(T^{{\bm{c_0}}/c}_{\Delta x}, T^{{\bm{c_1}}/c}_{\Delta x}, T^{{\bm{c_2}}/c}_{\Delta x},T^{{\bm{c_3}}/c}_{\Delta x},T^{{\bm{c_4}}/c}_{\Delta x}\big)\bm{M}^{-1}$, $\bm{P}:=\bm{T}\big(\bm{I}-\bm{S}\big)$ and $\bm{Q}:=\bm{TS}$. For the last term of Eq. (\ref{lb0}), based on the fact that $\bm{S},\bm{I},\bm{M}$ and $\bm{R}$ are the constant matrices (vector), we have $\tilde{\bm{R}}:=\bm{T}(\bm{I}-\bm{S}/2)\bm{MR}=(\bm{I}-\bm{S}/2)\bm{MR}$ because the constants are position-independent, and are not influenced by shift operator $T_{\Delta x}^{\bm{c_k}/c}$.

For the sake of brevity, we define $T_i:=T_{\Delta x}^{\bm{c_i}/c}$, and obtain the matrix $\bm{T}$ for the MRT-LB model (\ref{distributionf}) with the D2Q5 lattice structure (\ref{velocity}),
\begin{equation}\label{MT}
	\bm{T}=\left (
	\begin{matrix}
		\frac{\sum_{i=0}^4T_i}{5}&\frac{T_1-T_3}{2c}&\frac{T_2-T_4}{2c} &\frac{\sum_{i=1}^4T_i-4T_0}{20c^2}&\frac{T_1-T_2+T_3-T_4}{4c^2}\\
		c\frac{T_1-T_3}{5}&\frac{T_1+T_3}{2}&0&\frac{T_1-T_3}{20c}&\frac{T_1-T_3}{4c}\\
		c\frac{T_2 - T_4}{5}& 0&\frac{T_2+T_4}{2}&\frac{T_2-T_4}{20c}&\frac{T_4-T_2}{4c}\\
		c^2\frac{\sum_{i=1}^4T_i-4T_0}{5}&c\frac{T_1-T_3}{2}&c\frac{T_2 - T_4}{2}& \frac{\sum_{i=1}^4 T_i+16T_0}{20}&\frac{T_1-T_2+T_3-T_4}{4}\\
		c^2\frac{T_1 - T_2 + T_3 - T_4}{5}& c\frac{T_1 - T_3}{2}& c\frac{T_4-T_2}{2}&\frac{T_1-T_2+T_3-T_4}{20}&\frac{\sum_{i=1}^4T_i}{4} \\			
	\end{matrix}
	\right ).
\end{equation}
Based on the matrix $\bm{T}$ in Eq. (\ref{MT}), it is easy to obtain the corresponding matrices $\bm{P},\bm{Q}$ and vector $\tilde{\bm{R}}$. Then we have the characteristic polynomial $\mathcal{X}_{\bm{P}}$ of matrix $\bm{P}$ as
\begin{equation}
	\mathcal{X}_{\bm{P}}=\bm{X}^5+\upsilon_4 \bm{X}^4+\upsilon_3 \bm{X}^3+\upsilon_2\bm{X}^2+\upsilon_1\bm{X}+\upsilon_0\bm{I},
\end{equation}
where the characteristic polynomial coefficient $\upsilon_i\:(i=0,1,2,3,4)$ can be found in \ref{A1}. 

Following the Proposition 4 in Ref. \cite{Bellotti2022}, we have the following finite-difference scheme,
\begin{equation}\label{fd1}
	\phi^{n+1}_{i,j} = -\sum_{k=0}^{4} \upsilon_k \phi ^{n+k-4}_{i,j}+\big[\sum_{k=0}^{5}\big(\sum_{l=0}^{k} \upsilon_{5+l-k}\bm{P}^{l}\big)\big(\bm{Q}\bm{m}^{eq|n-k}(i\Delta x,j\Delta x)+\Delta t\tilde{\bm{R}}\big)\big]_1,
\end{equation}
where $\phi_{i,j}^{n}$ denotes $\phi(i\Delta x,j\Delta x, t_n),i,j\in\mathbbm{Z}$. $\bm{m}^{eq|n}(i\Delta x,j\Delta x)$ and $\tilde{\bm{R}}$ are given by
\begin{subequations}\label{meqr}
	\begin{align}
		\bm{m}^{eq}(i\Delta x, j\Delta x, t_n)&=\bm{M}\bm{f}^{eq}(i\Delta  x,j\Delta x, t_n)=[1,0,0,(1-5w_0)c^2,0]^{T}\phi^n_{i,j},\\
		\tilde{\bm{R}}&=[\frac{2-s_1}{2},0,0,\frac{1-5w_0}{2}(2-s_4)c^2,0]^TR.
	\end{align}
\end{subequations}
With the help of Eq. (\ref{meqr}) and substituting Eq. (\ref{xishu}) (see \ref{A1}) into Eq. (\ref{fd1}), we can obtain the following macroscopic six-level finite-difference (SLFD) scheme,
\begin{equation}\label{fd}
	\begin{aligned}
		\phi_{i,j}^{n+1}&=\alpha_1\phi_{i,j}^{n}+\alpha_2(\phi_{i-1,j}^{n}+\phi_{i+1,j}^{n}+\phi_{i,j-1}^{n}+\phi_{i,j+1}^{n})+\beta_1\phi_{i,j}^{n-1}+\beta_2(\phi_{i-1,j}^{n-1}+\phi_{i+1,j}^{n-1}+\phi_{i,j-1}^{n-1}+\phi_{i,j+1}^{n-1})\\
		&\quad+\beta_3(\phi_{i-1,j-1}^{n-1}+\phi_{i-1,j+1}^{n-1}+\phi_{i+1,j-1}^{n-1}+\phi_{i+1,j+1}^{n-1})
		+\gamma_1\phi_{i,j}^{n-2}+\gamma_2(\phi_{i-1,j}^{n-2}+\phi_{i+1,j}^{n-2}+\phi_{i,j-1}^{n-2}+\phi_{i,j+1}^{n-2})\\
		&\quad+\gamma_3(\phi_{i-1,j-1}^{n-2}+\phi_{i-1,j+1}^{n-2}+\phi_{i+1,j-1}^{n-2}+\phi_{i+1,j+1}^{n-2})+\zeta_1\phi_{i,j}^{n-3}
		+\zeta_2(\phi_{i-1,j}^{n-3}+\phi_{i+1,j}^{n-3}+\phi_{i,j-1}^{n-3}+\phi_{i,j+1}^{n-3})\\
		&\quad+\eta \phi_{i,j}^{n-4}
		+\Delta t\delta R,
	\end{aligned}
\end{equation}
where the parameters $\alpha_k, \zeta_k\:(k=1,2)$, $\beta_l,\gamma_l\:(l=1,2,3)$, $\eta$ and $\delta$ are given by
\begin{equation}\label{fdparameter}
	\begin{aligned}
		&\alpha_1=(w_0-1)s_4+1,\alpha_2=1-\frac{s_5}{4}-\frac{s_2}{2}-\frac{s_4w_0}{4},	\beta_1=(s_4w_0+s_5-2)(1-s_2),\\
		&\beta_2=\frac{(s_5+2s_2-4)(1-s_4)+w_0s_4(s_5+2s_2-3)}{4},\beta_3=\frac{(2-s_2-s_4w_0)(s_2+s_5-2)}{4},		\\
		&\gamma_1=s_4w_0(s_2-1)(s_5-1)+(1-s_4)(s_5-2)(s_2-1),\\
		&\gamma_2=\frac{(s_2-1)(s_4w_0(3-s_2-2s_5)+(2s_2+3s_5-s_2s_5-4))}{4},\\
		&\gamma_3=\frac{w_0s_4(s_2-1)(s_2+s_5-2)+(1-s_4)(s_2-2)(s_2+s_5-2)}{4},	\\
		&\zeta_1=(w_0s_4-1)(1-s_5)(s_2-1)^2,\zeta_2=\frac{s_4w_0(s_2-1)^2(s_5-1)+(1-s_4)(1-s_2)(2s_2+3s_5-s_2s_5-4)}{4},\\
		&\eta =(1-s_2)^2(1-s_4)(1-s_5),\:\delta =s_2^2s_4s_5.
	\end{aligned}
\end{equation}

We would like to point out that Eq. (\ref{fd}) is exactly equivalent to the general MRT-LB model (\ref{distributionf}) for the diffusion equation (\ref{equde}). In fact, simliar to the above discussion, we can also obtain the corresponding macroscopic six-level finite-difference scheme from a proper MRT-LB model for the diffusion equation with a linear source term (See \ref{A5} for details). Specially, if we set $s_4=1$ or $s_5=1$ in Eq. (\ref{fdparameter}), the SLFD scheme (\ref{fd}) will reduce to a five-level one, while here we only consider the case of $s_5=1$ (the reason will be provided in the next section).
\begin{equation}\label{fds5=1}
	\begin{aligned}
		\phi_{i,j}^{n+1}&=\alpha_1\phi_{i,j}^{n}+\alpha_2(\phi_{i-1,j}^{n}+\phi_{i+1,j}^{n}+\phi_{i,j-1}^{n}+\phi_{i,j+1}^{n})+\beta_1\phi_{i,j}^{n-1}+\beta_2(\phi_{i-1,j}^{n-1}+\phi_{i+1,j}^{n-1}+\phi_{i,j-1}^{n-1}+\phi_{i,j+1}^{n-1})\\
		&\quad+\beta_3(\phi_{i-1,j-1}^{n-1}+\phi_{i-1,j+1}^{n-1}+\phi_{i+1,j-1}^{n-1}+\phi_{i+1,j+1}^{n-1})+\gamma_1\phi_{i,j}^{n-2} +\gamma_2(\phi_{i-1,j}^{n-2}+\phi_{i+1,j}^{n-2}+\phi_{i,j-1}^{n-2}+\phi_{i,j+1}^{n-2})\\
		&		\quad+\gamma_3(\phi_{i-1,j-1}^{n-2}+\phi_{i-1,j+1}^{n-2}+\phi_{i+1,j-1}^{n-2}+\phi_{i+1,j+1}^{n-2})+\zeta_2(\phi_{i-1,j}^{n-3}+\phi_{i+1,j}^{n-3}+\phi_{i,j-1}^{n-3}+\phi_{i,j+1}^{n-3})+\Delta t\delta R,
	\end{aligned}
\end{equation}
where the parameters $\alpha_k\:(k=1,2),\beta_k,\gamma_k\:(k=1,2,3)$, $\zeta_2$ and $\delta$ are simplified by
\begin{equation}\label{fdcoe}
	\begin{aligned}
		&\alpha_1=(w_0-1)s_4+1,\alpha_2=\frac{3}{4}-\frac{s_2}{2}-\frac{s_4w_0}{4},\beta_1=(s_4w_0-1)(1-s_2),\\
		&\beta_2=\frac{(2s_2-3)(1-s_4)+w_0s_4(s_2-1)}{2},\beta_3=\frac{(2-s_2-s_4w_0)(s_2-1)}{4},\\
		&\gamma_1=(1-s_4)(1-s_2),\gamma_2=\frac{(1-s_4w_0)(1-s_2)^2}{4},\zeta_2=\frac{(s_4-1)(1-s_2)^2}{4},\\
		&\gamma_3=\frac{w_0s_4(s_2-1)^2+(1-s_4)(s_2-2)(s_2-1)}{4},\delta =s_2^2s_4.
	\end{aligned}
\end{equation}

	\section{The accuracy and stability analysis}\label{sec4}
Due to the equivalence between the MRT-LB model (\ref{distributionf}) and SLFD scheme (\ref{fd}), one can focus on one of them, while for simplicity, we only consider the equivalent SLFD scheme in this section. We first conduct a detailed theoretical analysis on the accuracy of the SLFD scheme (\ref{fd}). Then, with the von Neumann stability analysis method, we prove that both the MRT-LB model (\ref{distributionf}) and FLFD scheme (\ref{fds5=1}) are stable   if and only if $w_0\in(0,1),s_2,s_4\in(0,2)$, i.e., they are unconditionally stable.
\subsection{The accuracy of the macroscopic six-level finite-difference scheme}
We first apply the Taylor expansion to Eq. (\ref{fd}) at the position $\bm{r}=(i\Delta x,j\Delta x)$ and time $t_n$, and after some algebraic manipulations, one can obtain
\begin{equation}\label{tylor}
	\begin{aligned}
		&\big(1+\beta_1+4(\beta_2+\beta_3)+2\gamma_1+8(\gamma_2+\gamma_3)+3\zeta_1+12\zeta_2+4\eta\big)[\frac{\partial \phi}{\partial t}]_{i,j}^{n}\\
		&=\frac{2\alpha_2+2\beta_2+4\beta_3+2\gamma_2+4\gamma_3+2\zeta_2}{2}\frac{\Delta x^2}{\Delta t}\big\{[\frac{\partial ^ 2\phi}{\partial x^2}]_{i,j}^{n}+[\frac{\partial ^ 2\phi}{\partial y^2}]_{i,j}^{n}\big\}\\
		&\quad+\frac{\beta_1+4(\beta_2+\beta_3)+4\gamma_1+16(\gamma_2+\gamma_3)+9\zeta_1+36(\zeta_2+\zeta_3)+16\eta-1}{2}\Delta t [\frac{\partial ^ 2\phi}{\partial t^2}]_{i,j}^{n}\\
		&\quad+\frac{2\alpha_2+2\beta_2+4\beta_3+2\gamma_2+4\gamma_3+2\zeta_2}{24}\frac{\Delta  x^4}{\Delta t}\big\{[\frac{\partial ^4 \phi}{\partial x^4}]_{i,j}^{n}+[\frac{\partial ^4 \phi}{\partial y^4}]_{i,j}^{n}\big\}\\
		&\quad+(\beta_2+2\beta_3+2\gamma_2+4\gamma_3+3\zeta_2)\Delta x^2\{[\frac{\partial ^3 \phi}{\partial x^2\partial t}]_{i,j}^{n}+[\frac{\partial ^3 \phi}{\partial y^2\partial t}]_{i,j}^{n}\}\\
		&\quad+(\beta_3+\gamma_3)\frac{\Delta  x^4}{\Delta t}[\frac{\partial ^4 \phi}{\partial x^2\partial y^2}]_{i,j}^{n}+\delta R+\ldots \:.\\
	\end{aligned}
\end{equation}
Substituting Eq. (\ref{fdparameter}) into above equation and with the aid of the following relations derived from Eq. (\ref{tylor}),
\begin{equation}
	\begin{aligned}
		&[\frac{\partial ^2 \phi}{\partial t^2}]^{n}_{i,j}=\kappa^2\big\{[\frac{\partial ^4 \phi}{\partial x^4}]^{n}_{i,j}+[\frac{\partial ^4 \phi}{\partial y^4}]^{n}_{i,j}\big\}+2\kappa^2[\frac{\partial ^4 \phi}{\partial x^2\partial y^2}]^{n}_{i,j}+O(\Delta x^2+\Delta t), \\
		&[\frac{\partial ^3 \phi}{\partial t\partial x^2}]^{n}_{i,j}=\kappa\big\{[\frac{\partial ^4 \phi}{\partial x^4}]^{n}_{j}+[\frac{\partial ^4 \phi}{\partial x^2\partial y^2}]\big\}^{n}_{i,j}+O(\Delta x^2+\Delta t),\\
		&[\frac{\partial ^3 \phi}{\partial t\partial y^2}]^{n}_{i,j}=\kappa\big\{[\frac{\partial ^4 \phi}{\partial y^4}]^{n}_{j}+[\frac{\partial ^4 \phi}{\partial x^2\partial y^2}]\big\}^{n}_{i,j}+O(\Delta x^2+\Delta t),\\
	\end{aligned}
\end{equation}
we have
\begin{equation}
	\begin{aligned}
		[\frac{\partial \phi}{\partial t}]_{i,j}^{n}=&\kappa\big\{[\frac{\partial ^2\phi}{\partial x^2}]+[\frac{\partial ^2\phi}{\partial y^2}]\big\}_{i,j}^{n}+\frac{TR_{41}}{96s_2^3s_4s_5}\frac{\Delta x^4}{\Delta t}\big\{[\frac{\partial ^4\phi}{\partial x^4}]+[\frac{\partial ^4\phi}{\partial y^4}]\big\}_{i,j}^{n}+\frac{TR_{42}}{16s_2^3s_4s_5}\frac{\Delta x^4}{\Delta t}[\frac{\partial ^4\phi}{\partial x^2\partial y^2}]_{i,j}^n\\
		&+R+O(\Delta t^2+\Delta  x^4+\Delta x^2\Delta t).
	\end{aligned}
\end{equation}
At the diffusive scaling ($\Delta t \propto \Delta x^2$), one can obtain an explicit SLFD scheme with the second-order accuracy in time and fourth-order accuracy in space once the following conditions are satisfied,
\begin{subequations}\label{error}
	\begin{align}
		&	\varepsilon =\frac{1-w_0}{2}\big(\frac{1}{s_2}-\frac{1}{2}\big),\label{errorfirst}\\
		&	TR_{41}=(s_2 - 2)(w_0 - 1)\big(6s_2(s_4+s_5w_0)(2-s_2)  +s_2s_4s_5(s_2-6)(3w_0-1) +12s_4s_5(w_0-1)\big)
		=0,\label{errorsecond} \\
		&	TR_{42}=(2-s_2)(w_0 - 1)\big[2s_2(s_4+s_5w_0)(2-s_2) + s_4s_5(1-w_0)\big(4+s_2(s_2-6)\big) \big]=0,\label{errorthird}
	\end{align}
\end{subequations}
from which we can obtain the explicit expressions of $w_0$, $s_2$ and $s_4$ in terms of $s_5$ and $\varepsilon$
\begin{align}\label{solution}
	\left\{\begin{aligned}
		&w_0=\frac{s_5+(6s_5-12)\varepsilon}{s_5},\\
		&s_2=\frac{6s_5-12}{s_5-6},\\
		&s_4=\frac{12\varepsilon s_5^2-24\varepsilon s_5+2s_5^2}{2s_5-36\varepsilon+24\varepsilon s_5+\varepsilon s_5^2}.
	\end{aligned}\right.
\end{align}
This indicates that when the weight coefficient $w_0$ and the relaxation parameters $s_2$, $s_4$ and $s_5$ satisfy above Eq. (\ref{solution}) or Eq. (\ref{SW}), a fourth-order MRT-LB model can be derived.

In the following, we give some remarks on above results.

\begin{remark}
We note that the relaxation parameter $s_1$ corresponding to the zeroth-order moment of distribution functions $f_i(\bm{r},t)$ [the conservative variable $\phi(\bm{r},t)$] does not appear in the SLFD scheme (\ref{fd}), this means that it has no influence on the numerical results. This also explains why the relaxation parameter $s_1$ in the MRT-LB model (\ref{distributionf}) can be chosen arbitrarily, which is consistent with the previous works \cite{Lin2022, Chen2023}. However, unlike the relaxation parameter $s_1$, the relaxation parameter $s_2$ corresponding  to the first-order moment of the distribution function is related to the diffusion coefficient $\kappa$ [see Eq. (\ref{kappa}) or (\ref{errorfirst})], and the relaxation parameters $s_4$ and $s_5$ corresponding to the second-order moments of the distribution function have some important influences on the SLFD scheme (\ref{fd}) [see Eqs. (\ref{errorsecond}) and (\ref{errorthird})]. These similar results have also been reported on the MRT-LB model for one-dimensional diffusion equations \cite{Lin2022} and convection-diffusion equations \cite{Chen2023}.
\end{remark}
\begin{remark}
If we consider the case of $s_1=s_4=s_5$, i.e., the TRT-LB model, one can show that the condition $\Lambda^{eo} = 1/6$ and $\Lambda^{BGK} = 1/12$ used in the TRT-LB model has a fourth-order convergence rate \cite{Gin2012}, and it is only a special case of Eq. (\ref{SW}) or Eq. (\ref{solution}).
\end{remark}

\begin{remark} Now let us focus on the Eq. (\ref{solution}). It is obvious that the value of parameter $\varepsilon$ must be specified properly to make other parameters reasonable, i.e, the weight coefficients $w_0\in(0,1),w_{1-4}\in(0,1/4)$, the relaxation parameters $s_2$, $s_4$ and $s_5$ should be located in the range $(0,2)$. However, on the other hand, one can find that Eq. (\ref{fd}) is a six-level finite-difference scheme, and this will bring some liminations to the implementation of numerical experiments and applications to some extent. To partially solve this problem, we consider a simple five-level difference (FLFD) scheme under the condition of $s_5 = 1$ [see Eqs. (\ref{fds5=1}) and (\ref{fdcoe})] without reducing the accuracy order of the finite-difference scheme (\ref{fd}). In addition, one can also adjust the value of the parameter $\varepsilon$, lattice spacing $\Delta x$ and time step $\Delta t$ to give a specified diffusion coefficient $\kappa=\varepsilon \Delta x^2/\Delta t$ in the implementation of numerical simulations. Actually, we do not need to worry about the numerical instability problem, this is because the finite-difference scheme with $s_5=1$ is unconditionally stable, as shown in the part \ref{stability-analysis}. For this reason, the relaxation parameter $s_5=1$ is used in the following theoretical analysis and numerical
simulations, and then the relations between the weight coefficients $w_0$ and $w_1$, the relaxation parameters $s_2$ and $s_4$ and the discretiation parameter $\varepsilon$ can be obtained from Eq. (\ref{solution}) and are shown in Fig. \ref{parametersw}.
\begin{figure}
	\centering
	\includegraphics[width=10cm,height=7cm]{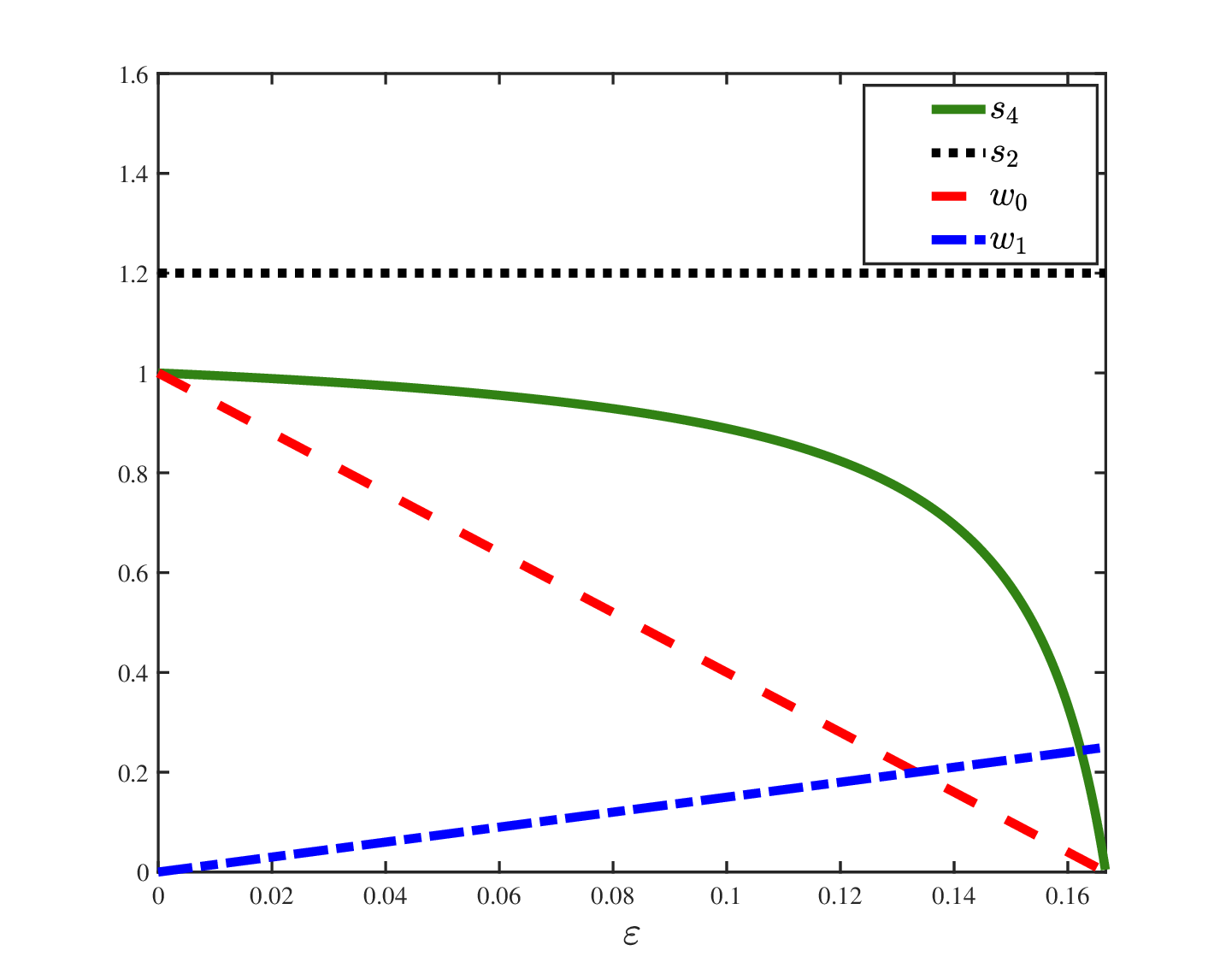}
	\caption{The weight coefficients $w_0$ and $w_1$, the relaxation parameters
		$s_2$ and $s_4$ as a function of the discretization parameter $\varepsilon$.}\label{parametersw}
\end{figure}
\end{remark}
\begin{remark}
To preserve the fourth-order accuracy of MRT-LB model (\ref{distributionf}) and FLFD scheme (\ref{fds5=1}), the initial condition of distribution function $f_i(\bm{r},t)$ in the MRT-LB model (\ref{distributionf}) and the initial value of $\phi(\bm{r},t)$ in FLFD scheme (\ref{fds5=1}) must be given properly. According to the previous work \cite{Chai2020}, we adopt the following expression to initialize the distribution function $f_i$ in implementing the MRT-LB model (\ref{distributionf}),
\begin{equation}\label{initial}
	\bm{f}=\bm{f}^{eq}-\Delta t\bm{\Lambda}^{-1}\mathcal{D}\bm{f}^{eq}+\Delta t
	\bm{\Lambda}^{-1}\big(\bm{I}-\frac{\bm{\Lambda}}{2}\big)\bm{R},
\end{equation}
where $\bm{\Lambda}=\bm{M}^{-1}\bm{S}\bm{M}$, $\mathcal D=\textbf{diag}\big(D_0,D_1,D_{2},D_3,D_4\big)$, $D_i=\partial _t+\bm{c_i}\cdot \nabla $. Additionally, it should be noted that with the help of Eq. (\ref{equde}), the term  $\partial_t \phi$ appeared in Eq. (\ref{initial}) can be replaced by the spatial derivatives of scalar variable $\phi$, which can be further determined from the initial condition. However, in the implementation of the FLFD scheme (\ref{fds5=1}), the values of the variable $\phi$ at the first four time levels are needed. Thus, in addition to the initial condition,  one must adopt some other numerical schemes to obtain the values of variable $\phi$ at the second, third and fourth time levels. The inclusion of the extra time levels would also bring a larger memory to store the variable $\phi$.
\end{remark}
\begin{remark}
The treatment on boundary condition is another important issue, and also affects the accuracies of the MRT-LB model (\ref{distributionf}) and FLFD scheme (\ref{fds5=1}). For the Dirichlet boundary condition in the MRT-LB model with the D2Q5 lattice structure, the problem is that there are two unknown distribution functions $f_i$ at the four concer points, and to overcome the problem, we adopt the following equation to obtain the unknown distribution functions $f_i$ \cite{Chai2020},
\begin{equation}\label{boundary}
	\bm{f}=\bm{f}^{eq}-\Delta t\bm{\Lambda}^{-1}\mathcal{D}\bm{f}^{eq}+\Delta t
	\bm{\Lambda}^{-1}\big(\bm{I}-\frac{\bm{\Lambda}}{2}\big)\bm{R}	+\Delta t^2\bm{\Lambda}^{-1}\mathcal{D}\big[\big(\bm{I}-\frac{\bm{\Lambda}}{2}\big)\bm{\Lambda}^{-1}\big]\mathcal{D}\bm{f}^{eq}.
\end{equation}
where the spatial derivatives of $\bm{f}^{eq}$ in Eq. (\ref{boundary}) can be obtained from the Dirchlet boundary conditions (see Sec. \ref{secNumer} for details).
\end{remark}

\subsection{The stability analysis}\label{stability-analysis}
We now discuss the stability of the MRT-LB model (\ref{distributionf}) and FLFD scheme (\ref{fds5=1}). With the help of the von Neumann stability analysis and the Routh-Hurwitz stability criterion \cite{Routh, Hurwitz, Gantmacher}, we can show the following Theorem.

\textbf{Theorem.}
The MRT-LB model (\ref{distributionf}) and FLFD scheme (\ref{fds5=1}) are stable if and only if $0<w_0<1,0<s_2,s_4<2$, i.e., the MRT-LB model (\ref{distributionf}) and FLFD scheme (\ref{fds5=1}) are unconditionally stable.

\textbf{Proof.} Due to the equivalence between the MRT-LB model (\ref{distributionf}) and FLFD scheme (\ref{fds5=1}) \cite{Bellotti2022}, we only consider the stability condition of latter. 

If we take the discrete Fourier transform to Eq. (\ref{fds5=1}), one can obtain the amplification matrix $\bm{G}(s_2,s_4,w_0,\theta_1,\theta_2)$ as 
\begin{equation}
	\bm{G}=\left (
	\begin{matrix}
		\alpha_1+2\alpha_2(\chi_1+\chi_2)&1&0&0\\
		\beta_1+2\beta_2(\chi_1+\chi_2)+4\beta_3\chi_1\chi_2&0&1&0\\
		2\gamma_2(\chi_1+\chi_2)+4\gamma_3\chi_1\chi_2&0&0&1\\
		2\zeta_2(\chi_1+\chi_2)&0&0&0 \\			
	\end{matrix}
	\right ),
\end{equation}
where $\chi_{1,2}=\cos\theta_{1,2}$ with $\theta_{1,2}\in[-\pi,\pi]$. Then the
characteristic polynomial of $\bm{G}$ can be expressed as
\begin{equation}\label{charac}
	p(\lambda)=\lambda^4+a_1\lambda^3+a_2\lambda^2+a_3\lambda+a_4,
\end{equation}
the coefficients $a_i\:(i=1,2,3,4)$ are given by
\begin{equation}
	\begin{aligned}
		&a_1=-[\alpha_1+2\alpha_2(\chi_1+\chi_2)],\\
		&a_2=-[\beta_1+2\beta_2(\chi_1+\chi_2)+4\beta_3\chi_1\chi_2],\\
		&a_3=-[2\gamma_2(\chi_1+\chi_2)+4\gamma_3\chi_1\chi_2],\\
		&a_4=-2\zeta_2(\chi_1+\chi_2).\\			
	\end{aligned}
\end{equation}

From the von Neumann stability analyisis, the FLFD scheme (\ref{fds5=1}) is stable if and only if the roots of the characteristic polynomial $p(\lambda)$ denoted by $\lambda_k$ ($k=1,2,3,4$) are inside the unit disk and those on the unit disk are simple, i.e., the characteristic polynomial $p(\lambda)$ is a von Neumann polynamial. In the following, we will divide the proof process into four steps.\\
\textbf{Step 1. }
We first define two polynominals as
\begin{subequations}
	\begin{align}
		&p^*(\lambda)=a_4\lambda^4+a_3\lambda^3+a_2\lambda^2+a_1\lambda+1,\\
		&q(\lambda)=\frac{p^*(0)p(\lambda)-p(0)p^*(\lambda)}{\lambda}=p_1\lambda^3+p_2\lambda^2+p_3\lambda+p_4,\label{qcharc}
	\end{align}
\end{subequations}
where $p_i
\:(i=1,2,3,4)$ are given by
\begin{subequations}\label{pp}
	\begin{align}
		&p_1=1-a_4^2>0,\label{p1leq0}\\
		&p_2=a_1-a_3a_4,\\
		&p_3=a_2-a_2a_4,\\
		&p_4=a_3-a_1a_4.
	\end{align}
\end{subequations}
It is obvious that $|p^*(0)|>|p(0)|$ under the conditions of $0<w_0<1,0<s_2,s_4<2$ and $|\chi_{1,2}|\leq1$. Following the Theorem 6.1 in Ref. \cite{Miller}, the polynomial $p(\lambda)$ (\ref{charac}) is a von Neumann polynomial if and only if $q(\lambda)$ (\ref{qcharc}) is a von Neumann polynomial. 

Now we focus on the polynomial $q(\lambda)$ (\ref{qcharc}). With the linear fractional transformation
\begin{equation}\label{tran}
	\lambda=\frac{1+z}{1-z},z\in\mathbbm{Z},
\end{equation}
the unit circle $|\lambda| = 1$ and the field $|\lambda| < 1$ are mapped to the imaginary axis [Re($z$)=0] and left-half plane [Re($z$) $<$ 0], and vise verse. Here $\mathbbm{Z}$ and Re denote the complex-number field and the real part of a complex number. Substituting
Eq. (\ref{tran}) into Eq. (\ref{qcharc}), we can obtain
\begin{equation}
	\begin{aligned}
		f(z)&=(1-z)^3\big(\frac{1+z}{1-z}\big)^3q\big(\frac{1+z}{1-z}\big)\\
		&=(p_1-p_2+p_3-p_4)z^3+(3p_1+3p_4-p_3-p_2)z^2
		+(3p_1-3p_4-p_3+p_2)z+(p_1+p_2+p_3+p_4).\\		
	\end{aligned}
\end{equation}
To ensure $q(\lambda)$ (\ref{qcharc}) to be the von Neumann polynomial, the following conditions must be satisfied,	
\begin{subequations}\label{2neq}
	\begin{align}
		p_1-p_2+p_3-p_4>0,\label{neq21}\\
		p_1+p_2+p_3+p_4>0.\label{neq22}\\	
		p_1-p_4>0,\label{neq23}\\
		p_1+p_4>0,\label{neq24}\\	
		p_1^2 - p_1p_3 - p_4^2 + p_2p_4>0,\label{neq25}
	\end{align}
\end{subequations}
which is equivalent with the Routh-Hurwitz stability criterion \cite{Routh, Hurwitz, Gantmacher}, the details can be found in \ref{stability}.

To simplify following analysis, we now introduce a domain of the parameters as
\begin{equation}
	\Omega=\{(s_2,s_4,w_0,\chi_1,\chi_2)\in(0,2)\times(0,2)\times(0,1)\times(-1,1)\times(-1,1)\}.
\end{equation}
\textbf{Step 2. The proof of Eq. (\ref{neq21})}

From Eq. (\ref{pp}), we have
\begin{equation}\label{step21}
	\begin{aligned}
		&p_1-p_2+p_3-p_4=\underbrace{-\frac{1}{4}[(s_4-1)(s_2-1)^2(\chi_1+\chi_2)+2]}_{<0}\times F(s_2,s_4,w_0,\chi_1,\chi_2)\:\:{\rm in}\:\: \Omega,\\
	\end{aligned}
\end{equation}
where
\begin{equation}\label{F}
	\begin{aligned}
		F=&(\chi_1+\chi_2)\big[4(s_2-1)(2-s_4)+s_2^2\big(s_4(w_0+1)-2\big)\big]+\chi_1\chi_2[4(1-s_2)(s_2+s_4-2)+2s_2s_4(s_2-1)(1-w_0)]\\
		&+2s_2s_4(w_0-1)+4(s_2+s_4-2).
	\end{aligned}
\end{equation} 
Here the aim is to prove $F<0$ in the open region $\Omega$, it means that the maximum value of $F$ is strictly less than zero in $\Omega$. Considering that the function $F$ may have no maximum value in the open region $\Omega$, we consider the function $F$ in the closed region $\overline{\Omega}$, and find that the function $F$ is continuously differentiable over the region $\Omega$. Thus, the maximum value of the function $F$ in the closed region $\overline{\Omega}$ can only be obtained at the stationary  points in $\Omega$ and boundary points at $\partial \overline{\Omega}$. Based on above discussion, the task is to prove $F<0$ at the stationary points in the open region $\Omega$ and $F\leq0$ at $\partial \overline{\Omega}$. 

First, we introduce the set of the stationary points and the corresponding values of the function $F$ as $A_{F}$ and $G_{F}$. After some algebraic manipulations, we get $A_F=\emptyset$ and $G_F=\emptyset$, thus the maximum value of the function $F$ in the closed region $\overline{\Omega}$ is obtained at $\partial \overline{\Omega}$. It should be noted that the function $F$ is symmetric with respect to $\chi_1$ and $\chi_2$, thus we only need to consider the following eight cases of the function $F$ (\ref{F}) at the boundary points,
\begin{subequations}\label{F2}
	\begin{align}
		&	F(s_2=0)= \underbrace{4(s_4-2)(1+\chi_1)(1+\chi_2)}_{<0}\:\:{\rm in}\:\:\Omega_{s_2=0},\label{Fs20} \\ 
		&	F(s_2=2)=\underbrace{-4s_4w_0(1-\chi_1)(1-\chi_2)}_{<0}\:\:{\rm in}\:\:\Omega_{s_2=2},\label{Fs22}\\
		&	F(s_4=0)=\underbrace{2(s_2 - 2)[s_2+(2-s_2)(1+\chi_1)(1+\chi_2)]}_{<0}\:\:{\rm in}\:\:\Omega_{s_4=0},\label{Fs40}\\
		&	F(s_4=2)=\underbrace{2s_2w_0[s_2(\chi_1-1)(1-\chi_2) +(s_2-2)(1-\chi_1\chi_2)]}_{<0}\:\:{\rm in}\:\:\Omega_{s_4=2},\label{Fs42}\\
		&	F(w_0=0)=\underbrace{-(2-s_2)(2-s_4)[(1+\chi_1)(1+\chi_2)(2-s_2) +s_2(1-\chi_1\chi_2) ]}_{<0}\:\:{\rm in}\:\:\Omega_{w_0=0},\label{Fw00}\\
		&F(w_0=1)=4(1+\chi_1\chi_2+\chi_1+\chi_2)\big[s_2(2-s_2)-\big((1-s_2)(1-s_4)+1\big)-\frac{1}{4}\big]\notag\\
		&\quad\quad\quad\quad\quad+(\chi_1+\chi_2)[1-2s_2(2-s_2)+2s_2^2s_4] +\chi_1\chi_2+[1-4s_2(2-s_2)]\:\:{\rm in}\:\:\Omega_{w_0=1},\label{Fw01}\\
		&	F(\chi_1=-1)=-s_2[(1+\chi_2)(4-2s_2-2s_4+s_2s_4+2s_4w_0)+ s_2s_4w_0(1-3\chi_2) ]\:\:{\rm in}\:\:\Omega_{\chi_1=-1},\label{Fchim1}\\
		&	F(\chi_1=1)=(2-s_2 )[(1+\chi_2)(-8+6s_2+4s_4-3s_2s_4+s_2s_4w_0)-4s_2 +2 s_2s_4(1-w_0)\:\:{\rm in}\:\:\Omega_{\chi_1=1},\label{Fchip1}
	\end{align}
\end{subequations}
where the open region $\Omega_{s_2=0}$ is defined as $\{(s_4,w_0,\chi_1,\chi_2)\in (0,2)\times(0,1)\times(-1,1)\times(-1,1)\}$, and the other cases in Eq. (\ref{F2}) can be defined similarly. 

Actually, it is difficult to determine directly whether these functions, i.e., Eqs. (\ref{Fw01}), (\ref{Fchim1}) and (\ref{Fchip1}), are no more than zero in the corresponding regions. Similar to the way we analyze the function $F$ (\ref{F}), we first consider the stationary points of functions $F(w_0=1),F(\chi_1=-1)$ and $F(\chi_1=1)$ in their open regions $\Omega_{w_0=1},\Omega_{\chi_1=-1}$ and $\Omega_{\chi_1=1}$, and have the following results,
\begin{subequations}\label{AF}
	\begin{align}
		&A_{F(w_0=1)}=\{(s_2=1,s_4=1,z_{\chi_1},-z_{\chi_1})\},G_{F(w_0=1)}={-4},\label{AFw01}\\ 
		&A_{F(\chi_1=-1)}=\emptyset,G_{F(\chi_1=-1)}=\emptyset,\label{AFchim1}\\
		&A_{F(\chi_1=1)}=\emptyset,G_{F(\chi_1=1)}=\emptyset,\label{AFchip1}
	\end{align}
\end{subequations}
where $z_{\upsilon}$ in Eq. (\ref{AFw01}) represents any value in the range of parameter $\upsilon$. For the function $F(w_0=1)$, we know that the maximum value may be -4 at the point $(s_2=1,s_4=1,w_0=1,z_{\chi_1},-z_{\chi_1})$ or larger at the boundary points $\partial\overline{\Omega}_{w_0=1}$, while for the functions $F(\chi_1=-1)$ and $F(\chi_1=1)$, the maximum value must be obtained at the boundary points $\partial\Omega_{\chi_1=-1}$ and $\partial\Omega_{\chi_1=1}$. We continue to consider the expressions of $F(w_0=1),F(\chi_1=-1)$ and $F(\chi_1=1)$ at the boundary points [Note: from $F\big(\upsilon=a\big)<0$ in $\Omega_{\upsilon=a}$, we have $F\big(w_0=1,\upsilon=a\big)\leq0$ in $\Omega_{w_0=1,\upsilon=a}$ \big(see Eqs.  (\ref{Fs20}-\ref{Fw00})\big), thus in the following, we do not need to list the case of $F(w_0=1,\upsilon=a)$],
\begin{subequations}
	\begin{align}
		&F(w_0=1,\chi_1=-1)= \underbrace{2s_2(1+\chi_2)(s_2-2) +2s_2^2s_4(\chi_2-1)} _{<0}\:\:{\rm in}\:\:\Omega_{w_0=1,\chi_1=-1},\label{Fw01chim1}\\ 
		&F(w_0=1,\chi_1=1)=\underbrace{2(2-s_2)[(1+\chi_2)(2-s_2)(s_4-2)+s_2(\chi_2-1)]}_{<0}\:\:{\rm in}\:\:\Omega_{w_0=1,\chi_1=1},\label{Fw01chip1}\\
		&F(\chi_1=-1,\chi_2=-1)=\underbrace{-4s_2^2s_4w_0}_{<0}\:\:{\rm in}\:\:\Omega_{\chi_1=-1,\chi_2=-1},\label{Fchim1chi2m1}\\
		&F(\chi_1=-1,\chi_2=1)=F(\chi_1=1,\chi_2=-1)=\underbrace{2s_2(s_2 - 2)[2+s_4(w_0-1)]}_{<0}\:\:{\rm in}\:\:\Omega_{\chi_1/\chi_2=-1,\chi_2/\chi_1=1},\label{Fchim1chi2p1}\\
		&F(\chi_1=1,\chi_2=1)=\underbrace{4(s_2 - 2)^2(s_4 - 2)}_{<0}\:\:{\rm in}\:\:\Omega_{\chi_1=1,\chi_2=1}.\label{Fchip1chi2p1}
	\end{align}
\end{subequations}
Then we have
\begin{itemize}
	\item [(i)]Under the condition of Eqs. (\ref{Fw01}), (\ref{AFw01}), (\ref{Fs20}-\ref{Fs42}) and (\ref{Fw01chim1}-\ref{Fw01chip1}), we have $F(w_0=1)\leq0$ in $\Omega_{w_0=1}$.
	\item [(ii)] Under the condition of Eqs. (\ref{Fchim1}), (\ref{AFchim1}), (\ref{Fs20}-\ref{Fw00}), (\ref{Fw01chim1}) and (\ref{Fchim1chi2m1}-\ref{Fchim1chi2p1}), we have $F(\chi_1=-1)\leq0$ in $\Omega_{\chi_1=-1}$.
	\item [(iii)] Under the condition of Eqs. (\ref{Fchip1}), (\ref{AFchip1}), (\ref{Fs20}-\ref{Fw00}), (\ref{Fw01chip1}) and (\ref{Fchim1chi2p1}-\ref{Fchip1chi2p1}), we have $F(\chi_1=1)\leq0$ in $\Omega_{\chi_1=1}$. 
\end{itemize}
From above discussion, we find that $F(s_2,s_4,w_0,\chi_1,\chi_2)$ (\ref{F}) has no stationary points in the open region $\Omega$, and is no more than zero at $\partial \overline{\Omega}$, this means that the condition (\ref{neq21}) holds.\\
\textbf{Step 3. The proof of Eq. (\ref{neq22})}

From Eq. (\ref{pp}), we can obtain
\begin{equation}\label{neq222}
	\begin{aligned}
		p_1+p_2+p_3+p_4=&\frac{s_4}{2}(2-s_2)(1-w_0)[(2-s_2)(1-\chi_1)(1-\chi_2)+s_2(1-\chi_1\chi_2)]\times\\
		&[2+(\chi_1+\chi_2)(s_4-1)(s_2-1)^2],
	\end{aligned}
\end{equation}
it is obvious that $p_1+p_2+p_3+p_4>0$ in the open region $\Omega$.\\
\textbf{Step 3. The proof of Eqs. (\ref{neq23}) and (\ref{neq24})}

To prove Eqs. (\ref{neq23}) and (\ref{neq24}), we only need to show $p_1^2-p_4^2>0$ [see $p_1>0$ in Eq. (\ref{p1leq0})]. Actually, the expression of $p_1^2-p_4^2$ can be given by
\begin{equation}
	\begin{aligned}
		p_1^2-p_4^2=&\Big[\frac{\big(\chi_1+\chi_2)^2(s_2 - 1)^4(s_4 - 1)^2}{4} - 1\Big]^2 -\big[(s_2 - 1)(s_4 - 1) - \frac{(\chi_1+ \chi_2)(s_2 - 1)(s_4w_0(s_2 - 1) - s_2 + 1)}{2} \\
		&+ \chi_1\chi_2\big(s_4w_0(s_2 - 1)^2 - (s_2 - 1)(s_2 - 2)(s_4 - 1)\big)+\frac{ (\chi_1+\chi_2)(s_2 - 1)^2(s_4 - 1)}{4}\times\\
		&\big(2s_4(w_0 - 1) - (\chi_1+\chi_2)(2s_2 + s_4w_0 - 3) + 2\big)\big]^2\:\:{\rm in}\:\:\Omega.
	\end{aligned}
\end{equation}
Due to the complexity of the expression of $p_1^2-p_4^2$, we define the function $H(s_2,s_4,w_0,\chi_1,\chi_2)=p_1^2-p_4^2$, and find that it is symmetric with respect to
the parameters $\chi_{1}$ and $\chi_2$. Simliar to above discussion on the function $F$ (\ref{F}), in order to show the function $H>0$ in the open region $\Omega$, we only need to consider the function $H$ in the closed region $\overline{\Omega}$, and analyze whether the function in the stationary point set $A_{H}$ is strictly larger than zero and that at $\partial\overline{\Omega}$ is no less than zero.
To this end, we consider the sets $A_{H}$ and $G_H$ in the open region $\Omega$, and obtain
\begin{equation}\label{AHG}
	A_H=\{(z_{s_2},s_4=1,z_{w_0},\chi_1=0,\chi_2=0),(s_2=1,z_{s_4},z_{w_0},z_{\chi_1},z_{\chi_2})\},G_H=\{1\}.
\end{equation}
It is obvious that the value of function $H$ is strictly larger than zero at stationary points. For the values of function $H$ at $\partial \overline{\Omega}$, we will obtain some functions with at most four, three, two, one and zero variable(s) when fixing one, two, three, four and five of the parameters $s_2,s_4,w_0,\chi_1$ and $\chi_2$.  Thus,  we only need to analyze the stationary points of function $H$ with one, two, three, four and five of these parameters $s_2,s_4,w_0,\chi_1$ and $\chi_2$ fixed, and prove that the corresponding function values at the stationary points are no less than zero.

Fixing one of the variables $s_2,s_4,w_0$ and $\chi_1$ of the function $H$, we have the following sets (due to the symmetry of function $H$ with respect to $\chi_1$ and $\chi_2$ , we only need to consider one case),
\begin{equation}
	\begin{aligned}
		&A_{H(s_2=0)}=\emptyset,G_{s_2=0}=\emptyset,A_{H(s_2=2)}=\emptyset,G_{H(s_2=2)}=\emptyset,A_{H(s_4=0)}=\{(s_2=1,z_{w_0},z_{\chi_1},z_{\chi_2}),\},
		G_{H(s_4=0)}=\{1\},\\
		&A_{H(s_4=2)}=\{(s_2=1,z_{w_0},z_{\chi_1},z_{\chi_2})\},G_{H(s_4=2)}=\{1\},\\	
		&A_{H(w_0=0)}=\{(s_2=1,z_{s_4},z_{\chi_1},z_{\chi_2}),(z_{s_2},s_4=1,z_{\chi_1},-z_{\chi_1})\},G_{H(w_0=0)}=\{1\},\\
		&A_{H(w_0=1)}=\{(s_2=1,z_{s_4},z_{\chi_1},z_{\chi_2}),(z_{s_2},s_4=1,\chi_1=0,z_{\chi_2}),(z_{s_2},s_4=1,z_{\chi_1},\chi_2=0)\},
		G_{H(w_0=1)}=\{1\},\\
		&A_{\chi_1=-1}=\{(s_2=1,z_{s_4},z_{w_0},z_{\chi_2})\},G_{H(\chi_1=-1)}=\{1\},A_{H(\chi_1=1)}=\{(s_2=1,z_{s_4},z_{w_0},z_{\chi_2})\},G_{H(\chi_1=1)}=\{1\}.
	\end{aligned}
\end{equation}
From above results, we can find that the values of the function $H$ with one variable fixed at the stationary points in the corresponding open region are larger than zero.

In addition, one can see from Table \ref{TAfix2} that the values of the function $H$ at stationary points are all larger zero  when fixing two
of these variables $s_2,s_4,w_0,\chi_1$ and $\chi_2$. Similarly, as shown in Table \ref{TAfix3}, the function values at stationary points are no less than zero when three variables are fixed. If we fix four variables of $s_2, s_4, w_0,\chi_1$ and $\chi_2$, the function $H$ becomes a linear function, and the function values at stationary points are still no less than zero, as seen from Table. \ref{TAfix4}. When fixing five parameters of $s_2, s_4,w_0, \chi_1$ and $\chi_2$, we can obtain 32 values of the function $H$, and it is easy to show that the corresponding values are all equal to zero. Therefore, based on the results in Eq. (\ref{AHG}), Tables \ref{TAfix2}, \ref{TAfix3} and \ref{TAfix4}, one can conclude that the values of function $H$ are strictly larger than zero in the open region $\Omega$ and no less than zero at $\partial\overline{\Omega}$, i.e., Eqs. (\ref{neq23}) and (\ref{neq24}).

	\begin{table}
	\caption{The stationary points and corresponding values of function $H$ when fixing two variables of $s_2$, $s_4$, $w_0$, $\chi_1$ and $\chi_2$.}
	\centering
	\label{TAfix2}  
	\begin{threeparttable}   
		\begin{tabular}{ccc}
			\hline \hline
			Fixed variable(s)&Stationary ponit(s)&Function value(s)\\
			\hline
			$s_2=0(2),s_4=0\tnote{1}$&-\tnote{2}&-\tnote{2}\\
			$s_2=0(2),s_4=2$&-&-\\
			$w_0=0,s_2=0(2)$&$(1,z_{\chi_1},z_{\chi_2}),z_{\chi_1}=-z_{\chi_2}$&1\\
			$w_0=1,s_2=0(2)$&$(s_4=1,\chi_1=0,z_{\chi_2}),  (s_4=1,z_{\chi_1},\chi_2=0)$&1\\
			$\chi_1=-1,s_2=0(2)$&-&-\\
			$w_0=0,s_4=0(2)$&$(s_2=1,z_{\chi_1},z_{\chi_2})$&1\\
			$\chi_1=-1,s_4=0(2)$&$(s_2=1,z_{w_0},z_{\chi_2})$&1\\
			$w_0=0(1),\chi_1=-1$&$(s_2=1,z_{s_4},z_{\chi_2})$&1\\
			$w_0=0(1),\chi_1=1$&$(s_2=1,z_{s_4},z_{\chi_2})$&1\\
			$\chi_1=\{-1,1\},\chi_2=\{-1,1\}\tnote{3}$&$(s_2=1,z_{s_4},z_{w_0})$&1\\
			\hline\hline
		\end{tabular}
		\begin{tablenotes}    
			\footnotesize               
			\item[1] the fixed variables are $s_2=0,s_4=0$ and $s_2=2,s_4=0$.
			\item[2] there are no stationary points (the corresponding values) for this case.
			\item[3] the fixed variables are any 4 cases where $\chi_1=-1(1)$ and $\chi_2=-1(1)$.
		\end{tablenotes}            
	\end{threeparttable}       
\end{table} 
	\begin{table}
	\caption{The stationary points and corresponding values of function $H$ when fixing three variables of $s_2$, $s_4$, $w_0$, $\chi_1$ and $\chi_2$.}
	\centering
	\label{TAfix3}  
	\begin{tabular}{ccc}
		\hline \hline
		Fixed variable(s)&Stationary ponit(s)&Function value(s)\\
		\hline
		$\makecell{
			w_0=0(1),\chi_1=-1,\chi_2=-1
			\\w_0=0(1),\chi_1=1,\chi_2=1,
			\\w_0=1,\chi_1=-1(1),\chi_2=1(-1)
		}$&$(s_2=1,z_{s_4})$&1\\
		\hline
		$w_0=0,\chi_1=-1(1),\chi_2=1(-1)$&$(s_2=1,z_{s_4}),(z_{s_2},1)$&1\\
		\hline
		$\makecell{s_4=0,\chi_1=1(-1),\chi_2=1(-1)\\s_4=2,\chi_1=1(-1),\chi_2=1(-1)\\s_4=0,\chi_1=1(-1),\chi_2=-1(1)}$&$(s_2=1,z_{w_0})$&1\\
		\hline
		$s_4=2,\chi_1=1(-1),\chi_2=-1(1)$&$(s_2=1,z_{w_0}),(z_{s_2},w_0=1/2)$&1\\
		\hline
		$s_4=\{0,2\},\chi_1=\{-1,1\},\chi_2=\{-1,1\}$&$(s_2=1,z_{w_0})$&1\\
		\hline
		$s_2=0(2),\chi_1=-1(1),\chi_2=-1(1)$&$(s_4=1,w_0=1/2)$&1\\
		\hline
		$\makecell{s_2=0,\chi_1=-1,\chi_2=1,\\s_2=2,\chi_1=-1(1),\chi_2=1(-1)}$&$(s_4=1/(1-z_{w_0}),z_{w_0})$&1\\
		\hline
		$s_2=0(2),\chi_1=1,\chi_2=1$&-&0\\
		\hline
		$\makecell{s_2=0(2),\chi_1=-1,\chi_2=1,\\s_2=2,\chi_1=1(-1),\chi_2=-1(1)}$&$s_4=1/(1-z_{w_0}),z_{w_0}$&1\\
		\hline
		$s_2=0(2),\chi_1=1,\chi_2=1$&-&0\\
		\hline
		$\makecell{s_2=0,w_0=1,\chi_2=1(-1),\\s_2=2,w_0=1,\chi_2=1(-1)}$&$(s_4=1,\chi_1=0)$&1\\
		\hline
		$\makecell{s_2=0,w_0=0,\chi_2=1(-1),\\s_2=2,w_0=0,\chi_2=-1(1)}$&-&-\\
		\hline
		$s_2=0,s_4=0,\chi_2=1(-1)$&$(z_{w_0},\chi_1=-1),(z_{w_0},\chi_1=1/2)$&0(27/32)\\
		\hline
		$\makecell{s_2=0(2),s_4=2,\chi_2=1(-1),\\s_2=0(2),s_2=2,\chi_2=-1(1)}$&-&-\\
		\hline
		$s_2=2,s_4=0,\chi_2=-1(1)$&-&0\\
		\hline
		$\makecell{s_2=0,s_4=0,w_0=0(1),\\s_2=2(0),s_4=2,w_0=1(0) }$
		&$(z_{\chi_1},z_{\chi_2})$&0\\
		\hline
		$\makecell{s_2=2,s_4=0,w_0=0(1),\\s_2=0(2),s_4=2,w_0=1(0)} $
		&-&0\\
		
		\hline \hline  
	\end{tabular}
\end{table}
\begin{table}
	\caption{The stationary points and corresponding values of function $H$ when fixing four variables of $s_2$, $s_4$, $w_0$, $\chi_1$ and $\chi_2$.}
	\centering
	\label{TAfix4}  
	\begin{threeparttable}  
		\begin{tabular}{ccc}
			\hline \hline
			Fixed variable(s)&Stationary ponit(s)&Function value(s)\\
			\hline
			$s_4=\{0,2\},w_0=\{0,1\},\chi_1=\{-1,1\},\chi_2=\{-1,1\}$&$s_2=1$&1\\		
			\hline
			$\makecell{s_2=0,w_0=0,\chi_1=1(-1),\chi_2=-1(1),\\s_2=2,w_0=0,\chi_1=1(-1),\chi_2=-1(1)}$&$s_4=1$&1\\
			\hline
			$\makecell{s_2=0,w_0=0,\chi_1=1(-1),\chi_2=1(-1),\\s_2=2,w_0=0,\chi_1=1(-1),\chi_2=1(-1),\\
				s_2=\{0,2\},w_0=1,\chi_1=\{-1,1\},\chi_2=\{-1,1\}\tnote{4}}$&$z_{{s_4}}$&0\\
			\hline
			$\makecell{s_2=0,s_4=2,\chi_1=1(-1),\chi_2=-1(1),\\s_2=2,s_4=2,\chi_1=1(-1),\chi_2=-1(1)}$&$w_0=1/2$&1\\
			\hline
			$\makecell{s_2=0,s_4=2,\chi_1=1(-1),\chi_2=1(-1),\\s_2=2,s_4=2,\chi_1=1(-1),\chi_2=1(-1),\\
				s_2=\{0,2\},s_4=0,\chi_1=\{-1,1\},\chi_2=\{-1,1\}}$&$z_{w_0}$&0\\
			\hline
			$\makecell{s_2=0,s_4=0,w_0=0(1),\chi_1=-1\\s_2=0(2),s_4=2,w_0=0(1),\chi_1=-1}$&$\chi_2=1/2,-1$&$\makecell{27/32\\0}$\\
			\hline
			$\makecell{s_2=0,s_4=2,w_0=1,\chi_1=-1(1),\\s_2=2,s_4=0,w_0=0,\chi_1=-1(1),\\s_2=2,s_4=\{0,2\},w_0=1,\chi_1=\{-1,1\}}$&$\chi_2=0$&0\\
			\hline
			$\makecell{s_2=2(0),s_4=2(0),w_0=1,\chi_1=1,\\s_2=0,s_4=0(2),w_0=0,\chi_2=1}$&$\chi_1=-1/2,1$&\makecell{27/32\\0}\\
			\hline \hline  
		\end{tabular}
		\begin{tablenotes}    
			\footnotesize               
			\item[4] the fixed variables are $w_0=1$ with any 8 cases where $s_2=0(2),\chi_1=-1(1)$ and $\chi_2=-1(1)$.
		\end{tablenotes}            
	\end{threeparttable}       
\end{table}
\textbf{Step 4. The proof of Eq. (\ref{neq25})}

We first introduce the function $K(s_2,s_4,w_0,\chi_1,\chi_2)=p_1^2-p_1p_3-p_4^2+p_2p_4$ and express it as
\begin{equation}\label{FK}
	\begin{aligned}
		&K=\Big[
		\frac{(\chi_1+\chi_2)^2(s_2 - 1)^4(s_4 - 1)^2}{4}- 1\Big]^2 - \Bigg[(s_2 - 1)(s_4 - 1)-\frac{(\chi_1+\chi_2)(s_2 - 1)(s_4w_0(s_2 - 1) - s_2 + 1)}{2}\\
		&\qquad+ \chi_1\chi_2\Big(s_4w_0(s_2 - 1)^2 - (s_2 - 1)(s_2 - 2)(s_4 - 1)\Big)+\frac{(\chi_1+\chi_2)(s_2 - 1)^2(s_4 - 1)}{4}\times \\
		&\qquad\Big(2s_4(w_0 - 1) - (\chi_1+\chi_2)(2s_2 + s_4w_0 - 3) + 2\Big)
		\Bigg]^2 + \Big[\frac{(\chi_1+\chi_2)^2(s_2 - 1)^4(s_4 - 1)^2}{4} - 1\Big]\times\\
		&\qquad\Bigg[(s_4w_0 - 1)(s_2 - 1) +\Big (\frac{(2s_2 - 3)(s_4 - 1)}{2} - \frac{s_4w_0(2s_2 - 2)}{2}\Big)
		(\chi_1+\chi_2) +(\chi_1+\chi_2)(s_2 - 1)^2(s_4 - 1)\times\\
		&\qquad\quad\Big[\frac{(s_4w_0 - 1)(s_2 - 1) }{2}+\Big( \frac{(2s_2 - 3)(s_4 - 1)}{4} - \frac{s_4w_0(2s_2 - 2)}{4}\Big )(\chi_1+\chi_2)\\
		&\qquad\quad+\frac{\chi_1\chi_2(s_2 - 1)(s_2 + s_4w_0 - 2)}{2}\Big]+\chi_1\chi_2(s_2 - 1)(s_2 + s_4w_0 - 2)\Bigg]+ \Bigg[(s_2 - 1)(s_4 - 1) \\
		&\qquad- \frac{(\chi_1+\chi_2)(s_2 - 1)(s_4w_0(s_2 - 1) - s_2 + 1)}{2}+\chi_1\chi_2\Big(s_4w_0(s_2 - 1)^2 - (s_2 - 1)(s_2 - 2)(s_4 - 1)\Big)\\
		&\qquad+ \frac{(\chi_1+\chi_2)(s_2 - 1)^2(s_4 - 1)\Big(2s_4(w_0 - 1) - (\chi_1+\chi_2)(2s_2 +s_4w_0 - 3) + 2\Big)}{4}\Bigg]\times\\
		&\qquad\Bigg[s_4(w_0 - 1) - (\chi_1+\chi_2)\frac{2s_2 + s_4w_0 - 3}{2} + (\chi_1+\chi_2)(s_2 - 1)^2(s_4 - 1)\times\\
		&\qquad\bigg(\frac{(s_2 - 1)(s_4 - 1)}{2} - \frac{(\chi_1+\chi_2)(s_2 - 1)\Big(s_4w_0(s_2 - 1) - s_2 + 1\Big)}{2}\\
		&\qquad\qquad\quad+ \chi_1\chi_2\frac{s_4w_0(s_2 - 1)^2 - (s_2 - 1)(s_2 - 2)(s_4 - 1)}{2}\bigg)+ 1\Bigg]\:\:{\rm in}\:\:\Omega.
	\end{aligned}
\end{equation}
To prove that the function $K$ (\ref{FK}) is strictly larger than zero in the open region $\Omega$, similarly, we first show that there are no stationary points of the function $K$ (\ref{FK}) in the open region $\Omega$ through some manipulations. We then consider the points at $\partial \overline{\Omega}$, and have
\begin{equation}
	\begin{aligned}
		&K(s_4=0)=\underbrace{(2-s_2)\big(\frac{\chi_1+\chi_2}{2}(s_2-1)^2+1\big)}_{>0}\times
		\underbrace{[\chi_1\chi_2(1-s_2)+\frac{\chi_1+\chi_2}{2}(s_2-2)+1]}_{>0}\times K_{1}\:\:{\rm in}\:\:\Omega_{s_4=0},
	\end{aligned}
\end{equation}
where
\begin{equation}\label{K1}
	K_1=\chi_1\chi_2(s_2-1)(2-s_2)+s_2+\frac{(\chi_1+\chi_2)^2}{4}\big((2-s_2)(2-s_2+s_2^3)-2\big).
\end{equation} 
Now we focus on the function $K_1(s_2,\chi_1,\chi_2)$ that has no stationary points in the open region $\Omega_{s_4=0}$, and the maximum value must be obtained at the boundary points $\partial \overline{\Omega}_{s_4=0}$. Based on this fact, we have
\begin{subequations}
	\begin{align}
		&K_1(s_2=0)=\underbrace{2(\chi_1\chi_2+1-\chi_1-\chi_2)}_{> 0}\:\:{\rm in}\:\:\Omega_{s_4=0,s_2=0},\\
		&K_1(s_2=2)=\underbrace{2(1-\chi_1\chi_2)(1-\frac{(\chi_1+\chi_2)^2}{4})}_{>0}\:\:{\rm in}\:\:\Omega_{s_4=0,s_2=2},\\
		&K_1(\chi_1=-1)=\underbrace{(2-s_2)\frac{\big(s_2+1+\chi_2(1-s_2)\big)\big(s_2(2-s_2)+\chi_2(w_2-1)^2+1\big)}{16}}_{>0}\times\notag\\
		&\quad\quad\quad\quad\quad\quad\underbrace{\big(s_2-4+\chi_2(4-3s_2)\big)}_{>0}K_{11}\:\:{\rm in}\:\:\Omega_{s_4=0,\chi_1=-1},\\
		&K_1(\chi_1=1)=\underbrace{\frac{\big(s_2(1+\chi_2)+1-\chi_2\big)\big((1+\chi_2)(s_2-1)^2+2\big)}{16}}_{>0}\times\underbrace{s_2(2-s_2)(\chi_2-1)}_{<0}K_{12}\:\:{\rm in}\:\:\Omega_{s_4=0,\chi_1=1},
	\end{align}
\end{subequations}
where 
\begin{subequations}
	\begin{align}
		&K_{11}=2s_2-3s_2^2+s_2^3-2+\chi_2(2s_2-2+s_2^2-s_2^3),\\
		&K_{12}=\chi_2(2-s_2^2)(1-s_2)+2(s_2-1)+s_2^2(s_2-3).
	\end{align}
\end{subequations}
It is easy to show $K_{11}\leq0$ in $\overline{\Omega}_{s_4=0,\chi_1=-1}$, $K_{12}\leq0$ in $\overline{\Omega}_{s_4=0,\chi_1=1}$ through some simple analysis, then we can obtain that the function $K_1$ (\ref{K1}) has no stationary points in the open region $\Omega_{s_4=0}$ and is no less than zero at $\partial\overline{\Omega}$, i.e., $K(s_4=0)\geq0$ at $\partial\overline{\Omega}_{s_4=0}$.

We would also like to point out that the discussion on the other boundary points is similar, i.e., $K(s_4=2), K(s_2=0)$, $K(s_2=2)$, $K(w_0=0$), $K(w_0=1)$, $K(\chi_1=-1)$, $K(\chi_1=1)$, and the details can be found in \ref{A4}. Therefore, we prove the condition (\ref{neq25}), and also Eq. (\ref{2neq}) in the region $\Omega$. 

For the case of $\chi_{1,2} =\pm 1$, one can adopt the reductive approach \cite{Miller} to obtain the roots of the polynomial $q(\lambda)$ (\ref{qcharc}) satisfying $|\lambda_k|\leq1$ ($k$=1,2 and 3). In another way, it is also obvious that the roots of the polynomial $q(\lambda)$ are continuous functions of $\chi_{1,2}$, and hence the roots of polynomial $q(\lambda)$ will satisfy the condition $|\lambda_k|\leq1$ in $\chi_{1,2}\in[-1,1]$ once the roots satisfy the condition $|\lambda_k|\leq1$ in $\chi_{1,2}\in(-1,1)$.

From the above discussion, one can see that the roots of the characteristic polynomial $p(\lambda)$  satisfy the condition $|\lambda_k|\leq 1 $ ($k$ = 1, 2,3 and 4), thus the present MRT-LB model (\ref{distributionf}) and FLFD scheme (\ref{fds5=1}) are both unconditionally stable.

	\section{Numerical results and discussion}\label{secNumer}
In this section, we conduct some simulations to test the present MRT-LB model (\ref{distributionf}) and FLFD scheme (\ref{fds5=1}) for the two-dimensional diffusion equation (\ref{equde}) , and the weight coefficient $w_0$, the relaxation parameters $s_2$ and $s_4$ are determined from Eq. (\ref{solution}). To measure the accuracy of the MRT-LB model (\ref{distributionf}) and FLFD scheme (\ref{fds5=1}), the following root-mean-square error ($RMSE$) is adopted \cite{Kruger2017},
\begin{equation}\label{RMSE}
	RMSE = \sqrt{\frac{\sum_{i=1}^{Nx}\sum_{j=1}^{N_y}[\phi(i\Delta x, j\Delta x,n\Delta t)-\phi^*(i\Delta x, j\Delta x,n\Delta t)]^2}{N_x\times N_y}},
\end{equation}
where $N_x\times N_y$ is the number of grid points, $\Delta x=\Delta y = L_x/N_x=L_y/N_y$ is lattice spacing ($L_x$ and $L_y$ are length  in $x$ and $y$ directions), $\phi$ and $\phi^{*}$ are the numerical and analytical solutions, respectively. Based on the definition of $RMSE$, one can estimate the convergence rate $(CR)$ of numerical scheme with the following formula,
\begin{equation}\label{CR}
	CR = \frac{ \log( RMSE_{\Delta x} /RMSE_{\Delta x/2} ) } {\log2}.
\end{equation}

\textbf{Example 1.}
We first consider the diffusion equation (\ref{equde}) with the following inhomogeneous initial and periodic boundary conditions,
\begin{equation}\label{Ex1}
	\begin{aligned}
		\left\{
		\begin{array}{lr}
			\partial_t\phi = \kappa(\partial^2_{x}\phi+\partial^2_{y}\phi)+\pi^2 ,  (x,y)\in {\rm D},&  \\
			\phi(x,y,t=0)=\sin(\pi x)\sin(\pi y) ,  &  \\
		\end{array}
		\right.
	\end{aligned}
\end{equation}
and obtain the solution as
\begin{equation}\label{ExSo1}
	\phi(x,y,t) = \sin(\pi x)\sin(\pi y)\exp(-2\kappa\pi^2 t)+\pi^2 t.
\end{equation}

\begin{figure}
	\centering
	\subfigure[$\:$MRT-LB model]
	{
		\begin{minipage}{0.40\linewidth}
			\centering
			\includegraphics[width=2.5in]{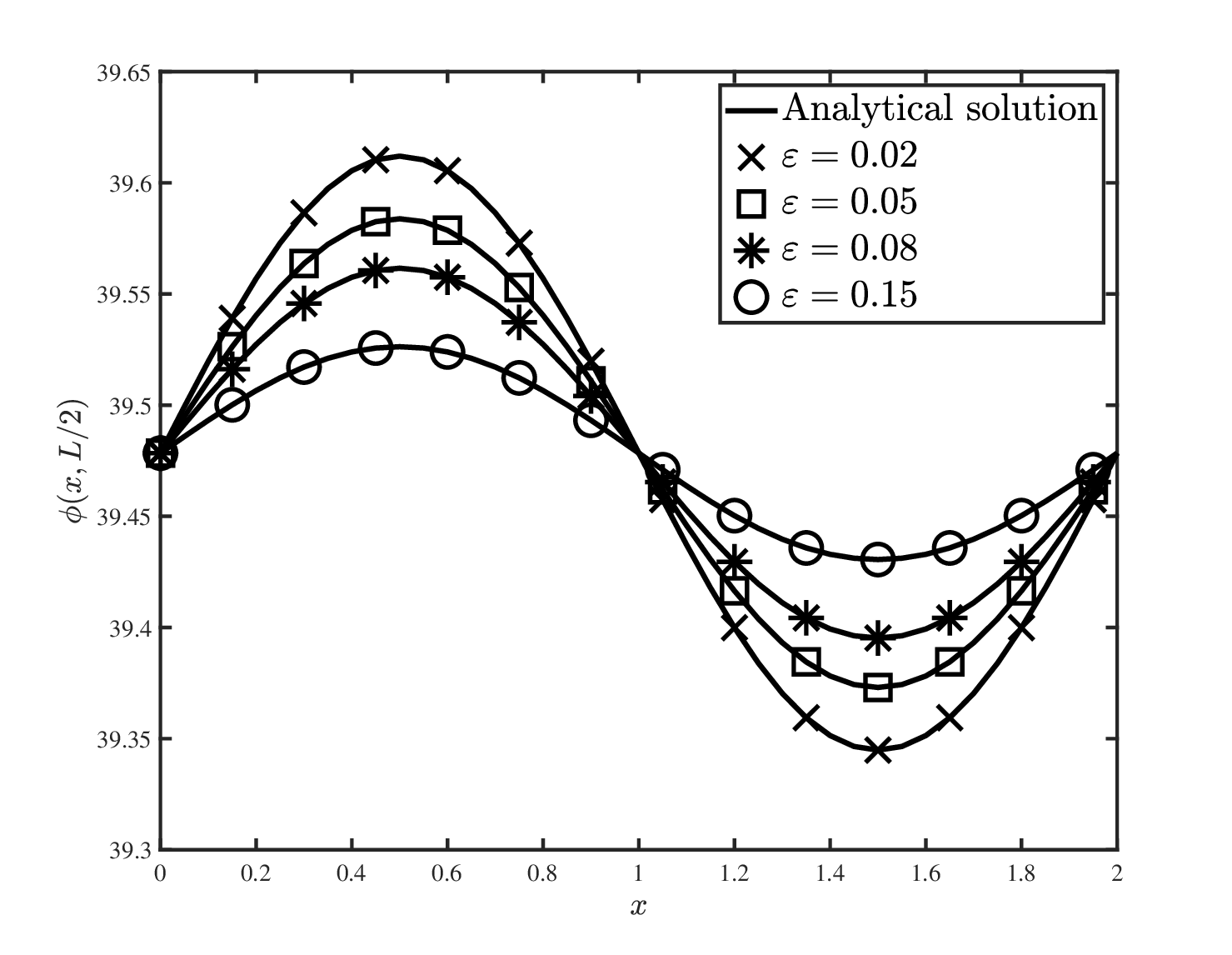}
		\end{minipage}
	}
	\subfigure[$\:$FLFD scheme]
	{
		\begin{minipage}{0.40\linewidth}
			\centering
			\includegraphics[width=2.5in]{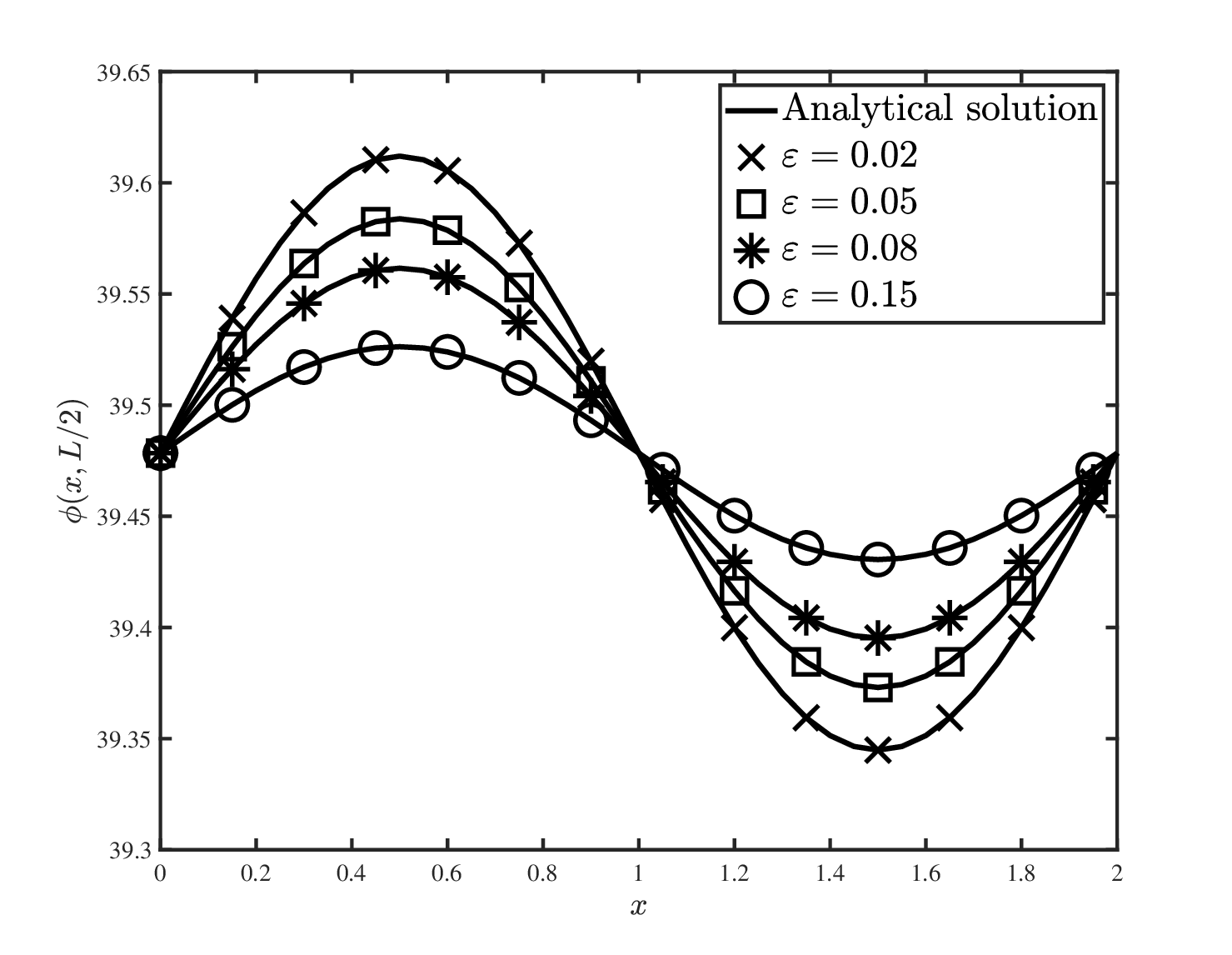}
		\end{minipage}
	}
	\caption{The numerical and analytical solutions under different values of discretization parameter $\varepsilon$ ($T =4.0$).}
	\label{Ex1Figure1}
	\setlength{\belowcaptionskip}{-1cm}
\end{figure}

\begin{figure}
	\centering
	\subfigure[$\:$MRT-LB model]
	{
		\begin{minipage}{0.40\linewidth}
			\centering
			\includegraphics[width=2.5in]{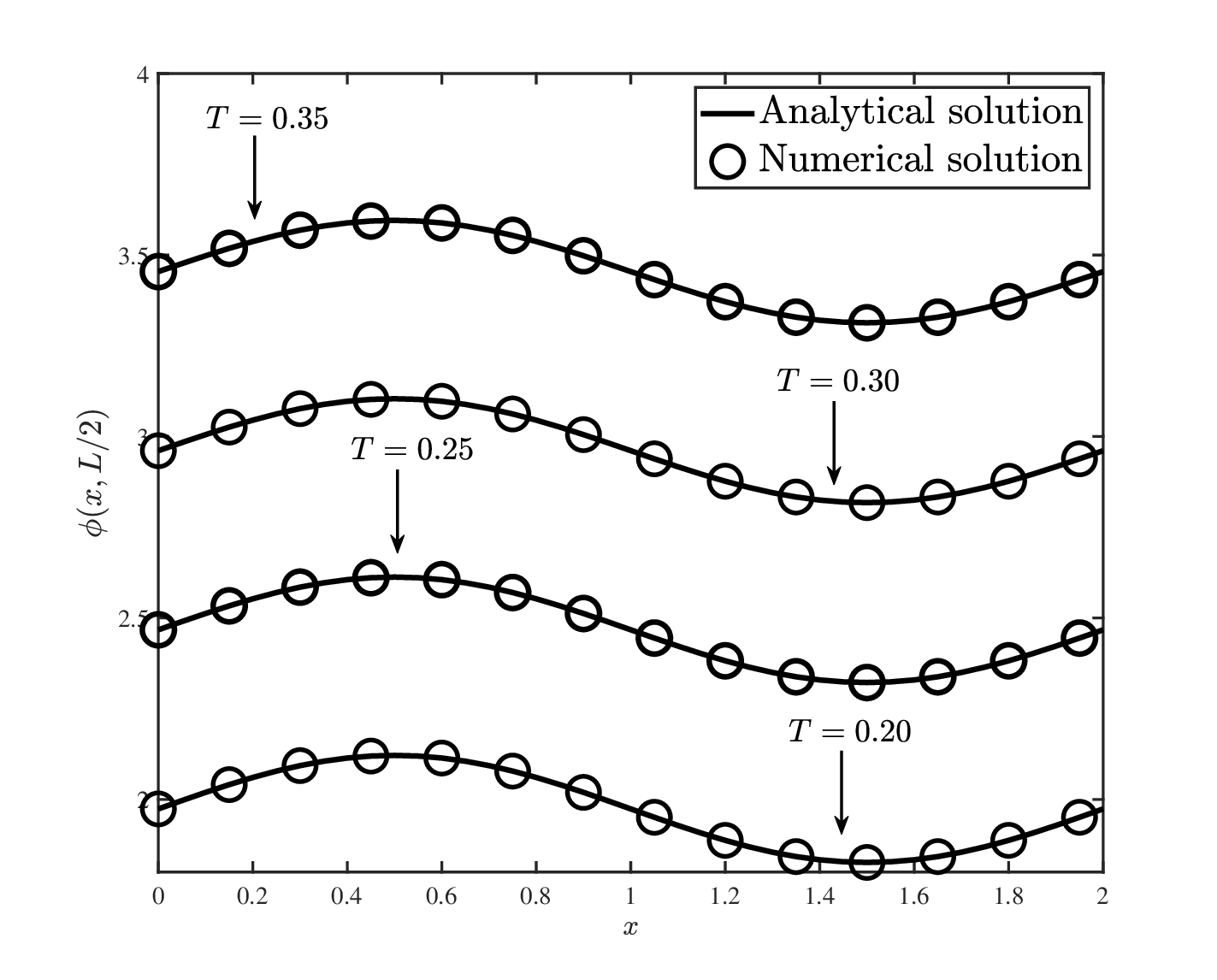}
		\end{minipage}
	}
	\subfigure[$\:$FLFD scheme]
	{
		\begin{minipage}{0.40\linewidth}
			\centering
			\includegraphics[width=2.5in]{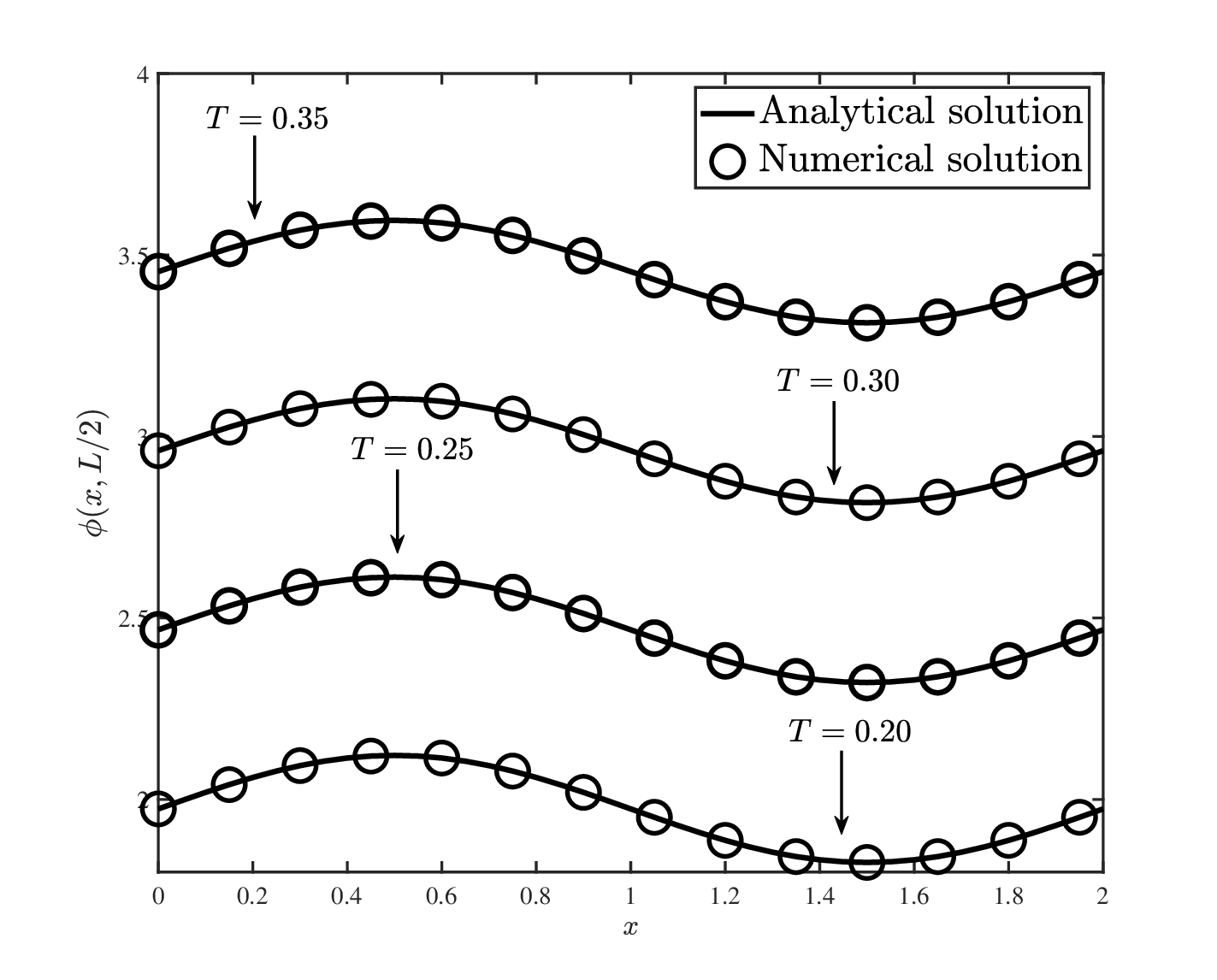}
		\end{minipage}
	}
	\caption{The numerical and analytical solutions under different values of time $T$ ($\varepsilon=0.15$).}
	\label{Ex1Figure2}
	\setlength{\belowcaptionskip}{-1cm}
\end{figure}

\begin{figure}
	\centering
	\subfigure[$\:$Analytical solution]
	{
		\begin{minipage}{0.40\linewidth}
			\centering
			\includegraphics[width=2.5in]{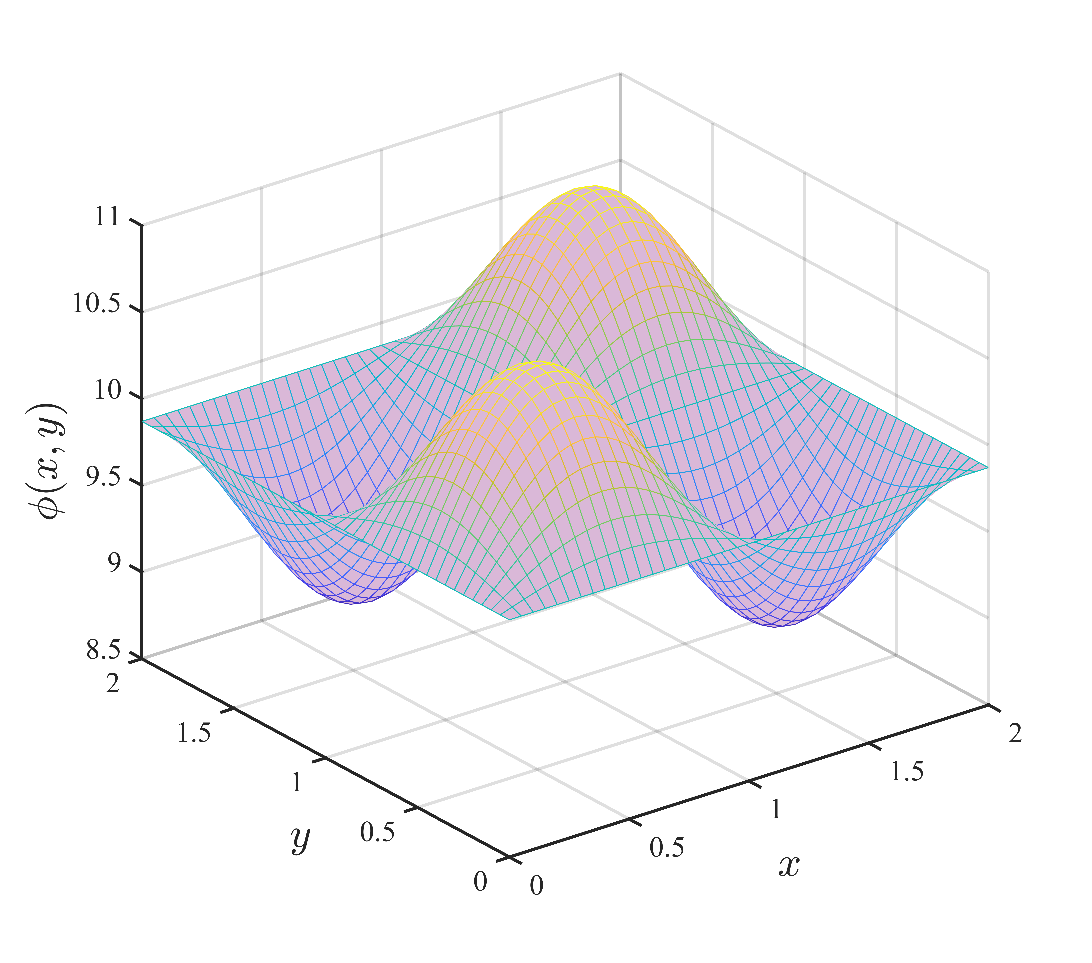}
		\end{minipage}
	}

	\subfigure[$\:$MRT-LB model]
	{
		\begin{minipage}{0.40\linewidth}
			\centering
			\includegraphics[width=2.5in]{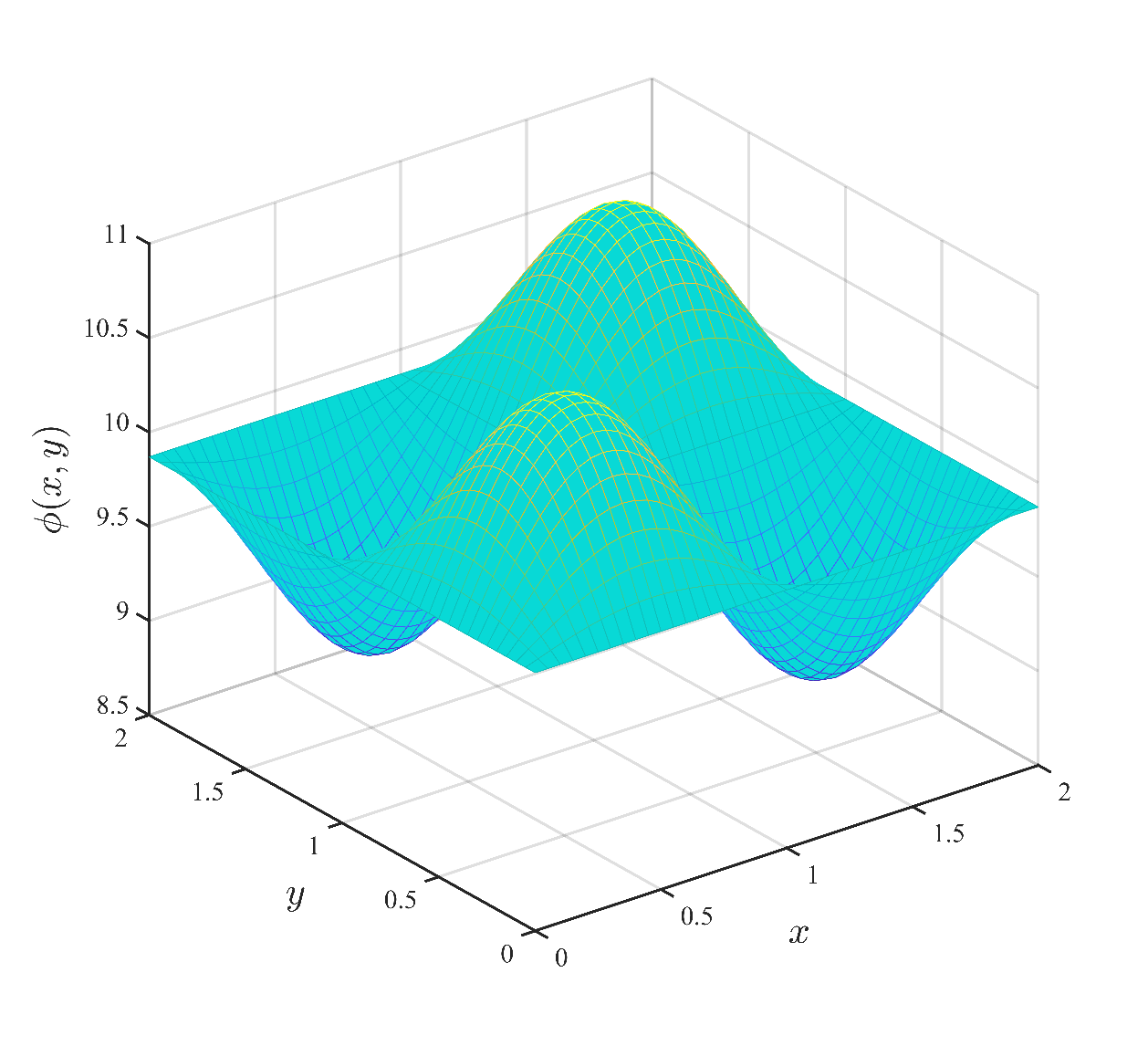}
		\end{minipage}
	}
	\subfigure[$\:$FLFD scheme]
	{
		\begin{minipage}{0.40\linewidth}
			\centering
			\includegraphics[width=2.5in]{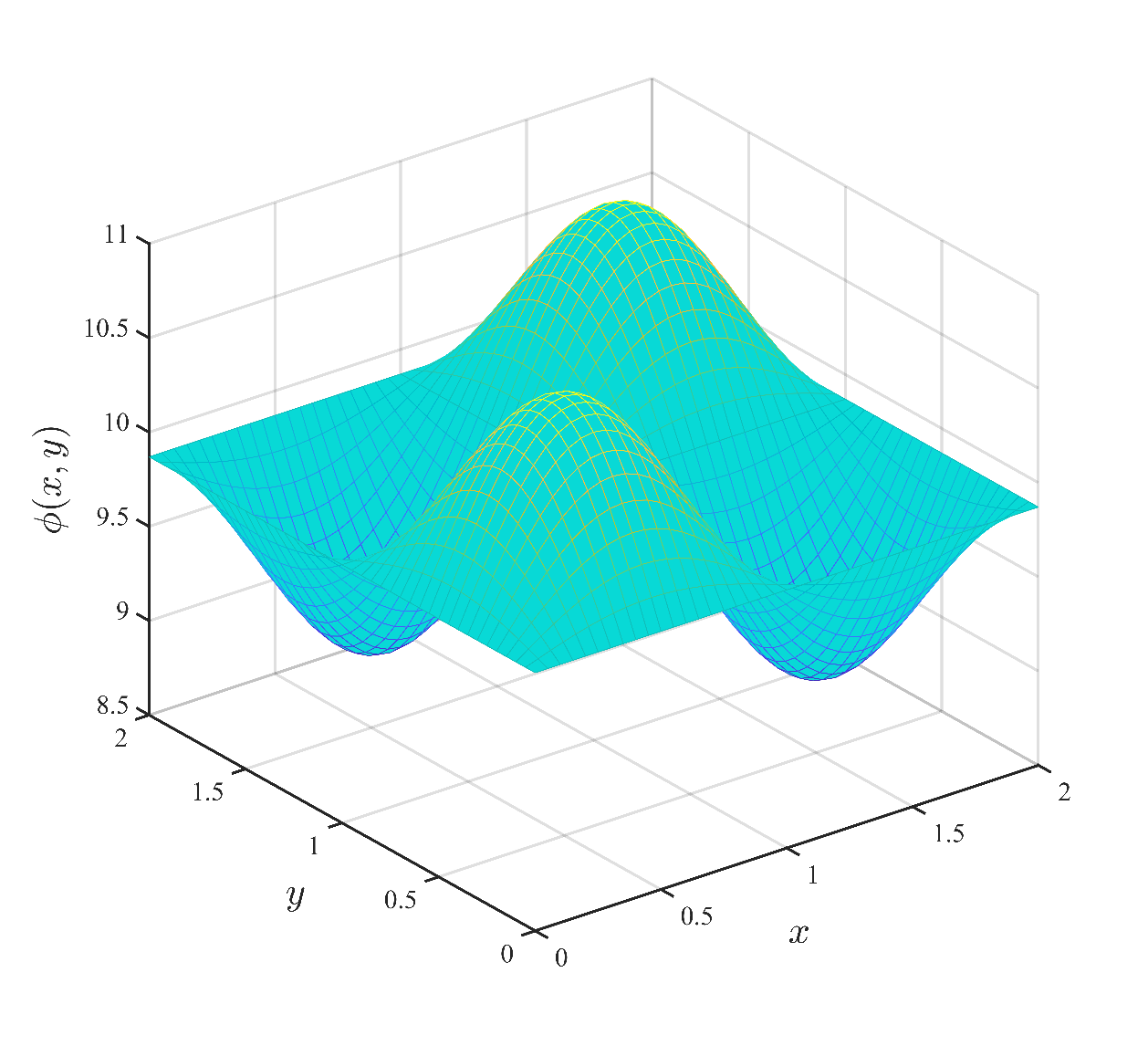}
		\end{minipage}
	}
	\caption{Distributions of the variable $\phi$ at the discretization parameter $\varepsilon=0.15$ and time $T=1.0$.}
	\label{Ex1Figure3}
	\setlength{\belowcaptionskip}{-1cm}
\end{figure}

\begin{figure}
	\centering
	\subfigure[$\:$MRT-LB model]
	{
		\begin{minipage}{0.40\linewidth}
			\centering
			\includegraphics[width=2.5in]{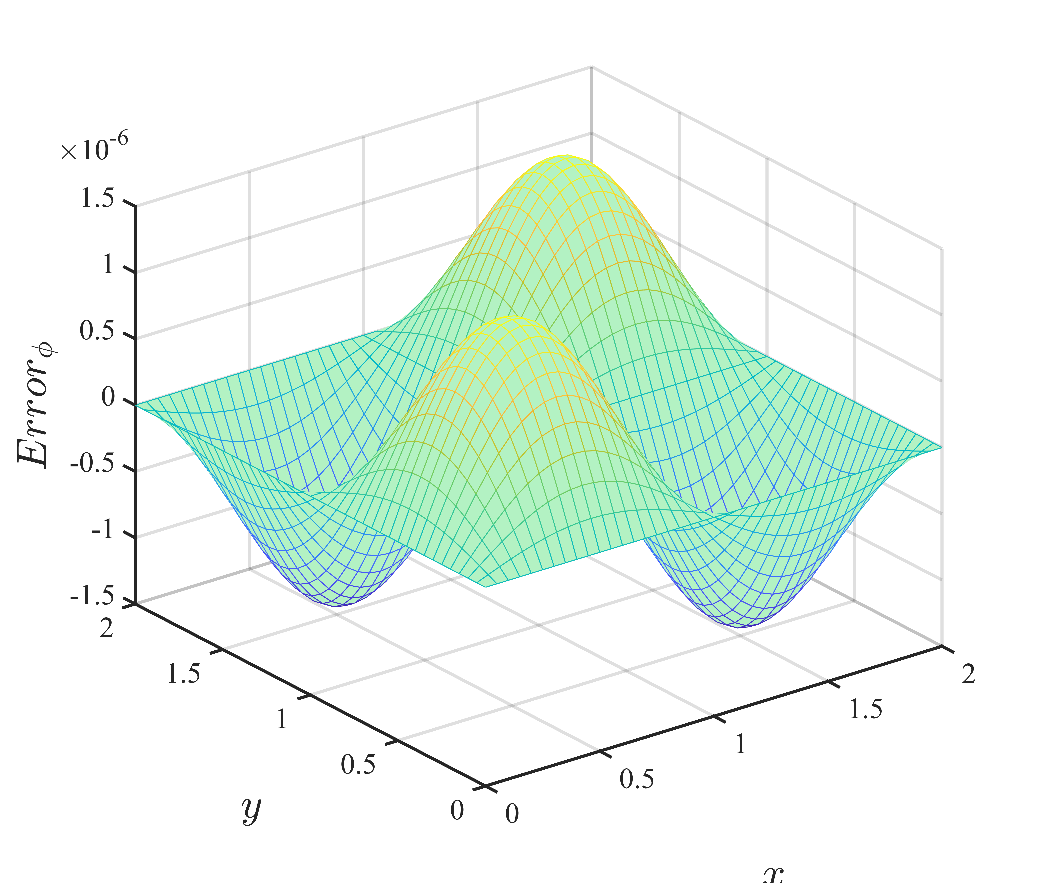}
		\end{minipage}
	}
	\subfigure[$\:$FLFD scheme]
	{
		\begin{minipage}{0.40\linewidth}
			\centering
			\includegraphics[width=2.5in]{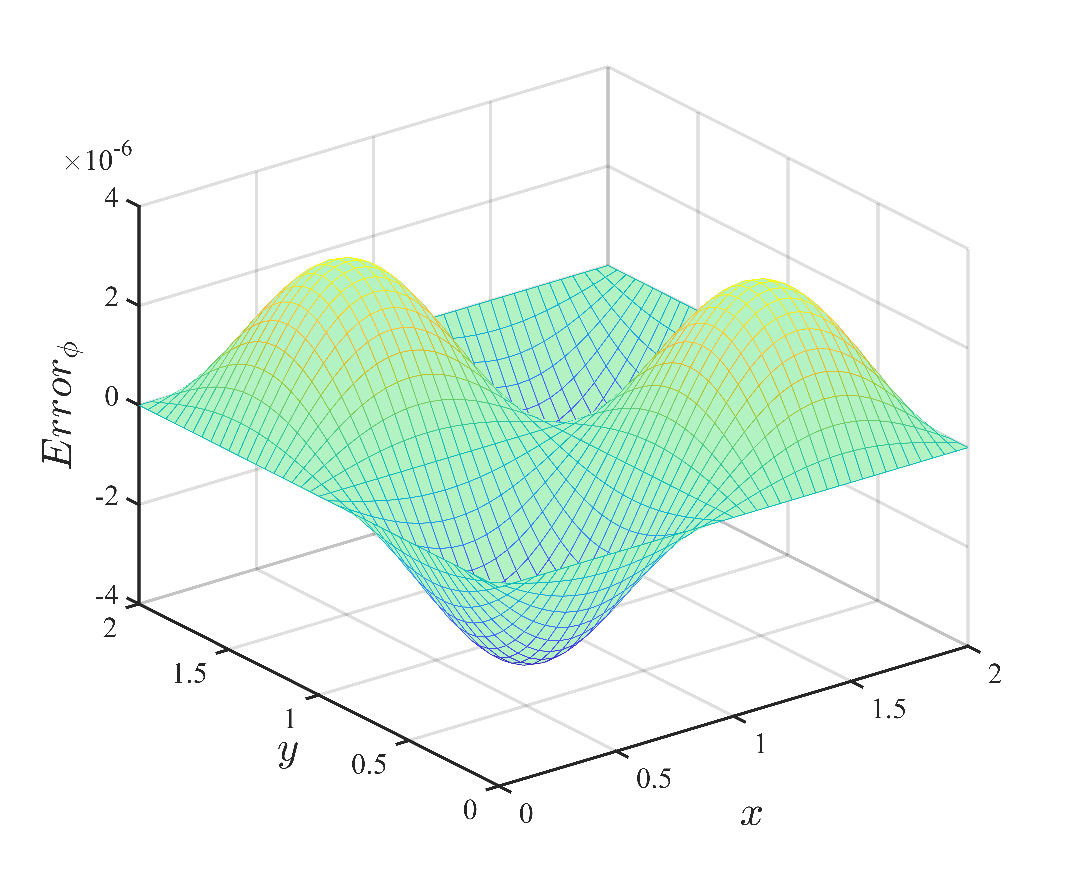}
		\end{minipage}
	}
	\caption{Distributions of the error (${\rm Error}_{\phi}$) at the discretization parameter $\varepsilon=0.15$ and time $T=1.0$.}
	\label{Ex1Figure4}
	\setlength{\belowcaptionskip}{-1cm}
\end{figure}

In our simulations, the computational domain is set to be ${\rm D}=[0,L]\times [0,L]$ with $L=2$, the specified lattice spacing $\Delta x=1/20$ and time step $\Delta t=1/40$. We first carry out some simulations under different values of discretization parameter $\varepsilon$ ($\varepsilon =0.02,0.05,0.08,0.15)$, and plot the profiles of the scalar variable $\phi$ in Fig. \ref{Ex1Figure1} where $T = 4.0$. To see the evolution of the variable $\phi$ in time, we also conduct some simulations under different values of time $T$ ($T=0.20,0.25,0.30,0.25$), and present the results in Fig. \ref{Ex1Figure2} where the discretization parameter is $\varepsilon = 0.15$. From these two figures, one can observe that the numerical results of MRT-LB model and FLFD scheme are in good agreement with the corresponding analytical solutions. In addition, we also 
present the distributions of variable $\phi$ and the absolute errors (${\rm Error}_{\phi}(x_i,y_j):=\phi(x_i,y_j)-\phi^*(x_i,y_j)$, $i,j=0,1,\ldots, L/\Delta x$) in Figs. \ref{Ex1Figure3} and \ref{Ex1Figure4} where $\varepsilon = 0.15$ and $T=1.0$, and find the maximum error of the MRT-LB model is less than $1.5\times 10^{-6}$, while that of the FLFD scheme is less than $4.0\times 10^{-6}$. 

\begin{figure}
	\centering
	\subfigure[$\:$MRT-LB model]
	{
		\begin{minipage}{0.40\linewidth}
			\centering
			\includegraphics[width=2.5in]{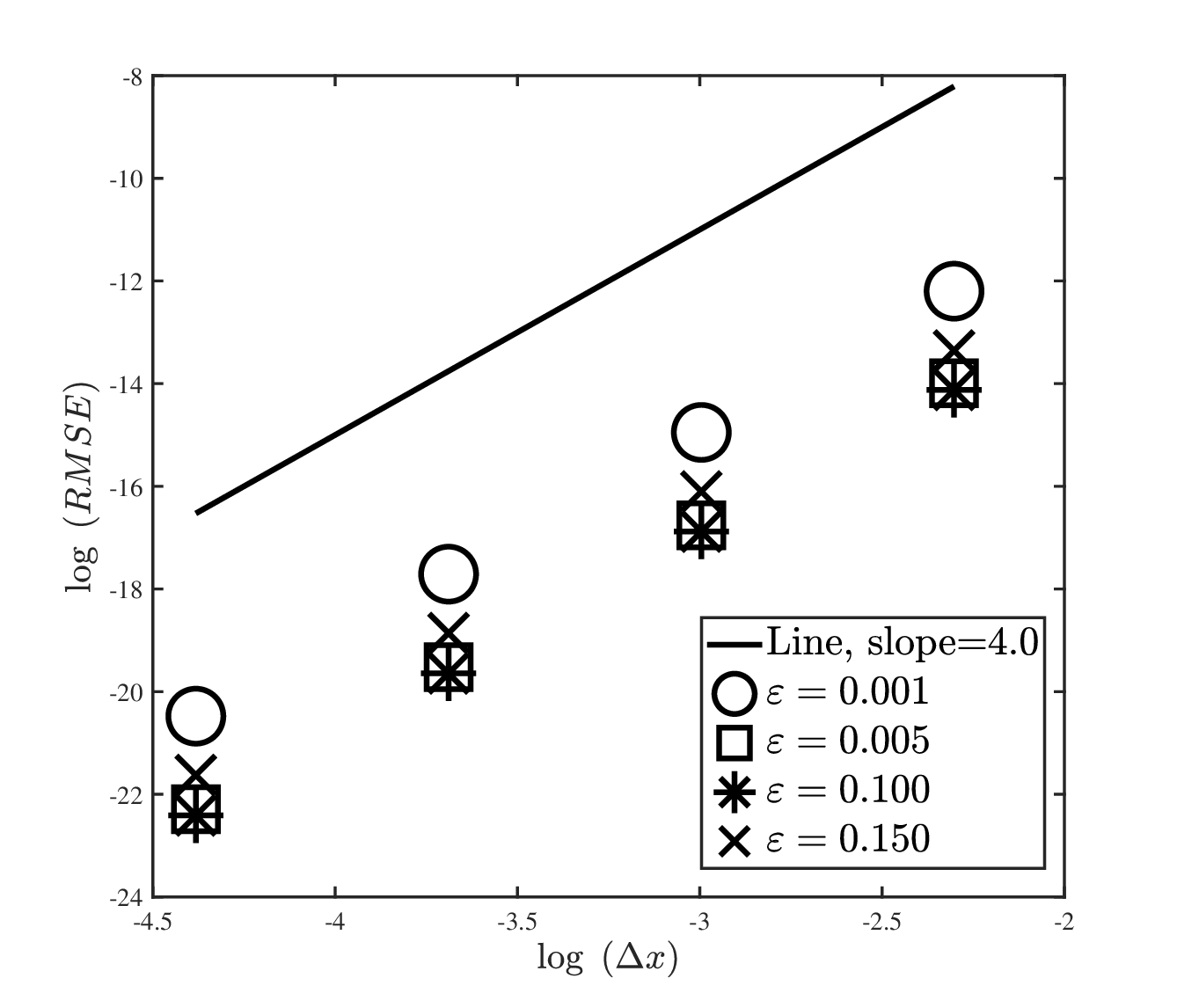}
		\end{minipage}
	}
	\subfigure[$\:$FLFD scheme]
	{
		\begin{minipage}{0.40\linewidth}
			\centering
			\includegraphics[width=2.5in]{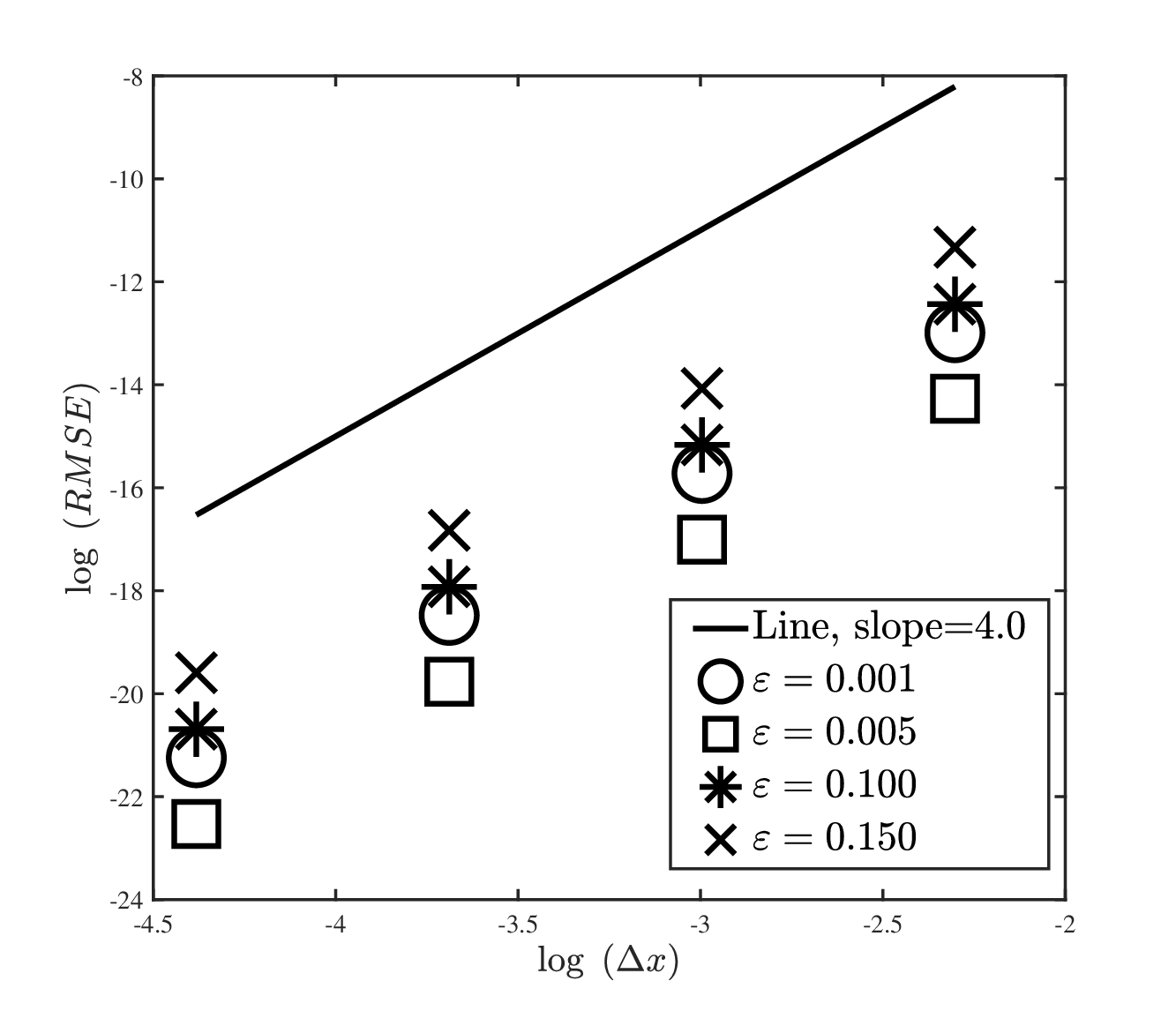}
		\end{minipage}
	}
	\caption{The convergence rates of MRT-LB model and FLFD scheme with different values of discretization parameter $\varepsilon$ ($T=10.0$).}
	\label{Ex1Slope1}
\end{figure}

Finally, to test the convergence rates of the MRT-LB model and FLFD scheme, some simulations under different values of lattice spacing $\Delta x$ are performed, and the time step $\Delta t$ is determined by a fixed constant $\Delta  x^2/\Delta t$. We measure the $RMSEs$ between the numerical and analytical solutions, and cauclate the $CRs$ at the time $T=10.0$ and different values of discretization parameter $\varepsilon$ in Fig. \ref{Ex1Slope1}, Tables \ref{TEx1CRLB} and \ref{TEx1CRFD}. From these results, one can observe that both the MRT-LB model and FLFD scheme have a fourth-order convergence rate in space.
\begin{table}
	\begin{center}
		\caption{The $RMSEs$ and $CRs$ of MRT-LB model under different values of discretization parameter $\varepsilon$ ($T=10.0$, $\Delta x=1/10$, $\Delta t=1/10$).}
		\label{TEx1CRLB}  
		\begin{tabular}{ccccccccc}
			\hline \hline
			$\varepsilon$&$RMSE_{\Delta x,\Delta t}$&$RMSE_{\Delta x/2,\Delta t/4}$&
			$RMSE_{\Delta x/4,\Delta t/16}$&
			$RMSE_{\Delta x/8,\Delta t/64}$&$CR$\\
			\midrule[1pt]
			0.001 &  $6.7249\times 10^{-7}$ &$4.3010\times 10^{-8}$ &$2.7205\times 10^{-9}$&$1.7092\times 10^{-10}$&$\sim 3.9807$ \\
			
			0.005&$2.7793\times 10^{-6}$&$1.7768\times 10^{-7}$&$1.1238\times 10^{-8}$&$7.0653\times 10^{-10}$&$\sim 3.6472$\\
			
			0.100&$7.3601\times 10^{-7}$&$4.6588\times 10^{-8}$&$2.9392\times 10^{-9}$&$1.8473\times 10^{-10}$&$\sim 3.9867$\\
			
			0.150&$1.5830\times 10^{-6}$&$1.0176\times 10^{-7}$&$6.4425\times 10^{-9}$&$4.0519\times 10^{-10}$&$\sim 3.9773$\\
			\hline \hline
		\end{tabular}
	\end{center}
\end{table}
\begin{table}
	\begin{center}
		\caption{The $RMSEs$ and $CRs$ of FLFD scheme under different values of discretization parameter $\varepsilon$ ($T=10.0$, $\Delta x=1/10$, $\Delta t=1/10$).}
		\label{TEx1CRFD}  
		\begin{tabular}{ccccccccc}
			\hline \hline
			$\varepsilon$&$RMSE_{\Delta x,\Delta t}$&$RMSE_{\Delta x/2,\Delta t/4}$&
			$RMSE_{\Delta x/4,\Delta t/16}$&
			$RMSE_{\Delta x/8,\Delta t/64}$&$CR$\\
			\midrule[1pt]
			0.001 &  $3.5341\times 10^{-7}$ &$2.3026\times 10^{-8}$ &$1.4631\times 10^{-9}$&$9.2142\times 10^{-11}$&$\sim 3.9684$ \\
			
			0.005&$1.4207\times 10^{-6}$&$9.2523\times 10^{-8}$&$5.8786\times 10^{-9}$&$3.7000\times 10^{-10}$&$\sim 3.9689$\\
			
			0.100&$3.9773\times 10^{-6}$&$2.5877\times 10^{-7}$&$1.6437\times 10^{-8}$&$1.0348\times 10^{-9}$&$\sim 3.9684$\\
			
			0.150&$1.2018\times 10^{-5}$&$7,7574\times 10^{-7}$&$4.9177\times 10^{-8}$&$3.0941\times 10^{-9}$&$\sim 3.9744$\\
			\hline \hline
		\end{tabular}
	\end{center}
\end{table}

\textbf{Example 2.}
In the second example, we consider the two-dimensional equation (\ref{equde}) with the following initial and the Dirichlet boundary conditions,
\begin{subequations}\label{Ex2}
	\begin{align}
		&\phi(x,y,0)=\exp[-(x+y)],0\leq x,y\leq 1,\\
		&\phi(x,0,t)=\exp(-x+2\kappa t),\phi(1,y,t)=\exp[-(1+y)+2\kappa t],t\geq 0\label{Ex2b}\\
		&	\phi(x,1,t)=\exp[-(1+x)+2\kappa t],\phi(0,y,t)=\exp(-y+2\kappa t),t\geq 0.\label{Ex2c}
	\end{align}
\end{subequations}
and derive the analytical solution as
\begin{equation}
	\phi(x,y,t)=\exp[-(x+y)+2\kappa t].
\end{equation}

We now give a special discussion on how to implement the Dirichlet boundary conditions in the MRT-LB model with the D2Q5 lattice structure. Actually, one can first determine the unique unknown distribution function $f_i$ on four boundaries (excluding four concer points) accurately with the Dirichlet boundary conditions (\ref{Ex2b}) and (\ref{Ex2c}). However, there are two unknown distributions $f_i$ at the four concer points (e.g., $f_1$ and $f_4$ are unknown at the upper-left concer point) that can be calculated with the aid of Eq. (\ref{boundary}) and the Dirichlet boundary conditions (\ref{Ex2b}) and (\ref{Ex2c}). For instance, for the upper-left corner point $\bm{r}_{ul}=[0,1]^T$, we can compute the spatial derivatives of $\bm{f}^{eq}$ ($\partial_x\bm{f}^{eq}|_{\bm{r}_{ul},t}$ and $\partial_y\bm{f}^{eq}|_{\bm{r}_{ul},t}$) with the the Dirichlet condition (\ref{Ex2c}).

\begin{figure}
	\centering
	\subfigure[$\:$MRT-LB model]
	{
		\begin{minipage}{0.40\linewidth}
			\centering
			\includegraphics[width=2.5in]{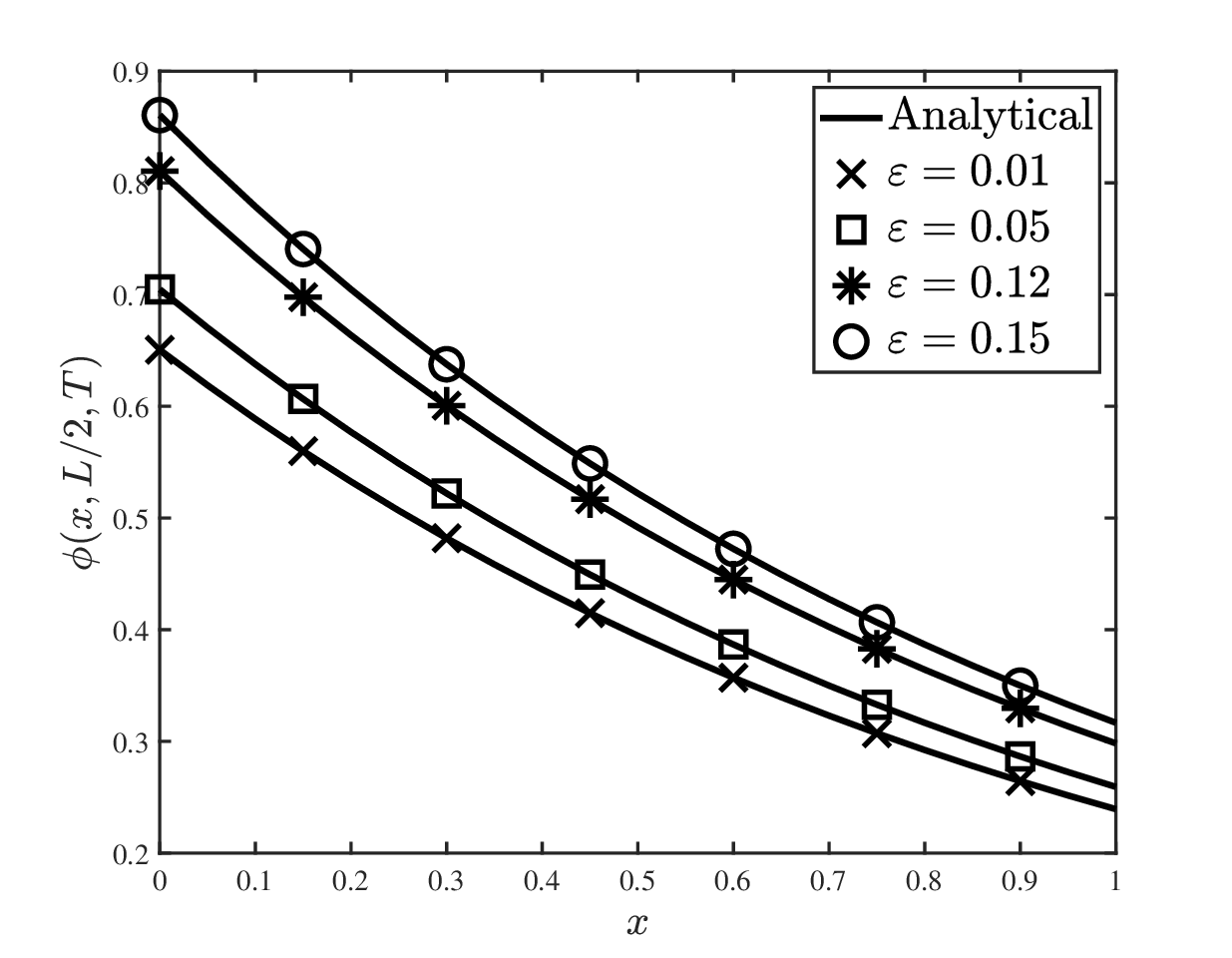}
		\end{minipage}
	}
	\subfigure[$\:$FLFD scheme]
	{
		\begin{minipage}{0.40\linewidth}
			\centering
			\includegraphics[width=2.5in]{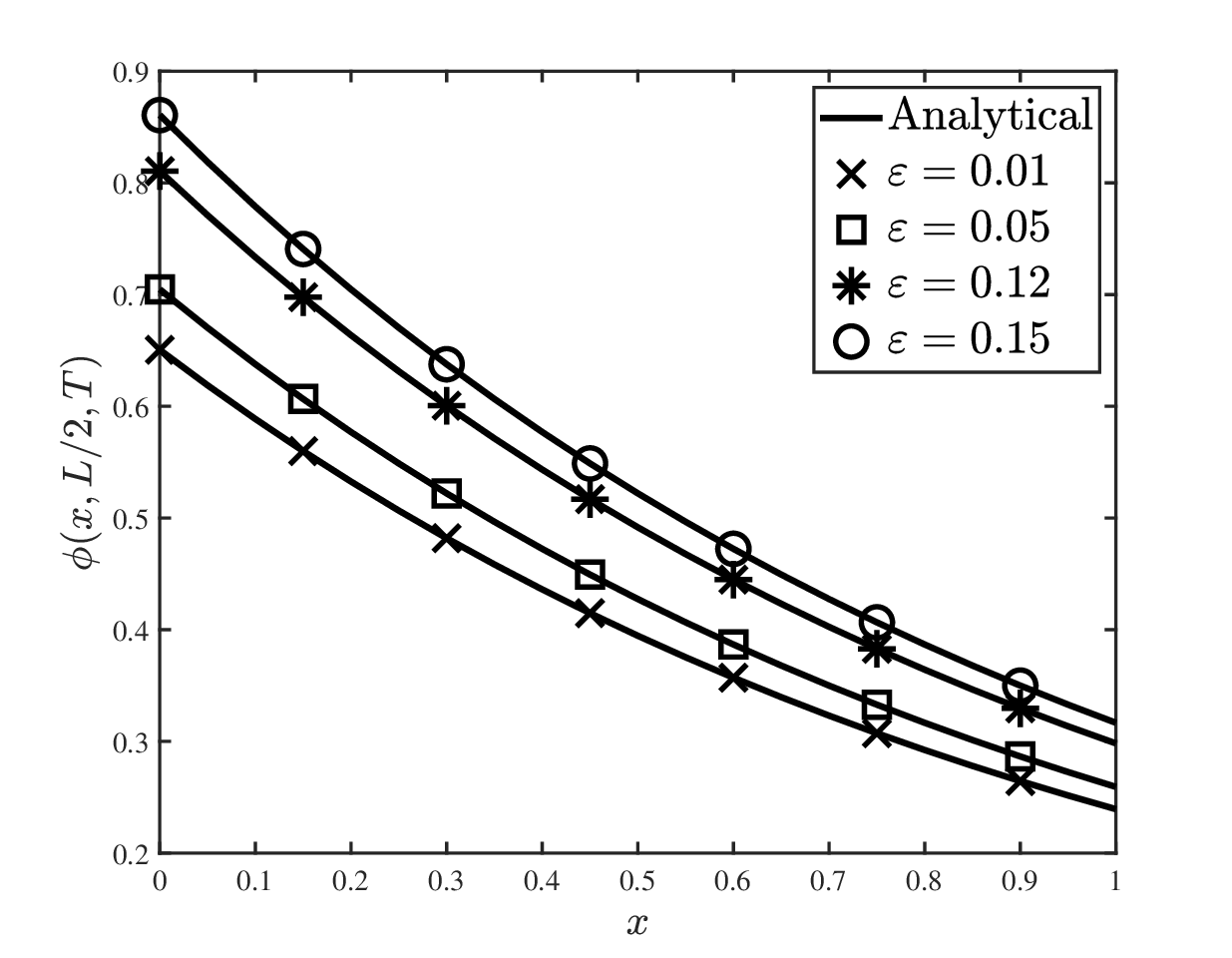}
		\end{minipage}
	}
	\caption{The numerical and analytical solutions under different values of discretization parameter $\varepsilon$ ($T =4.0$).}
	\label{Ex2Figure1}
	\setlength{\belowcaptionskip}{-1cm}
\end{figure}

\begin{figure}
	\centering
	\subfigure[$\:$MRT-LB model]
	{
		\begin{minipage}{0.40\linewidth}
			\centering
			\includegraphics[width=2.5in]{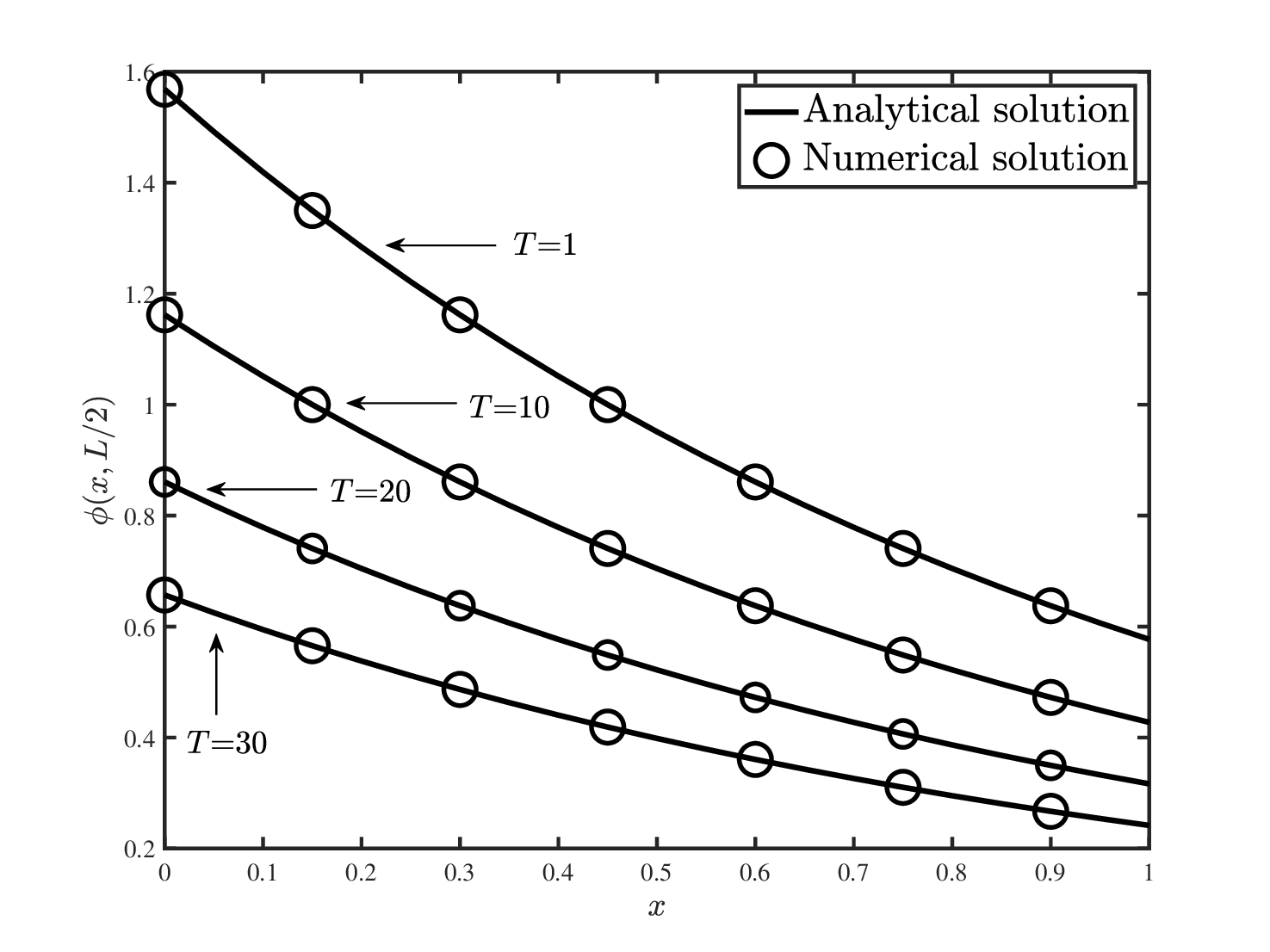}
		\end{minipage}
	}
	\subfigure[$\:$FLFD scheme]
	{
		\begin{minipage}{0.40\linewidth}
			\centering
			\includegraphics[width=2.5in]{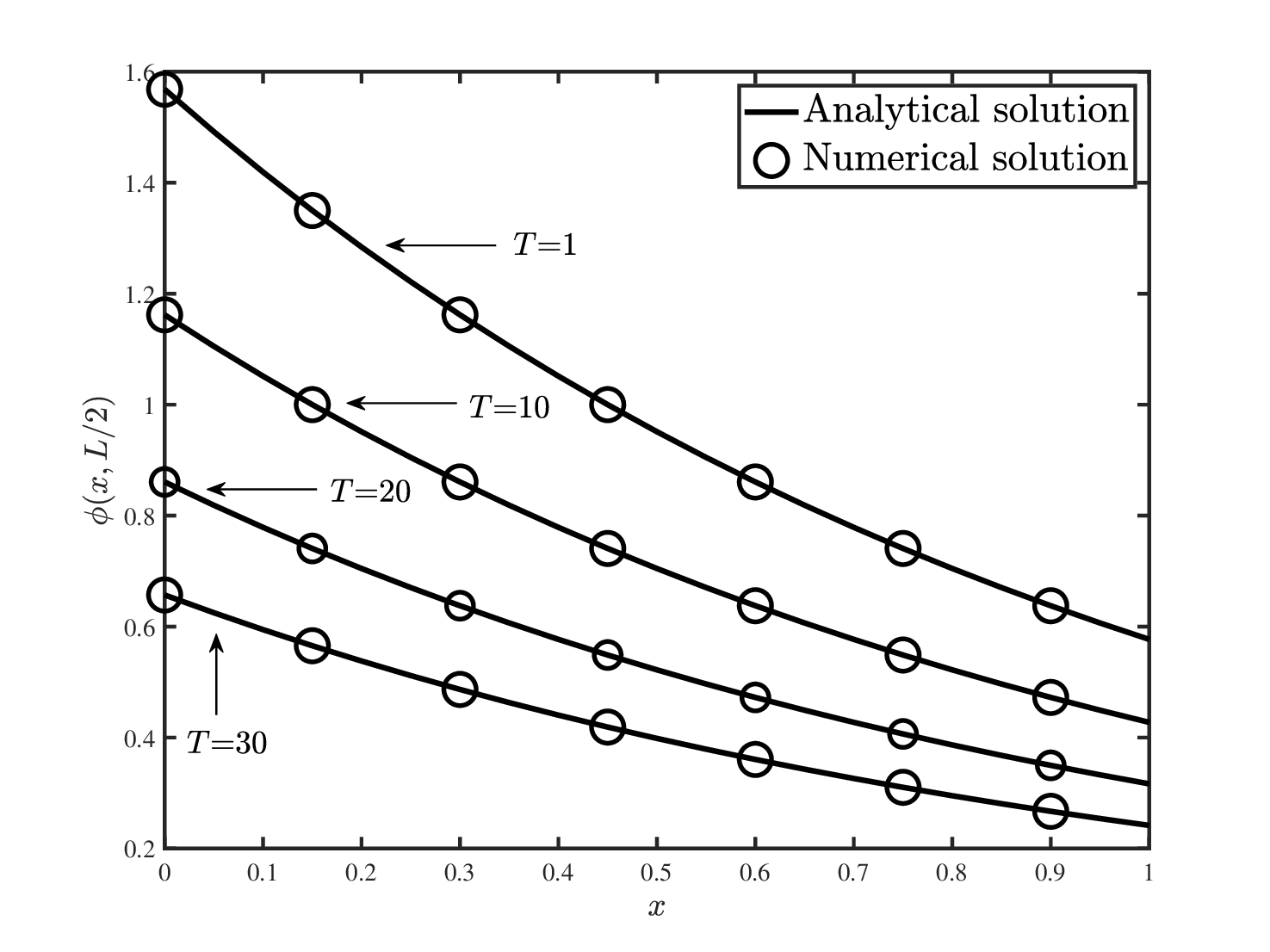}
		\end{minipage}
	}
	\caption{The numerical and analytical solutions under different values of time $T$ ($\varepsilon=0.15$).}
	\label{Ex2Figure2}
	\setlength{\belowcaptionskip}{-1cm}
\end{figure}

\begin{figure}
	\centering
	\subfigure[$\:$Analytical solution]
	{
		\begin{minipage}{0.40\linewidth}
			\centering
			\includegraphics[width=2.5in]{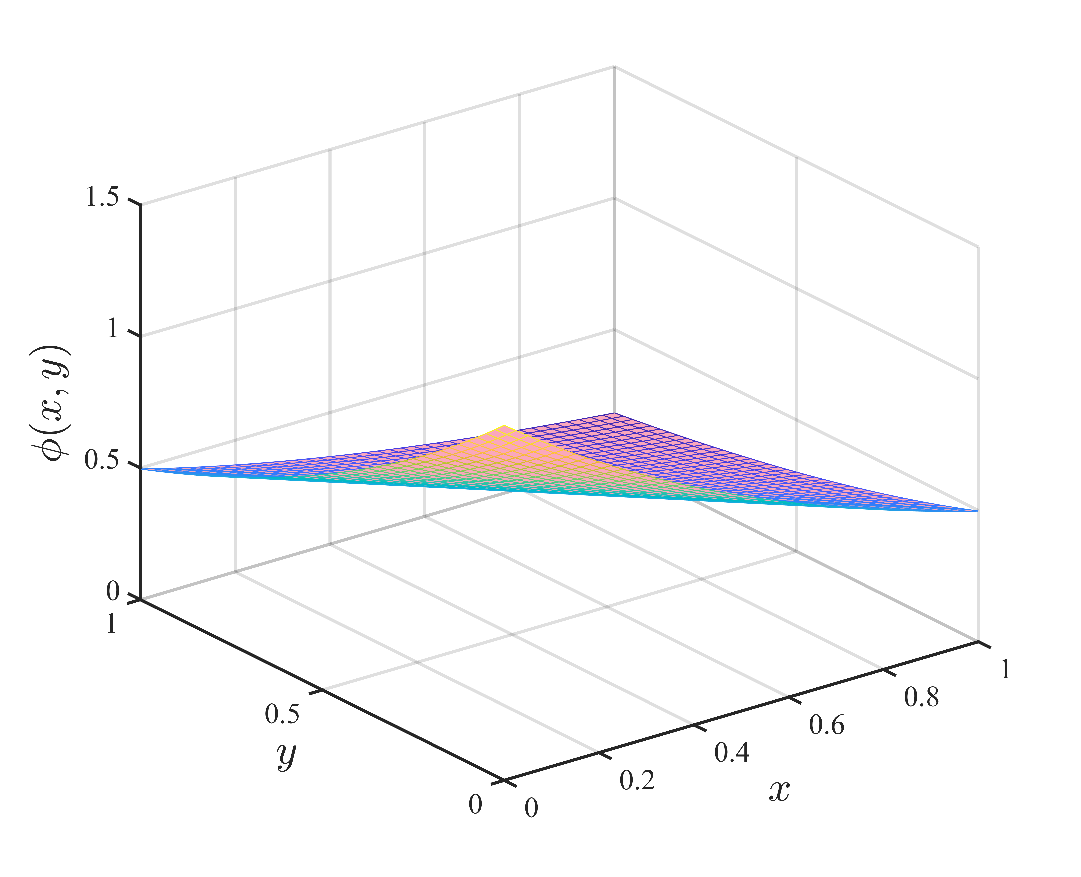}
		\end{minipage}
	}

	\subfigure[$\:$MRT-LB model]
	{
		\begin{minipage}{0.40\linewidth}
			\centering
			\includegraphics[width=2.5in]{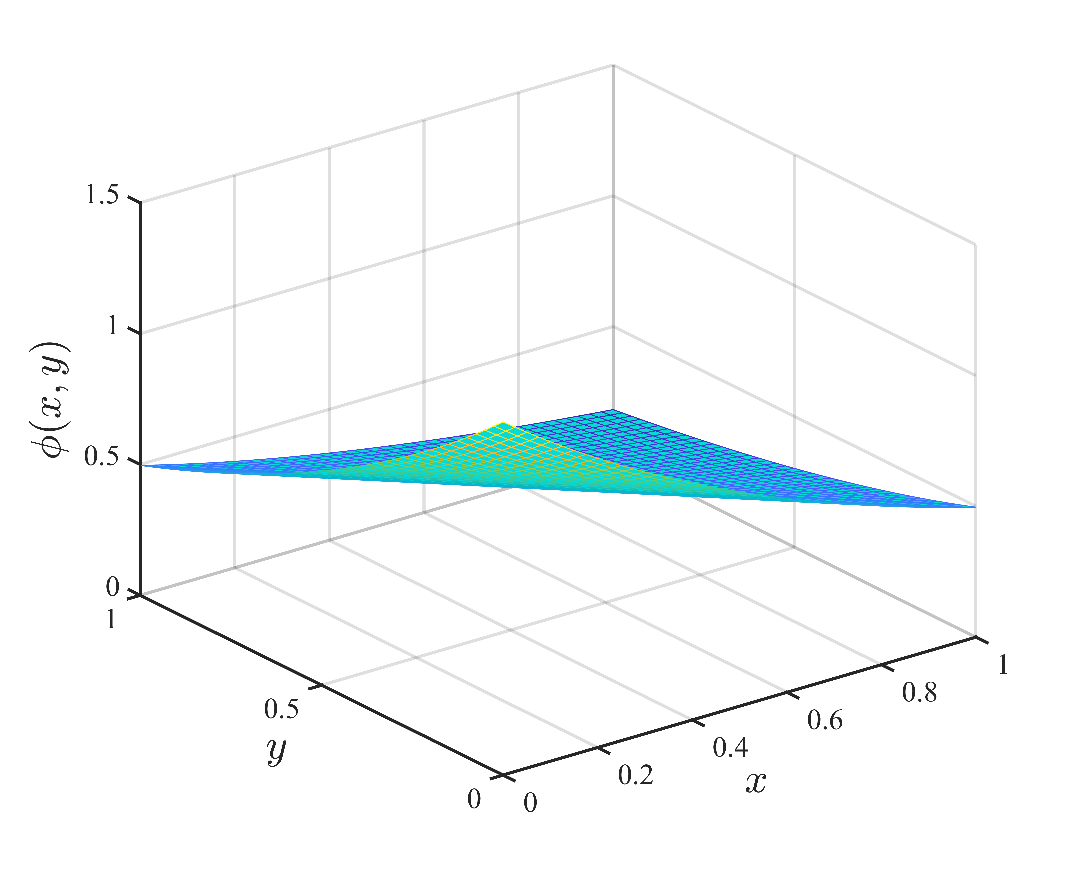}
		\end{minipage}
	}
	\subfigure[$\:$FLFD scheme]
	{
		\begin{minipage}{0.40\linewidth}
			\centering
			\includegraphics[width=2.5in]{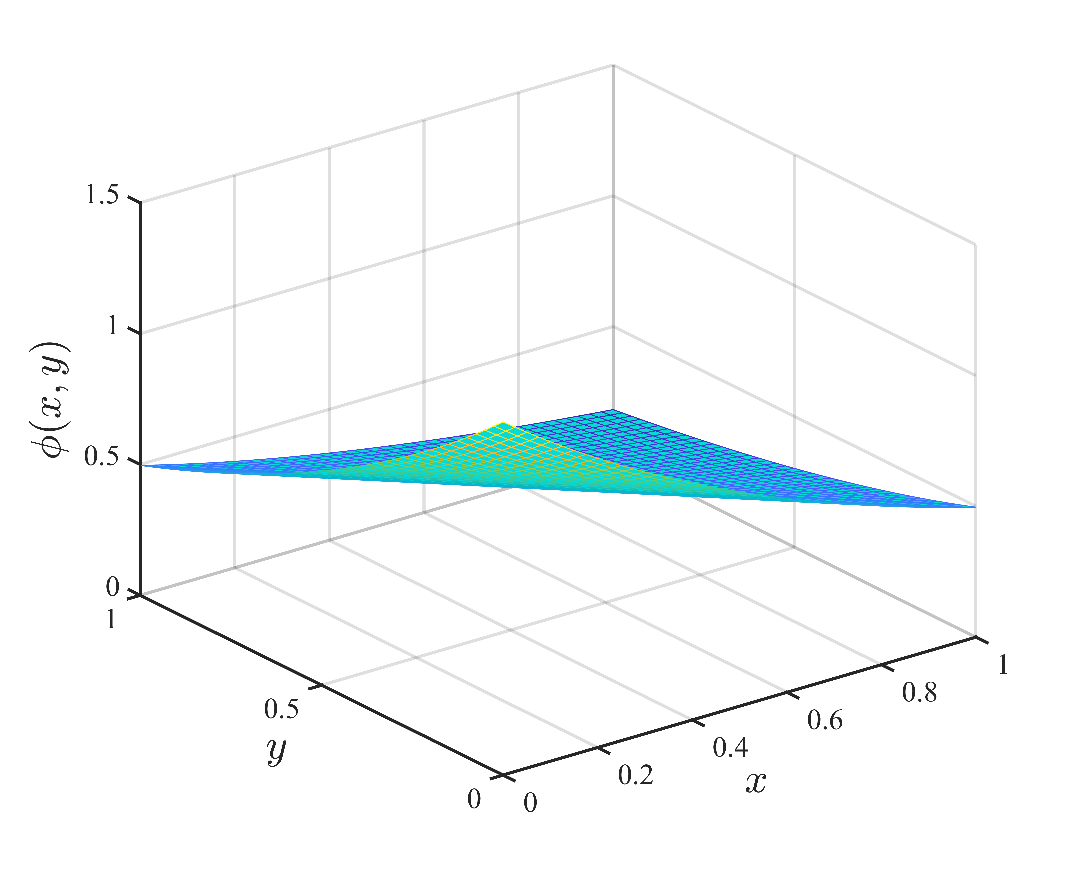}
		\end{minipage}
	}
	\caption{Distributions of the variable $\phi$ at the discretization parameter $\varepsilon=0.15$ and time $T=1.0$.}
	\label{Ex2Figure3}
	\setlength{\belowcaptionskip}{-1cm}
\end{figure}

\begin{figure}
	\centering
	\subfigure[$\:$MRT-LB model]
	{
		\begin{minipage}{0.40\linewidth}
			\centering
			\includegraphics[width=2.5in]{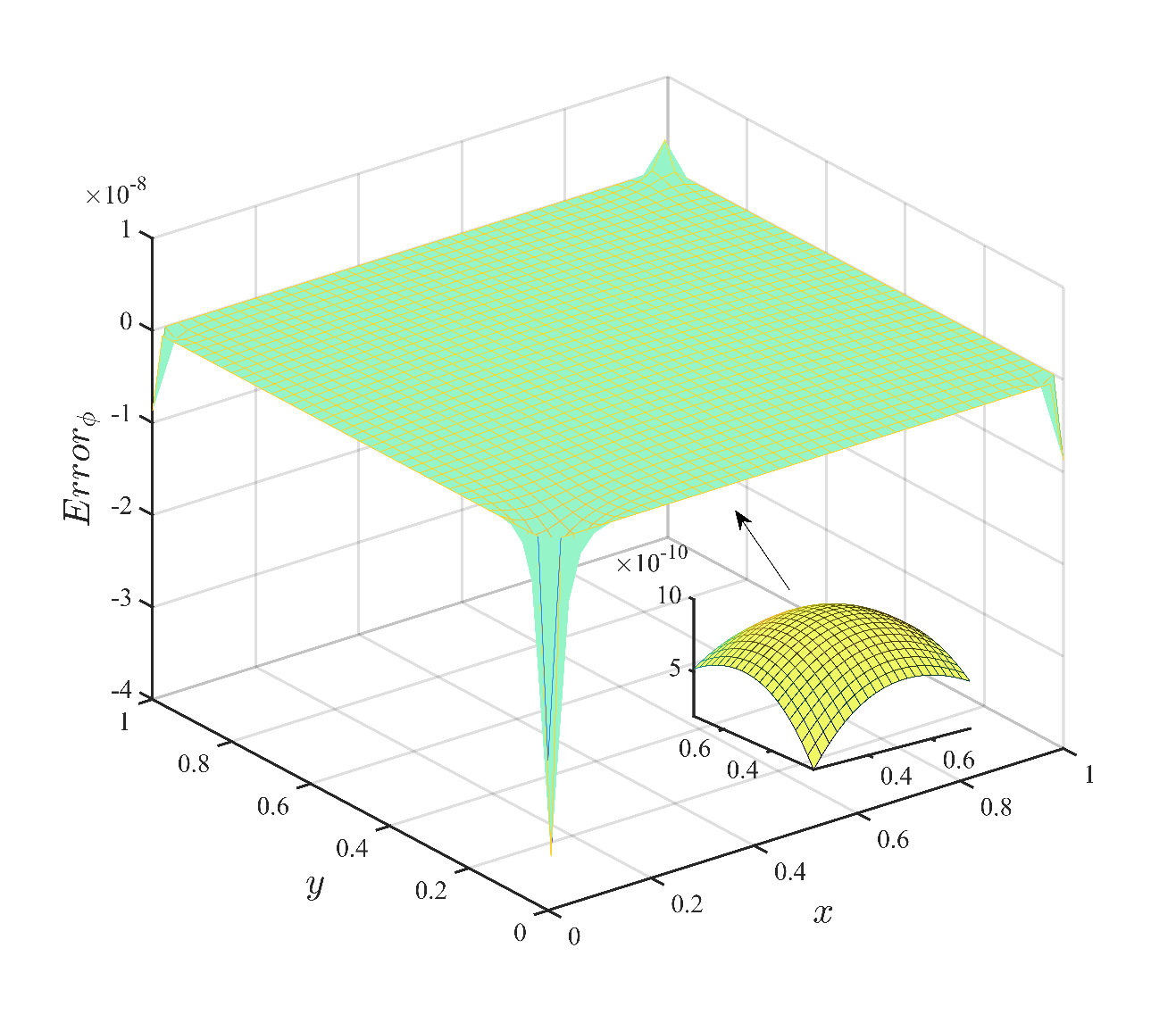}
		\end{minipage}
	}
	\subfigure[$\:$FLFD scheme]
	{
		\begin{minipage}{0.40\linewidth}
			\centering
			\includegraphics[width=2.5in]{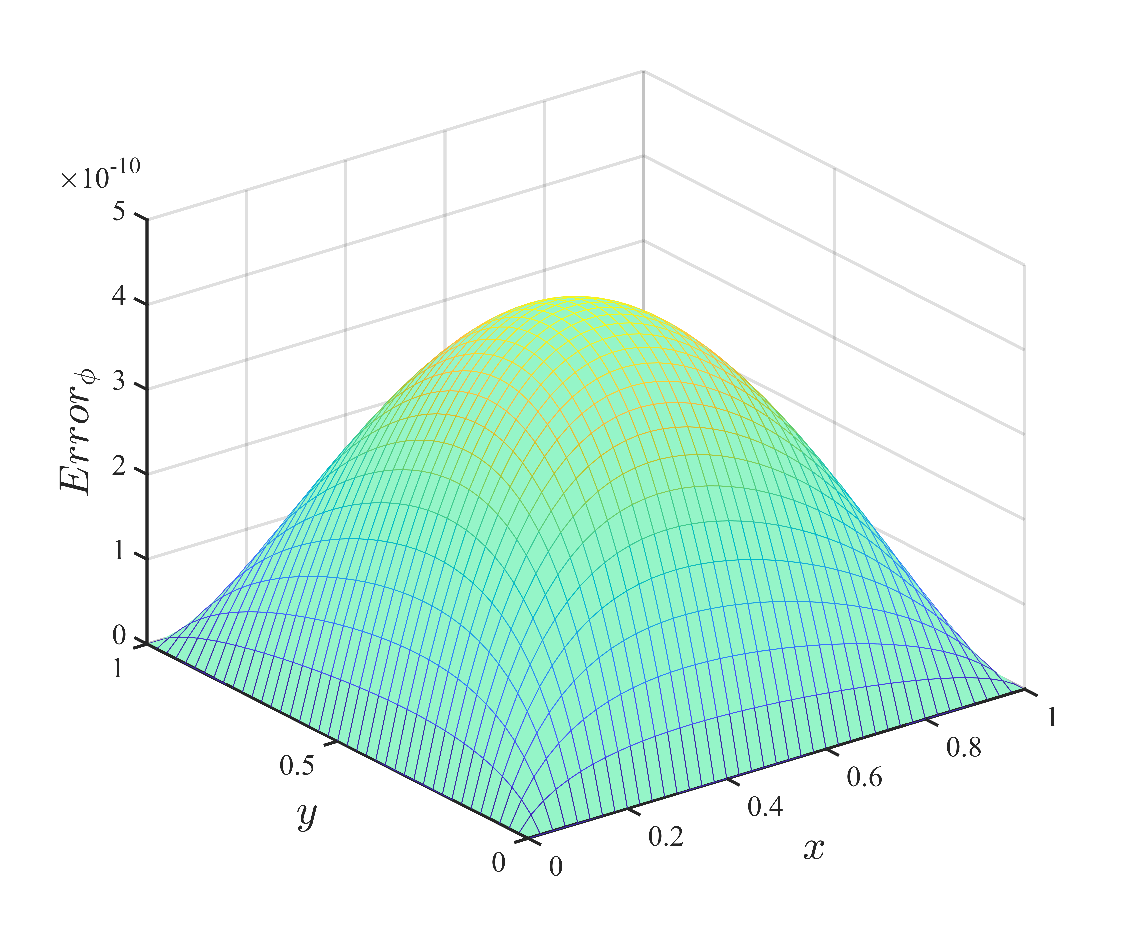}
		\end{minipage}
	}
	\caption{Distributions of the error (${\rm Error}_{\phi}$) at the discretization parameter $\varepsilon=0.15$ and time $T=1.0$.}
	\label{Ex2Figure4}
	\setlength{\belowcaptionskip}{-1cm}
\end{figure}

Similar to the previous example, we set the discretization parameter $\varepsilon$, which is corresponding to a specified different diffusion coefficient $\kappa$ for the given lattice spacing $\Delta x=1/20$ and time step $\Delta t=1/40$. We first conduct some simulations at different values of discretization parameter $\varepsilon$ ($\varepsilon =0.01,0.05,0.12,0.15)$ and $T=4.0$, and also at different values of time $T$ ($T=1,10,20,30$) and $\varepsilon = 0.15$. As shown in Figs. \ref{Ex2Figure1} and \ref{Ex2Figure2}, the numerical results are very close to the analytical solutions. Furthermore, to give a comparison between the analytical solution and numerical results of the MRT-LB model and FLFD scheme, we perform the simulations at the discretization parameter $\varepsilon=0.15$ and time $T=1.0$, and plot the results in Figs.  \ref{Ex2Figure3} and \ref{Ex2Figure4}. It is clear that the numerical results are in good agreement with the analytical solution (see Fig. \ref{Ex2Figure3}), while the
MRT-LB model gives the larger errors at four corner points [see Fig. \ref{Ex2Figure4}(a)]. This is because the treatment on the Dirichelet boundary conditions at four corner points in the MRT-LB model is not accurate enough. On the contrary, the FLFD scheme is more accurate, especially at the four corner points [see Fig. \ref{Ex2Figure4}(b)], which is caused by the fact that the values of $\phi$ at the corner points are given precisely by the Dirichelet boundary conditions.

\begin{figure}
	\centering
	\subfigure[$\:$MRT-LB model]
	{
		\begin{minipage}{0.40\linewidth}
			\centering
			\includegraphics[width=2.5in]{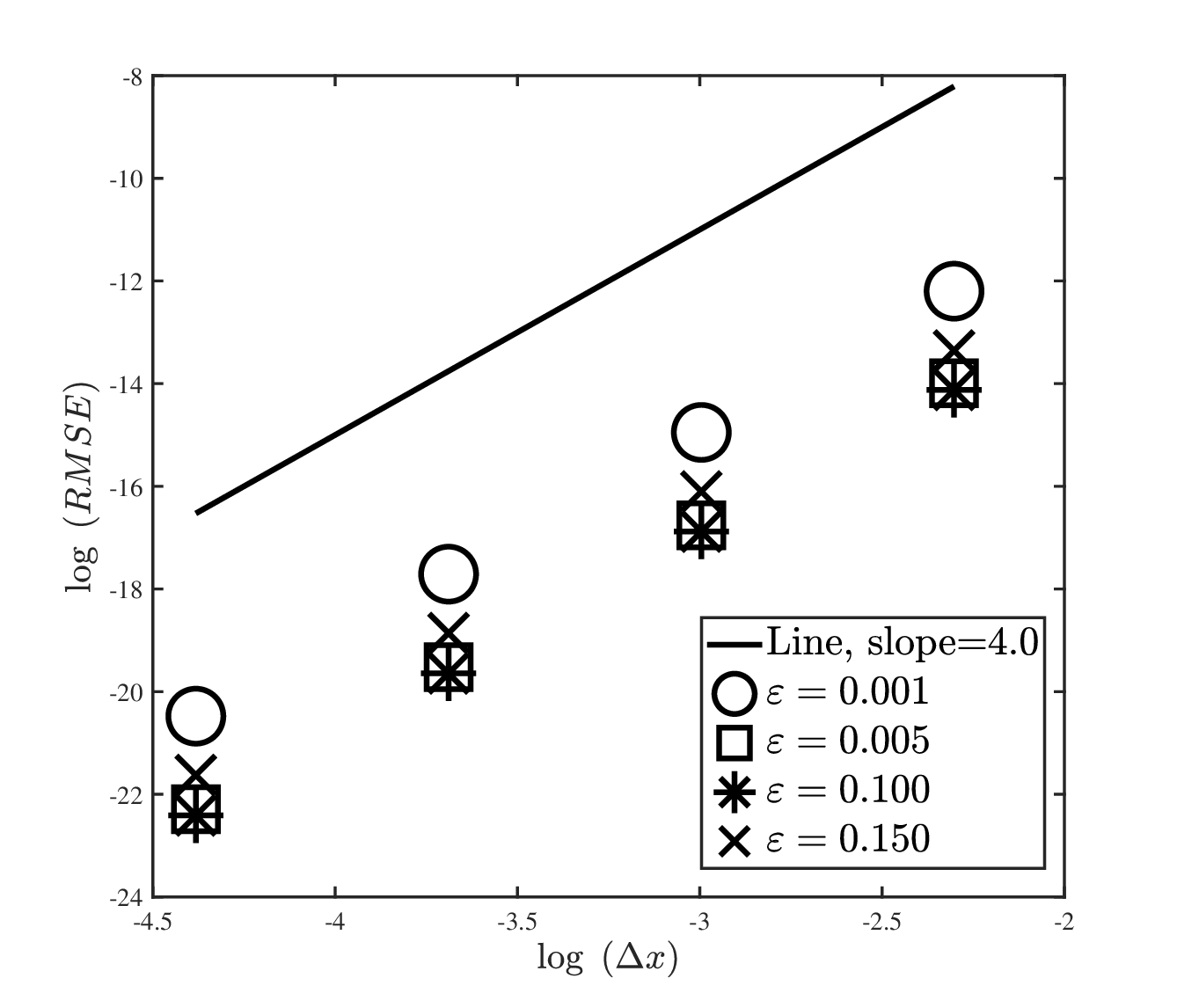}
		\end{minipage}
	}
	\subfigure[$\:$FLFD scheme]
	{
		\begin{minipage}{0.40\linewidth}
			\centering
			\includegraphics[width=2.5in]{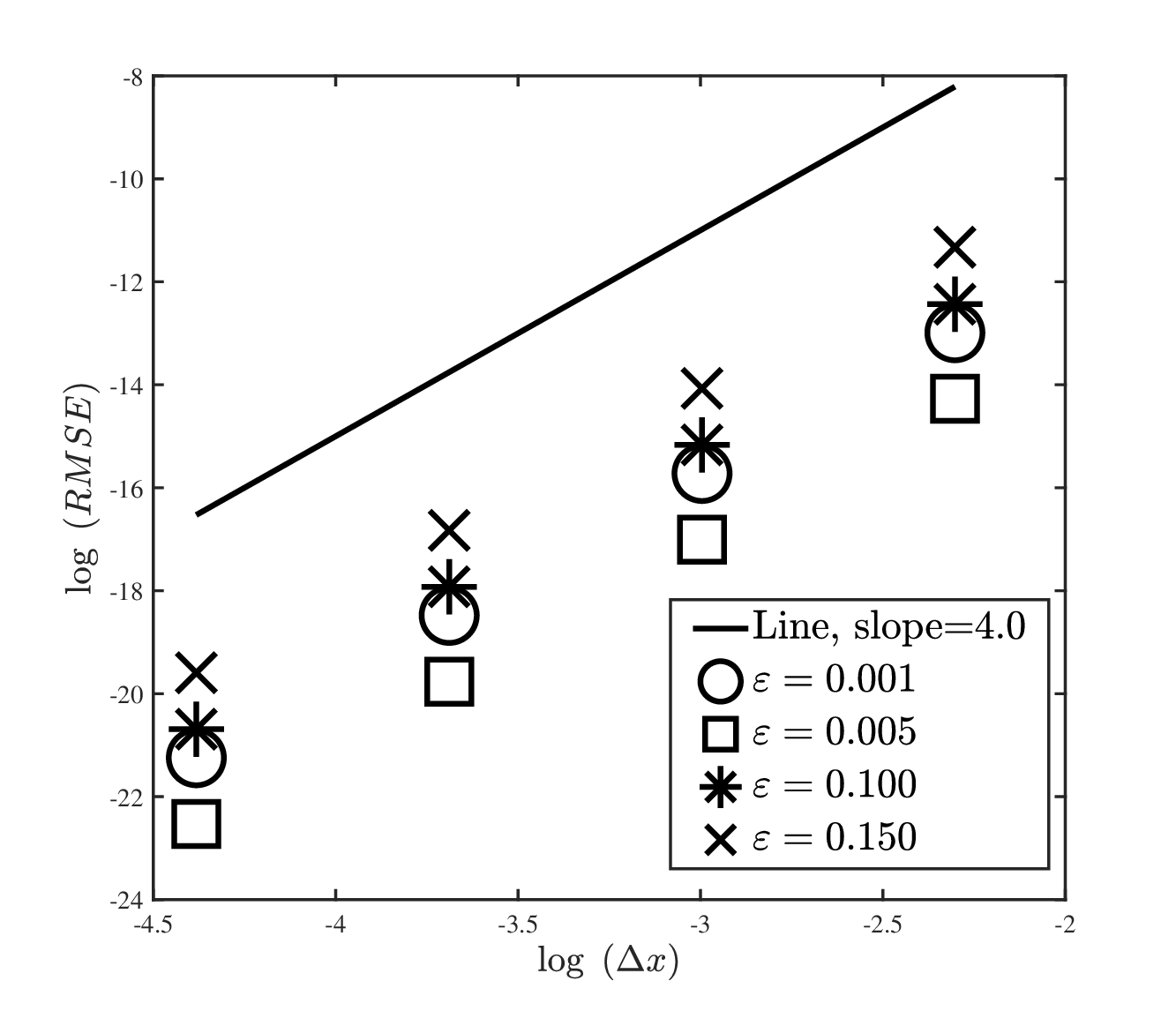}
		\end{minipage}
	}
	\caption{The convergence rates of MRT-LB model and FLFD scheme with different values of discretization parameter $\varepsilon$ ($T=10.0$).}
	\label{Ex2Slope1}
\end{figure}

In addition, this problem is also adopted to test the convergence rates of the present MRT-LB model and FLFD scheme. For this purpose, we perform a number of simulations with some specified values of discretization parameter $\varepsilon$, and calculate the $RMSEs$ and $CRs$ under different lattice sizes. As seen from Fig. \ref{Ex2Slope1}, Tables \ref {TEx2CRLB} and \ref{TEx2CRFD} where the simulations are suspended at time $T=10.0$ and the lattice spacing is varied from $\Delta x/8$ to $\Delta x=1/10$ with a fixed $\Delta x^2/\Delta t$, the present MRT-LB model and FLFD scheme have a fourth-order convergence rate in space.

\begin{table}
	\begin{center}
		\caption{The $RMSEs$ and $CRs$ of MRT-LB model under different values of discretization parameter $\varepsilon$ ($T=10.0$, $\Delta x=1/10$, $\Delta t=1/10$).}
		\label{TEx2CRLB}  
		\begin{tabular}{ccccccccc}
			\hline \hline
			$\varepsilon$&$RMSE_{\Delta x,\Delta t}$&$RMSE_{\Delta x/2,\Delta t/4}$&
			$RMSE_{\Delta x/4,\Delta t/16}$&
			$RMSE_{\Delta x/8,\Delta t/64}$&$CR$\\
			\midrule[1pt]
			0.04 &  $8,8311\times 10^{-8}$ &$6.6863\times 10^{-9}$ &$4.3672\times 10^{-10}$&$2.9422\times 10^{-11}$&$\sim  3.5475$ \\
			
			0.08&$1.3741\times 10^{-7}$&$7.2973\times 10^{-9}$&$4.6827\times 10^{-10}$&$2.9422\times 10^{-11}$&$\sim 3.7263$\\
			
			0.10&$2.1585\times 10^{-7}$&$9.3445\times 10^{-9}$&$6.1023\times 10^{-10}$&$3.9040\times 10^{-11}$&$\sim 4.1420$\\
			
			0.12&$3.1187\times 10^{-7}$&$1.6277\times 10^{-8}$&$1.0587\times 10^{-9}$&$6.7528\times 10^{-11}$&$\sim 4.0577$\\
			\hline \hline
		\end{tabular}
	\end{center}
\end{table}
\begin{table}
	\begin{center}
		\caption{The $RMSEs$ and $CRs$ of FLFD scheme under different values of discretization parameter $\varepsilon$ ($T=10.0$, $\Delta x=1/10$, $\Delta t=1/10$).}
		\label{TEx2CRFD}  
		\begin{tabular}{ccccccccc}
			\hline \hline
			$\varepsilon$&$RMSE_{\Delta x,\Delta t}$&$RMSE_{\Delta x/2,\Delta t/4}$&
			$RMSE_{\Delta x/4,\Delta t/16}$&
			$RMSE_{\Delta x/8,\Delta t/64}$&$CR$\\
			\midrule[1pt]
			0.04 &  $2.5276\times 10^{-10}$ &$5.8786\times 10^{-11}$ &$1.2537\times 10^{-12}$&$6.3751\times 10^{-14}$&$\sim 3.978$ \\
			
			0.08&$8.9627\times 10^{-10}$&$6.6733\times 10^{-11}$&$4.3916\times 10^{-12}$&$2.6318\times 10^{-13}$&$\sim 3.9689$\\
			
			0.10&$2.8156\times 10^{-9}$&$2.0894\times 10^{-10}$&$1.3743\times 10^{-11}$&$8.3952\times 10^{-13}$&$\sim 3.9039$\\
			
			0.12&$6.6352\times 10^{-9}$&$4.9170\times 10^{-10}$&$3.2335\times 10^{-11}$&$2.0162\times 10^{-12}$&$\sim 3.8948$\\
			\hline \hline
		\end{tabular}
	\end{center}
\end{table}

\section{Conclusions}\label{sec6}

In this paper, we first developed the general MRT-LB model for the two-dimensional diffusion equation with a constant source term, and then obtained an equivalent macroscopic SLFD scheme (\ref{fd}) and a simplified FLFD scheme (\ref{fds5=1}). It is found that both the MRT-LB model and SLFD scheme (\ref{fd}) can achieve a fourth-order accuracy in space once the weight coefficient $w_0$, the relaxation parameters $s_2$, $s_4$ and $s_5$ satisfy Eq. (\ref{solution}). In addition, through a detailed theoretical analysis, we also proved that the MRT-LB model and FLFD scheme (\ref{fds5=1}), are unconditionally stable. We further conducted some simulations, and the results show that the MRT-LB model and the equivalent finite-difference scheme indeed have a fourth-order convergence rate in space, which is consistent with our theoretical analysis. It is worth mentioning that the MRT-LB model may be more efficient than its equivalent finite-difference scheme from a computational point of view since it is only a two-level method. Finally, it should be noted that we can also develop a fourth-order MRT-LB model for the two-dimensional diffusion equation with a linear source term, and obtain the equivalent macroscopic six-level finite-difference scheme. 

\section*{Acknowledgments}

The computation is completed in the HPC Platform of Huazhong University of Science and Technology. This work was financially supported by the National Natural Science Foundation of China (Grants No. 12072127 and No. 51836003).
\appendix
\section{The discussion on the transform matrix $M$}\label{AM}
Actually, in the MRT-LB model with D2Q5 lattice structure, the transform matrix $\bm{M}$, unlike that in Eq. (\ref{ms}), can also be given in a natural way,
\begin{equation}\label{m2}
	\begin{aligned}
		&\bm{M_N}=\left (
		\begin{matrix}
			1& 1 &1&1&1\\
			0 & c & 0&-c &0\\
			0 & 0&c & 0&-c\\
			0 & c^2 & 0&c^2&0\\
			0 & 0 &c^2&0 &c^2\\
		\end{matrix}
		\right ),
	\end{aligned}
\end{equation}
where the non-orthogonal vectors are used. Following the approach in Sec. \ref{sec3}, one can derive the following equivalent macroscopic SLFD scheme from the MRT-LB model (\ref{distributionf}),
\begin{equation}\label{fd3}
	\begin{aligned}
		\phi_{i,j}^{n+1}&=\alpha_1\phi_{i,j}^{n}+\alpha_2(\phi_{i-1,j}^{n}+\phi_{i+1,j}^{n})+\alpha_3(\phi_{i,j-1}^{n}+\phi_{i,j+1}^{n})+\beta_1\phi_{i,j}^{n-1}+\beta_2(\phi_{i-1,j}^{n-1}+\phi_{i+1,j}^{n-1})+\beta_3(\phi_{i,j-1}^{n-1}+\phi_{i,j+1}^{n-1})\\
		&\quad+\beta_4(\phi_{i-1,j-1}^{n-1}+\phi_{i-1,j+1}^{n-1}+\phi_{i+1,j-1}^{n-1}+\phi_{i+1,j+1}^{n-1})+\gamma_1\phi_{i,j}^{n-2}+\gamma_2(\phi_{i-1,j}^{n-2}+\phi_{i+1,j}^{n-2})+\gamma_3(\phi_{i,j-1}^{n-2}+\phi_{i,j+1}^{n-2})\\
		&\quad+\gamma_4(\phi_{i-1,j-1}^{n-2}+\phi_{i-1,j+1}^{n-2}+\phi_{i+1,j-1}^{n-2}+\phi_{i+1,j+1}^{n-2})+\zeta_1\phi_{i,j}^{n-3}
		+\zeta_2(\phi_{i-1,j}^{n-3}+\phi_{i+1,j}^{n-3})+\zeta_3(\phi_{i,j-1}^{n-3}+\phi_{i,j+1}^{n-3})\\
		&\quad+\eta \phi_{i,j}^{n-4}
		+\Delta t\delta R,
	\end{aligned}
\end{equation}
where the parameters $\alpha_k,\zeta_k\:(k=1,2,3)$, $\beta_l,\gamma_l\:(l=1,2,3,4),\eta$ and $\delta$ are given by
\begin{equation}\label{parameter3}
	\begin{aligned}
		&\alpha_1=1-2w_1(s_4+s_5),\alpha_2=-\frac{s_2+s_4-2}{2}+w_1s_4,\alpha_3=-\frac{s_2+s_5-2}{2}+w_1s_5,\\
		&\beta_1=(s_4+s_5-2)(1-s_2)-2w_1(1-s_2)(s_4+s_5),\\
		&\beta_2=\frac{s_2+s_4-2}{2}+s_4w_1(1-s_2)-s_5w_1(s_2+s_4-2),\\
		&\beta_3=\frac{s_2+s_5-2}{2}+s_5w_1(1-s_2)-s_4w_1(s_2+s_5-2),\\
		&\beta_4=-\frac{(s_2+s_4-2)(s_2+s_5-2)}{4}+\frac{s_4w_1(s_2+s_5-2)+s_5w_1(s_2+s_4-2)}{2},\\
		&\gamma_1=(s_2-1)(s_4+s_5-2)-2w_1(s_2-1)(2s_4s_5-s_4-s_5),\\
		&\gamma_2=-\frac{(s_2-1)(s_5-1)(s_2+s_4-2)}{2}+w_1(s_2-1)[s_4(s_5-1)+s_5(s_2+s_4-2)],\\
		&\gamma_3=-\frac{(s_2-1)(s_4-1)(s_2+s_5-2)}{2}+w_1(s_2-1)[s_5(s_4-1)+s_4(s_2+s_5-2)],\\
		&\gamma_4=\frac{(s_2+s_4-2)(s_2+s_5-2)}{4}+w_1(1-s_2)\big[\frac{s_4(s_2+s_5-2)}{2}+\frac{s_5(s_2+s_4-2)}{2}\big],\\
		&\zeta_1=-(s_2-1)^2(s_4-1)(s_5-1)+2w_1(s_2-1)^2(2s_4s_5-s_4-s_5),\\
		&\zeta_2=\frac{(s_2-1)(s_5-1)(s_2+s_4-2)}{2}-w_1s_4(s_2-1)^2(s_5-1),\\
		&\zeta_3=\frac{(s_2-1)(s_4-1)(s_2+s_5-2)}{2}-w_1s_5(s_2-1)^2(s_4-1),\\
		&\eta=(s_2-1)^2(s_4-1)(s_5-1),\delta=s_2^2s_4s_5.
	\end{aligned}
\end{equation}
From above Eqs. (\ref{fd3}) and (\ref{parameter3}), one can find that the parameters $\alpha_2,\beta_2,\gamma_2$ and $\zeta_2$ before the difference terms in $x$ direction are same to the parameters $\alpha_3,\beta_3,\gamma_3$ and $\zeta_3$ before the difference terms in $y$ directions if and only if $s_4=s_5$.

Through the Taylor expansion anlysis, one can obtain the following conditions to ensure that the SLFD scheme has a fourth-order accuracy at the diffusive scaling ($\Delta t \propto \Delta x^2$),
\begin{align}\label{solutionm}
	\left\{\begin{aligned}
		&w_0= 1-4\sqrt{3} \varepsilon,\\
		&s_2=\frac{6}{3+\sqrt{3}},\\
		&s_4= s_5=\frac{6}{2\sqrt{3}+3},\\
	\end{aligned}\right.
\end{align} 
From Eq. (\ref{solutionm}) one can find that the MRT-LB model with the transform matrix $\bm{M_N}$ has a fourth-order accuracy in space, and it would reduce to the TRT-LB model \cite{Gin2012} if the relaxation parameters corresponding to the odd and even-order moments  are given by $s_o=6/(3+\sqrt{3})$ and $s_e=6/(2\sqrt{3}+3)$. 

	We now give a comparison of two different transform matrices.
\begin{itemize}
	\item [(i)] The choice of transform matrix $\bm{M}$ determines that the relaxation parameters $s_4$ and $s_5$ are corresponding to two linearly independent quantities, which are the linear combination of the zero-order moment, the second-order moments in the $x$ and $y$ directions. However, the choice of the transform matrix $\bm{M_N}$ determines that the relaxation parameter $s_4$ corresponding to the second-order moment in the $x$ direction, while the relaxation parameter $s_5$ corresponding to the second-order moment in the $y$ direction.
	\item [(ii)] In the accuray analysis of the MRT-LB model with the transform matrix $\bm{M}$, the coefficients before the fourth-order spatial derivative terms $\partial^4_x\phi$ and $\partial^4_y\phi$ are same. In this case, the number of the second-order truncation-error terms, i.e., $O(\Delta x^2+\Delta t+\Delta x^4/\Delta t)$, to be eliminated is two and the relation (\ref{errorfirst}) is satisfied, thus the number of the free parameters $w_0,s_2,s_4$ and $s_5$ is one more than that of the conditions to ensure that the MRT-LB model has a fourth-order accuracy. However, when the transform matrix $\bm{M_N}$ is considered, the coefficients before the fourth-order spatial derivative terms $\partial^4_x\phi$ and $\partial^4_y\phi$ are not same unless $s_4=s_5$, which indicates the number of the free parameters $w_1,s_2,s_4$ and $s_5$ is no more than that of the fourth-order conditions.
\end{itemize}
From above discussion, one can see that the relaxation parameters $s_2,s_4$ and $s_5$ in Eq. (\ref{solution}) can be changed more flexibly, while those in Eq. (\ref{solutionm}) cannot, which shows that the time level of the equivalent fourth-order finite-difference scheme derived by the MRT-LB model with transform matrix $\bm{M_N}$ and diagonal relaxation matrix $\bm{S}$ can only be six. However, if we consider the relaxation matrix as $\bm{S_N}=\bm{M_N}\bm{M}^{-1}\bm{SM}\bm{M_N}^{-1}$, i.e.,
\begin{equation}\label{m2}
	\begin{aligned}
		&\bm{S_N}=\left (
		\begin{matrix}
			s_1& 0 &0&0&0\\
			0 & s_2 & 0&0 &0\\
			0 & 0&s_3 & 0&0\\
			\frac{2c^2(s_1-s_4)}{5} & 0 & 0&\frac{s_4+s_5}{2}&\frac{s_4-s_5}{2}\\
			\frac{2c^2(s_1-s_4)}{5} & 0 &0&\frac{s_4-s_5}{2} &\frac{s_4+s_5}{2}\\
		\end{matrix}
		\right ),
	\end{aligned}
\end{equation}
it is obvious and easy to show that the MRT-LB model with transfrom matrix $\bm{M_N}$ and relaxation martix $\bm{S_N}$ is same as that with transfrom matrix $\bm{M}$ and diagonal relaxation martix $\bm{S}$ \cite{Chai2020}, thus the equivalent finite-differece schemes of them are same as each other. This indicates that both the transfrom matrix and relaxation matrix have some important influences on the accuracy analysis of the MRT-LB model. 

	\section{The coefficients of characteristic polynomial $\mathcal{X}_{P}$}\label{A1}
For brevity, we first introduce $Z_1:=\sum_{i=1}^4T_i/4$ and $Z_2:=(T_1T_2+T_1T_4+T_2T_3+T_2T_4)/4$, then through the Faddeev-Leverrier algorithm for the computation of the characteristic polynomial of a square matrix \cite{Hou}, one can express the coefficients of polynomial $\mathcal{X}_{\bm{P}}$ as
\begin{subequations}\label{xishu}
	\begin{align}
		\gamma_4&=\frac{4s_5}{5}-1+(2s_2+\frac{s_4}{5}+s_5-4)Z_1,\\
		\gamma_3&=
		\big(\frac{4s_4s_5+8s_2s_4-17s_4}{5}-s_5-2s_2+4
		\big)Z_1\notag\\
		&\quad+
		\big(\frac{s_4s_5+s_2s_4-2s_4}{5}+s_2s_5-2s_5+s_2^2-4s_2
		\big)Z_2+(s_5-2+\frac{s_4}{5})(s_2-1),\\
		\gamma_2&=(\frac{4s_2s_4s_5-2s_4s_5+s_2^2s_4-4s_2s_4+3s_4}{5}+s_2^2s_5-4s_2s_5+3s_5-2s_2^2+6s_2-4)Z_1\notag\\
		&\quad+(\frac{4s_2s_4s_5-9s_4s_5+4s_2^2s_5-17s_2s_4+18s_4}{5}-s_2s_5+2s_5-s_2^2+4s_2-4)Z_2\notag\\
		&\quad+
		\big(\frac{4s_2s_4s_5-4s_4s_5-9s_2s_4+9s_4}{5}-s_2s_5+s_5+2s_2-2\big),\\
		\gamma_1&=
		\big(\frac{4s_2s_4s_5(s_2-1)}{5}-56s_4s_5-36s_2s_4+76s_4-20s_2s_5+60s_5+40s_2-80\big)Z_1\notag\\
		&\quad+\big(\frac{s_2s_4s_5(s_2-1)}{5}-4s_4s_5-4s_2s_4+4s_4-20s_2s_5+20s_5+20s_2-40\big),\\
		\gamma_0&=-(1-s_2)^2(1-s_4)(1-s_5)T_0.
	\end{align}
\end{subequations}

	\section{The six-level finite-difference scheme for the diffusion equations with a linear source term}\label{A5}
The two-dimensional diffusion equation with a linear source term can be written as
\begin{equation}\label{equlinear}
	\partial_t\phi=\nabla\cdot(\kappa\nabla\phi)+\zeta\phi+R,
\end{equation}
where the parameters $\zeta$ and $R$ are two constants.

The evolution equation of the MRT-LB model for the Eq. (\ref{equlinear}) reads \cite{Chai2014, Cui2016}
\begin{equation}\label{distributionf2}
	\begin{aligned}
		f_i(\bm{r}+c_i\Delta t,t+\Delta t)&=f_i(\bm{r},t)-\big(\bm{M}^{-1}\bm{SM}\big)_{i,k}[f_k(\bm{r},t)-f_k^{eq}(\bm{r},t)]+\Delta t (\bm{I}-\frac{\bm{M}^{-1}\bm{SM}}{2})_{i,k}\tilde{R}_k.
	\end{aligned}
\end{equation}
Which is named as the MRT-LB-AL model hereafter. For the D2Q5 lattice structure, the discrete velocity $\bm{c_{i}}$, the matrices $\bm{M}$ and $\bm{S}$, the equilbrium distribution function $f_i^{eq}$ are same as Eqs. (\ref{velocity}) (\ref{ms}) and (\ref{feq}), while the discrete source term $\tilde{R}_i$ is given by
\begin{equation}
	\tilde{R}_i=w_i(\zeta \phi(\bm{r},t)+R),
\end{equation}
and simultaneously, the conservative moment $\phi(\bm{r},t)$ is calculated by
\begin{equation}\label{phi_fi}
	\phi(\bm{r},t)=\sum_if_i^{eq}
	=\frac{\sum_if_i+\Delta t R/2}{1-\zeta\Delta t/2}.
\end{equation}
Through the direct Taylor expansion analysis \cite{Holdych2004, Wagner2006, Dubois2008}, one can correctly recover the equation (\ref{equlinear}) from the MRT-LB-AL model (\ref{distributionf2}) with the relation (\ref{kappa}). 

Similarly, one can derive the following equivalent macroscopic six-level finite-difference scheme (SLFD-AL scheme) from the MRT-LB-AL model (\ref{distributionf2}),
\begin{equation}\label{fdal}
	\begin{aligned}
		\phi_{i,j}^{n+1}&=\overline{\alpha}_1\phi_{i,j}^{n}+\overline{\alpha}_2(\phi_{i-1,j}^{n}+\phi_{i+1,j}^{n}+\phi_{i,j-1}^{n}+\phi_{i,j+1}^{n})+\overline{\beta}_1\phi_{i,j}^{n-1}+\overline{\beta}_2(\phi_{i-1,j}^{n-1}+\phi_{i+1,j}^{n-1}+\phi_{i,j-1}^{n-1}+\phi_{i,j+1}^{n-1})\\
		&\quad+\overline{\beta}_3(\phi_{i-1,j-1}^{n-1}+\phi_{i-1,j+1}^{n-1}+\phi_{i+1,j-1}^{n-1}+\phi_{i+1,j+1}^{n-1})+\overline{\gamma}_1\phi_{i,j}^{n-2}+\overline{\gamma}_2(\phi_{i-1,j}^{n-2}+\phi_{i+1,j}^{n-2}+\phi_{i,j-1}^{n-2}+\phi_{i,j+1}^{n-2})\\
		&\quad+\overline{\gamma}_3(\phi_{i-1,j-1}^{n-2}+\phi_{i-1,j+1}^{n-2}+\phi_{i+1,j-1}^{n-2}+\phi_{i+1,j+1}^{n-2})+\overline{\zeta}_1\phi_{i,j}^{n-3}
		+\overline{\zeta}_2(\phi_{i-1,j}^{n-3}+\phi_{i+1,j}^{n-3}+\phi_{i,j-1}^{n-3}+\phi_{i,j+1}^{n-3})\\
		&\quad+\overline{\eta} \phi_{i,j}^{n-4}
		+\Delta t\overline{\delta} R,
	\end{aligned}
\end{equation}
where the parameters $\overline{\alpha}_k, \overline{\zeta}_k\:(k=1,2)$, $\overline{\beta}_l,\overline{\gamma}_l\:(l=1,2,3)$, $\overline{\eta}$ and $\overline{\delta}$ are expressed as
\begin{subequations}\label{fdparameter2}
	\begin{align}
		&\overline{\alpha}_1=\frac{2}{2-\zeta \Delta t}\big[\alpha_1-\frac{\zeta\Delta t}{2} \big(w_0(s_4-2)  +1-s_4		\big)\big],\\			
		&\overline{\alpha}_2=\frac{2}{2-\zeta \Delta t}\big[\alpha_2-\frac{\zeta \Delta t}{8}\big(2-2s_2-s_5+w_0(2-s_4)	\big)\big],		\\
		&\overline{\beta}_1=\frac{2}{2-\zeta \Delta t}\big[\beta_1+\frac{\zeta \Delta t}{8}\big(4w_0(s_2 - 1)(s_4 - 2)+4s_5(s_2-1)\big)\big],\\
		&\overline{\beta}_2=\frac{2}{2-\zeta \Delta t}\big[\beta_2+\frac{\zeta\Delta t}{8}\big((s_4 - 1)(2s_2 + s_5 - 2)-w_0(s_4 - 2)(2s_2 + s_5 - 3)\big)\big],\\
		&\overline{\beta}_3=\frac{2}{2-\zeta \Delta t}\big[\beta_3+\frac{\zeta\Delta t}{8}\big(w_0(s_4 - 2)(s_2 + s_5 - 2)+s_2(s_2 + s_5 - 2)\big)
		\big],		\\
		&\overline{\gamma}_1=\frac{2}{2-\zeta \Delta t}
		\big[\gamma_1+\frac{\zeta\Delta t}{8}\big(	
		4s_5(s_2 - 1)(s_4 - 1)-4w_0(s_2 - 1)(s_4 - 2)(s_5 - 1)	
		\big)\big],\\
		&\overline{\gamma}_2=\frac{2}{2-\zeta \Delta t}\big[\gamma_2+\frac{\zeta\Delta t}{8}\big(			(s_2 - 1)(s_5 + s_2s_5 - 2)+w_0(s_2 - 1)(s_4 - 2)(s_2 + 2s_5 - 3)\big)\big],\\
		&\overline{\gamma}_3=\frac{2}{2-\zeta \Delta t}\big[\gamma_3+\frac{\zeta\Delta t}{8}\big(	s_2(s_4 - 1)(s_2 + s_5 - 2)-w_0(s_2 - 1)(s_4 - 2)(s_2 + s_5 - 2)\big],	\\
		&\overline{\zeta}_1=\frac{2}{2-\zeta \Delta t}\big[\zeta_1-\frac{\zeta\Delta t(1-s_2)}{8}\big(
		4(1-s_2)(1-s_5) +4w_0(s_2 - 1)(s_4 - 2)(s_5 - 1)
		\big)\big],\\
		&\overline{\zeta}_2=\frac{2}{2-\zeta\Delta t}\big[\zeta_2-\frac{\zeta\Delta t(1-s_2)}{8}\big((s_4 - 1)(s_5 + s_2s_5 - 2)-w_0(s_2 - 1)(s_4 - 2)(s_5 - 1)\big)\big],\\
		&\overline{\eta} =\frac{2}{2-\zeta\Delta t}[\eta+\frac{\zeta\Delta t}{2}(s_2 - 1)^2(s_4 - 1)(s_5 - 1)],\overline{\delta} =\frac{2}{2-\zeta \Delta t}\delta.
	\end{align}
\end{subequations}
Here the parameters $\alpha_k, \zeta_k\:(k=1,2)$, $\beta_l,\gamma_l\:(l=1,2,3)$, $\eta$ and $\delta$ are same as those in Eq. (\ref{fdparameter}). Through the Taylor expansion anlysis, one can obtain an explicit SLFD-AL scheme with the second-order accuracy in time and fourth-order accuracy [at the diffusive scaling ($\Delta t \propto \Delta x^2$)] in space once the following relations are satisfied, 
\begin{align}\label{solution2}
	\left\{\begin{aligned}
		 &w_0=1-4\varepsilon,\\
		&s_2=1,\\
		&s_4=\frac{6-24\varepsilon }{5-12\varepsilon},\\
		&s_5=\frac{6}{5}.\\
	\end{aligned}\right.\end{align} 
In particular, if we introduce a modified relaxation parameter $\overline{s}_2$ \cite{Silva2023} to satisfy the following relation,
\begin{equation}\label{s2}
	\frac{1-w_0}{2}\big[\big(\frac{1}{\overline{s}_2}-
	\frac{1}{2}\big)-\big(\frac{1}{\overline{s}_2^2}-\frac{1}{\overline{s}_2}\big)\zeta\Delta t\big]=\frac{1-w_0}{2}(\frac{1}{s_2}-\frac{1}{2}),
\end{equation}
then one can determine $s_2$ from Eq. (\ref{s2}),
\begin{equation}\label{s22}
	s_2=\frac{2\Delta t\zeta}{\Delta t\zeta - \sqrt{(\overline{s}_2 - 4\Delta t\zeta + \Delta t^2\overline{s}_2\zeta^2 + 2\Delta t\overline{s}_2\zeta)/ \overline{s}_2} + 1}.
\end{equation}
Based on Eq. (\ref{s22}), one can obtain another one fourth-order SLFD-AL scheme when the following conditions are satisfied,
\begin{align}\label{solution3}
	\left\{\begin{aligned}
		&w_0=\frac{s_5+(6s_5-12)\varepsilon}{s_5},\\
		&\overline{s}_2=\frac{6s_5-12}{s_5-6},\\
		&s_4=\frac{12\varepsilon s_5^2-24\varepsilon s_5+2s_5^2}{2s_5-36\varepsilon+24\varepsilon s_5+\varepsilon s_5^2}.
	\end{aligned}\right.
\end{align}

Now, we give four remarks on above results.
\begin{remark} It should be noted that the parameter $\overline{s}_2$ is just a modified relaxation parameter, and the real relaxtion parameter $s_2$ in Eq. (\ref{s22}) is dependent on the paramter $\zeta$ and the time step $\Delta t$. Thus, in the implementation of the MRT-LB-AL model with the condition (\ref{solution3}), one should first determine the parameters $w_0,s_5$ and the modified relaxation parameter $\overline{s}_2$ in terms of $s_4$, then the relaxation parameter $s_2$ in the diagonal relaxation matrix $\bm{S}$ can be obtained from Eq. (\ref{s22}).
\end{remark}
\begin{remark} If $\zeta=0$, Eqs. (\ref{distributionf2}) and (\ref{fdparameter2}) would reduce to Eqs. (\ref{distributionf}) and (\ref{fdparameter}), and the MRT-LB-AL model and SLFD-AL scheme are the same to the MRT-LB model (\ref{distributionf}) and SLFD scheme (\ref{fd}), as expected. However, it should be noted that the conditions (\ref{solution2}), (\ref{solution3}) and (\ref{solution}) to ensure the fourth-order accuracy of the MRT-LB model do not have such a relation, i.e., the condition (\ref{solution2}) or (\ref{solution3}) cannot reduce to Eq. (\ref{solution}) when $\zeta=0$. Actually, one can show that the condition (\ref{solution2}) is a special case of the condition (\ref{solution}).
\end{remark}\begin{remark} In terms of the stability analyis, simliar to the proof in the part \ref{stability-analysis}, one can also show that both the MRT-LB-AL model and SLFD-AL scheme are unconditionally stable.
\end{remark}
\begin{remark}
If the transform matrix $\bm{M_N}$ and diagonal relaxation matrix $\bm{S}$ are considered, the number of the conditions to ensure the MRT-LB-AL to be fourth-order accurate is less than that of the free parameters $w_0,s_2,s_4$ and $s_5$. Thus, for the diffusion equation with a linear source term (\ref{equlinear}), the MRT-LB model with transform matrix $\bm{M_N}$ and diagonal relxation matrix $\bm{S}$ cannot achieve a fourth-order accuracy unless one introduces a modified relaxation parameter $\overline{s}_2$ like Eq. (\ref{s22}), and in this case, the fourth-order conditions are given by
\begin{align}\label{solution4}
	\left\{\begin{aligned}
				&w_0=1-4\sqrt{3}\varepsilon,\\
		&s_2=\frac{2\Delta t\zeta}{\Delta t\zeta - \sqrt{(\overline{s}_2 - 4\Delta t\zeta + \Delta t^2\overline{s}_2\zeta^2 + 2\Delta t\overline{s}_2\zeta)/\overline{s}_2} + 1},\\
		&s_4= s_5=\frac{6}{2\sqrt{3}+3}.\\
	\end{aligned}\right.
\end{align}
where $\overline{s}_2=6/(3+\sqrt{3})$. Actually, one can observe that the difference between Eqs. (\ref{solution4}) and (\ref{solutionm}) is just the expression of $s_2$.
\end{remark}
In the following, we name the MRT-LB-AL model and the SLFD-AL scheme with the conditions (\ref{solution2}) and (\ref{solution3}) as the MRT-LB-AL1 model, MRT-LB-AL2 model, SLFD-AL1 scheme and SLFD-AL2 scheme. To test the two MRT-LB-AL models and SFLD-AL schemes for the diffusion equation (\ref{equlinear}), we consider an example with the following initial and periodic boundary conditions,
\begin{equation}\label{Ex3}
	\begin{aligned}
		\left\{
		\begin{array}{lr}
			\partial_t\phi=\kappa(\partial_x^2\phi+\partial_y^2\phi)-\pi^2(1-\kappa)\phi ,  (x,y)\in [-1,1]\times[-1,1],&  \\
			\phi(x,y,t=0)=\sin(\pi x)\sin(\pi y) ,  &  \\
		\end{array}
		\right.
	\end{aligned}
\end{equation}
and obtain the analytical solution as
\begin{equation}\label{ExSo1}
	\phi(x,y,t)=\sin(\pi x)\sin(\pi y)\exp[-\pi^2(\kappa+1)t].
\end{equation}

We conduct a number of simulations with some specified values of  discretization parameter $\varepsilon$, and calculate the $RMSEs$ and $CRs$ under different lattice sizes. As shown in Tables \ref {TEx3CRLB}, \ref {TEx4CRLB}, \ref {TEx3CRFD} and \ref{TEx4CRFD} where the simulations are suspended at time $T=2.0$ and the lattice spacing is varied from $\Delta x/8$ to $\Delta x=1/10$ with a fixed $\Delta x^2/\Delta t$. The results show that the two MRT-LB-AL models and SLFD-AL schemes have a fourth-order convergence rate in space.

\begin{table}
	\begin{center}
		\caption{The $RMSEs$ and $CRs$ of MRT-LB-AL1 model under different values of discretization parameter $\varepsilon$ ($T=2.0$, $\Delta x=1/10$, $\Delta t=1/20$).}
		\label{TEx3CRLB}  
		\begin{tabular}{ccccccccc}
			\hline \hline
			$\varepsilon$&$RMSE_{\Delta x,\Delta t}$&$RMSE_{\Delta x/2,\Delta t/4}$&
			$RMSE_{\Delta x/4,\Delta t/16}$&
			$RMSE_{\Delta x/8,\Delta t/64}$&$CR$\\
			\midrule[1pt]
			0.05 &  $2.4698\times 10^{-10}$ &$2.5765\times 10^{-11}$ &$1.6454\times 10^{-12}$&$1.0354\times 10^{-13}$&$\sim 3.9035$ \\
			
			0.08&$3.0363\times 10^{-10}$&$2.2482\times 10^{-11}$&$1.4355\times 10^{-12}$&$9.0325\times 10^{-14}$&$\sim 3.9050$\\
			
			0.10&$2.7737\times 10^{-10}$&$2.0499\times 10^{-11}$&$1.3088\times 10^{-12}$&$8.2354\times 10^{-14}$&$\sim 3.9059$\\
			
			0.12&$2.5291\times 10^{-10}$&$1.8661\times 10^{-11}$&$1.1915\times 10^{-12}$&$7.4974\times 10^{-14}$&$\sim 3.9067$\\
			\hline \hline
		\end{tabular}
	\end{center}
\end{table}

\begin{table}
	\begin{center}
		\caption{The $RMSEs$ and $CRs$ of SLFD-AL1 scheme under different values of discretization parameter $\varepsilon$ ($T=2.0$, $\Delta x=1/10$, $\Delta t=1/20$).}
		\label{TEx3CRFD}  
		\begin{tabular}{ccccccccc}
			\hline \hline
			$\varepsilon$&$RMSE_{\Delta x,\Delta t}$&$RMSE_{\Delta x/2,\Delta t/4}$&
			$RMSE_{\Delta x/4,\Delta t/16}$&
			$RMSE_{\Delta x/8,\Delta t/64}$&$CR$\\
			\midrule[1pt]
			0.05 &  $3.3591\times 10^{-10}$ &$2.5532\times 10^{-11}$ &$1.6435\times 10^{-12}$&$1.0363\times 10^{-13}$&$\sim  3.8875$ \\
			
			0.08&$2.9279\times 10^{-10}$&$2.2234\times 10^{-11}$&$1.4316\times 10^{-12}$&$9.0274\times 10^{-14}$&$\sim 3.8878$\\
			
			0.10&$2.6655\times 10^{-10}$&$2.0245\times 10^{-11}$&$1.3040\times 10^{-12}$&$8.2232\times 10^{-14}$&$\sim 3.8875$\\
			
			0.12&$2.4192\times 10^{-10}$&$1.8402\times 10^{-11}$&$1.1859\times 10^{-12}$&$7.4800\times 10^{-14}$&$\sim 3.8864$\\
			\hline \hline
		\end{tabular}
	\end{center}
\end{table}

\begin{table}
	\begin{center}
		\caption{The $RMSEs$ and $CRs$ of MRT-LB-AL2 model under different values of discretization parameter $\varepsilon$ ($T=2.0$, $\Delta x=1/10$, $\Delta t=1/20$).}
		\label{TEx4CRLB}  
		\begin{tabular}{ccccccccc}
			\hline \hline
			$\varepsilon$&$RMSE_{\Delta x,\Delta t}$&$RMSE_{\Delta x/2,\Delta t/4}$&
			$RMSE_{\Delta x/4,\Delta t/16}$&
			$RMSE_{\Delta x/8,\Delta t/64}$&$CR$\\
			\midrule[1pt]
			0.05 &  $3.3623\times 10^{-10}$ &$2.4941\times 10^{-11}$ &$1.6198\times 10^{-12}$&$1.0859\times 10^{-13}$&$\sim 3.8655$ \\
			
			0.08&$3.0236\times 10^{-10}$&$2.2327\times 10^{-11}$&$1.4255\times 10^{-12}$&$8.9788\times 10^{-14}$&$\sim 3.9053$\\
			
			0.10&$2.8705\times 10^{-10}$&$2.1531\times 10^{-11}$&$1.3899\times 10^{-12}$&$9.0262\times 10^{-14}$&$\sim 3.8783$\\
			
			0.12&$2.7782\times 10^{-10}$&$2.1517\times 10^{-11}$&$1.4472\times 10^{-12}$&$1.0597\times 10^{-13}$&$\sim 3.7854$\\
			\hline \hline
		\end{tabular}
	\end{center}
\end{table}

\begin{table}
	\begin{center}
		\caption{The $RMSEs$ and $CRs$ of SLFD-AL2 scheme under different values of discretization parameter $\varepsilon$ ($T=2.0$, $\Delta x=1/10$, $\Delta t=1/20$).}
		\label{TEx4CRFD}  
		\begin{tabular}{ccccccccc}
			\hline \hline
			$\varepsilon$&$RMSE_{\Delta x,\Delta t}$&$RMSE_{\Delta x/2,\Delta t/4}$&
			$RMSE_{\Delta x/4,\Delta t/16}$&
			$RMSE_{\Delta x/8,\Delta t/64}$&$CR$\\
			\midrule[1pt]
			0.05 &  $3.1243\times 10^{-10}$ &$2.4367\times 10^{-11}$ &$1.6078\times 10^{-12}$&$1.0823\times 10^{-13}$&$\sim 3.8317$ \\
			
			0.08&$2.7863\times 10^{-10}$&$2.1798\times 10^{-11}$&$1.4165\times 10^{-12}$&$8.9619\times 10^{-14}$&$\sim 3.8674$\\
			
			0.10&$2.6390\times 10^{-10}$&$2.1049\times 10^{-11}$&$1.3848\times 10^{-12}$&$9.0349\times 10^{-14}$&$\sim 3.8374$\\
			
			0.12&$2.5489\times 10^{-10}$&$2.1091\times 10^{-11}$&$1.4472\times 10^{-12}$&$1.0641\times 10^{-13}$&$\sim 3.7420$\\
			\hline \hline
		\end{tabular}
	\end{center}
\end{table}	
\section{The derivation of Eq. (\ref{2neq})}\label{stability}
Based on the Routh-Hurwitz stability criterion \cite{Routh, Hurwitz, Gantmacher}, $q(\lambda)$ (\ref{qcharc}) is the von Neumann polynomial when the following conditions are satisfied,
\begin{subequations}\label{neq}
	\begin{align}
		p_1-p_2+p_3-p_4>0,\label{neq1}\\
		3p_1+3p_4-p_2-p_3>0,\label{neq2}\\
		3p_1-3p_4-p_3+p_2>0,\label{neq3}\\
		p_1+p_2+p_3+p_4>0,\label{neq4}\\		
		p_1^2 - p_1p_3 - p_4^2 + p_2p_4>0.
	\end{align}
\end{subequations}
Now, we prove the equivalence of Eqs. (\ref{neq}) and (\ref{2neq}). On the one hand, we can derive Eq. (\ref{2neq}) from Eq. (\ref{neq}) after the summations of Eqs. (\ref{neq1}) and (\ref{neq3}), Eqs. (\ref{neq2}) and (\ref{neq4}).  On the other hand, to derive Eq. (\ref{neq}) from Eq. (\ref{2neq}), we first have
\begin{equation}
	\begin{aligned}
		(p_1+p_2+p_3+p_4)+(p_1-p_2+p_3-p_4)=2p_1+2p_3>0,	
	\end{aligned}
\end{equation}
then under the condition of $p_1>0$ [see Eq. (\ref{p1leq0})] and Eq. (\ref{neq25}), one can obtain
\begin{equation}\label{a21}
	p_1\times(p_1+p_3)+(p_1^2-p_1p_3-p_4^2+p_2p_4)=2p_1^2-p_4^2+p_2p_4>0.
\end{equation} 
Now we consider following three cases,
\begin{subequations}
	\begin{align}
		&	{\rm Case\:1:}\:p_2<\frac{p_4^2-2p_1^2}{p_4},-p_1< p_4<0,\label{D41}\\
		&	{\rm Case\:2:}\:p_2>\frac{p_4^2-2p_1^2}{p_4},0<p_4< p_1,\label{D42}\\
		&	{\rm Case\:3:}\:p_4=0,-p_1<p_4< p_1,\label{D43}
	\end{align}
\end{subequations}
where the conditions (\ref{neq23}) and (\ref{neq24}) have been used. From Eqs. (\ref{D41}) and (\ref{D42}), one can derive
\begin{equation}\label{con}
	\begin{aligned}
		&p_1-p_2>0,\\
		&p_1+p_2>0,
	\end{aligned}
\end{equation}
then we have the following results,
\begin{subequations}\label{D6}
	\begin{align}
	3p_1+3p_4-p_2-p_3&=\frac{ p_1\big(3p_1+3p_4-p_2-(p_3)\big)}{p_1}>\frac{ 3p_1^2+3p_1p_4-p_1p_2+(p_4^2-p_1^2-p_2p_4)}{p_1},\\
		&=\frac{ 2p_1^2+3p_1p_4-p_1p_2-p_4^2-p_2p_4}{p_1}=\frac{(p_1+p_4)^2+ (p_1-p_2)(p_1+p_4)}{p_1}>0,\\
	3p_1-3p_4-p_3+p_2&=\frac{ p_1\big(3p_1-3p_4-(p_3)+p_2\big)}{p_1}>\frac{ 3p_1^2-3p_1p_4+p_1p_2+(p_4^2-p_1^2-p_2p_4)}{p_1},\\
		&	=\frac{ 2p_1^2-3p_1p_4+p_1p_2-p_4^2-p_2p_4}{p_1}=\frac{(p_1-p_4)^2+ (p_1+p_2)(p_1-p_4)}{p_1}>0.
	\end{align}
\end{subequations}
For the Case 3, under the condition of Eqs. (\ref{neq21}), (\ref{neq22}) and (\ref{neq25}), we have
\begin{subequations}
	\begin{align}
		&-(p_1+p_3)<p_2<p_1+p_3,\\
		&p_1-p_3>0,
		\end{align}
\end{subequations}
then one can obtain,
\begin{subequations}\label{D8}
	\begin{align}
		3p_1+3p_4-p_2-p_3&=3p_1-(p_2)-p_3>3p_1-(p_1+p_3)-p_3=2(p_1-p_3)>0,\\
		3p_1-3p_4-p_3+p_2&=3p_1-p_3+(p_2)>3p_1-p_3+(-p_1-p_3)=2(p_1-p_3)>0.
	\end{align}
\end{subequations}
Thus, from Eqs. (\ref{D6}) and (\ref{D8}), we prove the conditions (\ref{neq2}) and (\ref{neq3}).
\section{The discussion on Eq. (\ref{neq25}) }\label{A4}
Similar to the way to prove $K(s_4=0)\geq0$ in the closed region $\overline{\Omega}_{s_4=0}$ (see the part \ref{stability-analysis}), we now consider the other seven cases of $K(s_4=2),K(s_2=0),K(s_2=2),K(w_0=0),K(w_0=1),K(\chi_1=-1)$ and $K(\chi_1=1)$. 
Through some symbolic manipulations, one can show that they all have no stationary points in their open regions $\Omega_{s_4=2},\Omega_{s_2=0},\Omega_{s_2=2},\Omega_{w_0=0},\Omega_{w_0=1},\Omega_{\chi_1=-1}$ and $\Omega_{\chi_1=1}$. Therefore, we only need to test whether the values of the above seven functions at the boundaries (i.e.,  $\partial\overline{\Omega}_{s_4=2},
\partial\overline{\Omega}_{s_2=0},
\partial\overline{\Omega}_{s_2=2},
\partial\overline{\Omega}_{w_0=0},
\partial\overline{\Omega}_{w_0=1},
\partial\overline{\Omega}_{\chi_1=-1},
\partial\overline{\Omega}_{\chi_1=1}$) are no less than zero. \\
\textbf{Case 1. The proof of $\bm{K{(s_4=2)}}\geq0$ at $\partial\overline{\Omega}_{s_4=2}$}
\begin{subequations}
	\begin{align}
		&K(s_4=2,s_2=0)=\underbrace{(\chi_1+\chi_2)^2(w_0-1)(\chi_1+\chi_2-2)}_{\geq0}\times\underbrace{\frac{w_0+\chi_1+\chi_2+\chi_1\chi_2(1-w_0)+1}{2}}_{\geq0}\:\:{\rm in}\:\:\overline{\Omega}_{s_4=2,s_2=0},\\
		&K(s_4=2,s_2=2)=\underbrace{-w_0(\chi_1-\chi_2)^2(\chi_1+\chi_2-2)}_{\geq0}\times\underbrace{\frac{\chi_1+\chi_2+2+w_0(\chi_1\chi_2-1)}{2}}_{\geq0}\:\:{\rm in}\:\:\overline{\Omega}_{s_4=2,s_2=2},\\
		&K(s_4=2,w_0=0)=\underbrace{(2-s_2)\big(1-\frac{\chi_1+\chi_2}{2}(s_2-1)^2\big)}_{\geq0}\times\underbrace{\big(\chi_1\chi_2(1-s_2)+\frac{\chi_1+\chi_2}{2}(2-s_2)+1\big)}_{\geq0}K_{21}\:\:{\rm in}\:\:\overline{\Omega}_{s_4=2,w_0=0},\\
		&K(s_4=2,w_0=1)=s_2\underbrace{\big(\chi_1\chi_2(s_2-1)+\frac{\chi_1+\chi_2}{2}s_2+1\big)}_{\geq0}\times\underbrace{\big(\frac{\chi_1+\chi_2}{2}(s_2-1)^2-1\big)}_{\leq0}K_{22}\:\:{\rm in}\:\:\overline{\Omega}_{s_4=2,w_0=1},\\		
		&K(s_4=2,\chi_1=-1)=\underbrace{(1+\chi_2)\frac{3-s_2(2-s_2)-\chi_2(s_2-1)^2}{16})}_{\geq0}K_{23}K_{24}\:\:{\rm in}\:\:\overline{\Omega}_{s_4=2,\chi_1=-1},\\
		&K(s_4=2,\chi_1=1)=\underbrace{\frac{s_2(2-s_2)-\chi_2(s_2-1)^2+1}{16}}_{\geq0}K_{25}K_{26}\:\:{\rm in}\:\:\overline{\Omega}_{s_4=2,\chi_1=1},
	\end{align}
\end{subequations}
where 
\begin{subequations}\label{A41}
	\begin{align}
		&K_{21}=s_2+\chi_1\chi_2(s_2-2)(1-s_2)+\frac{(\chi_1+\chi_2)^2}{4}(2-4s_2-s_2^4+s_2^2+2s_2^3)\:\:{\rm in}\:\:\overline{\Omega}_{s_4=2,s_2=0},\\
		&K_{22}=s_2-2+\chi_1\chi_2s_2(s_2-1)+\frac{(\chi_1+\chi_2)^2}{4}(2-8s_2+s_2^4-6s_2^3+11s_2^2)\:\:{\rm in}\:\:\overline{\Omega}_{s_4=2,s_2=2},\\
		&K_{23}=(s_2-2w_0)(1-s_2)+2+\chi_2(1-s_2)(2-s_2-2w_0)\:\:{\rm in}\:\:\overline{\Omega}_{s_4=2,w_0=0},\\
		&K_{24}=\chi_2(2s_2+4w_0-2s_2^2-s_2^3+s_2^4+18s_2^2w_0-6s_2^3w_0-16s_2w_0)\notag\\
		&\quad\quad+2s_2+4w_0-6s_2^2w_0+2s_2^3w_0-2s_2^2+3s_2^3-s_2^4\:\:{\rm in}\:\:\overline{\Omega}_{s_4=2,w_0=1},\\
		&K_{25}=\chi_2^2(-3s_2^2-2s_2w_0-4+2s_0+7s_2)+\chi_2(6s_2-4w_0+4s_2w_0-4s_2^2)\notag\\
		&\quad\quad+3s_2+2w_0-2s_2w_0-s_2^2+4\:\:{\rm in}\:\:\overline{\Omega}_{s_4=2,\chi_1=-1},\\
		&K_{26}=\chi_2\big(4(w_0-1)+6s_2+s_2^4+6s_2^2w_0-2s_2^3w_0-8s_2w_0\big)+\chi_2^2(-6w_0-3s_2^2)
		\notag\\
		&\quad\quad+8s_2w_0-4w_0-6s_2+2s_2^3w_0+8s_2^2-5s_2^3+s_2^4+4\:\:{\rm in}\:\:\overline{\Omega}_{s_4=2,\chi_1=1}.
	\end{align}
\end{subequations}
\textbf{Case 2. The proof of $\bm{K{(s_2=0)}}\geq0$ at $\partial \overline{\Omega}_{s_2=0}$}
\begin{subequations}
	\begin{align}
		&K(s_2=0,w_0=0)=2\underbrace{(\frac{\chi_1+\chi_2}{2}-\chi_1\chi_2)}_{\geq0}\underbrace{\big(1+\frac{\chi_1+\chi_2}{2}(1-s_4)\big)}_{\geq0}K_{31}\:\:{\rm in}\:\:\overline{\Omega}_{s_2=0,w_0=0},\\
		&K(s_2=0,w_0=1)=(2-s_4)\underbrace{\big(1+\frac{\chi_1+\chi_2}{2}(1-s_4)\big)}_{\geq0}\underbrace{(1+\chi_1\chi_2-\chi_1-\chi_2)}_{\geq0}K_{32}\:\:{\rm in}\:\:\overline{\Omega}_{s_2=0,w_0=1},\\
		&K(s_2=0,\chi_1=-1)=(1+\chi_2)\underbrace{\frac{\big(2-s_4(1-w_0\big)\big(2s_4+(1-s_4)(1+\chi_2)\big)}{8}}_{\geq0}K_{33}\:\:{\rm in}\:\:\overline{\Omega}_{s_2=0,\chi_1=-1},\\
		&K(s_2=0,\chi=1)=s_4(w_0-1)(\chi_2-1)^2\underbrace{\frac{s_4-2+\chi_2(s_4-1)-1}{8}}_{\leq0}K_{34}\:\:{\rm in}\:\:\overline{\Omega}_{s_2=0,\chi_1=1},
	\end{align}
\end{subequations}
where 
\begin{subequations}\label{A42}
	\begin{align}
		&K_{31}=(1+2\chi_1\chi_2)(1-s_4)^2+\frac{\chi_1+\chi_2}{2}(s_4^3-3s_4^2+6s_4-4)+1\:\:{\rm in}\:\:\overline{\Omega}_{s_2=0,w_0=0},\\
		&K_{32}=s_4-\chi_1\chi_2(2-s_4)-\frac{(\chi_1+\chi_2)^2}{2}(s_4-1)\:\:{\rm in}\:\:\overline{\Omega}_{s_2=0,w_0=1},\\
		&K_{33}=\chi_2^2\big(-s_4^2(1+w_0)-4+s_4(w_0+5)\big)+2\chi_2s_4(1+2s_4w_0-5w_0)\notag\\
		&\quad\quad+s_4w_0(5-3s_4)+s_4(s_4-3)+4\:\:{\rm in}\:\:\overline{\Omega}_{s_2=0,\chi_1=-1},\\
		&K_{34}=\chi_2\big(s_4(3w_0-1)+(s_4-2+s_4^2(1-w_0)+s_4w_0\big)+s_4^2(1-w_0)+2\:\:{\rm in}\:\:\overline{\Omega}_{s_2=0,\chi_1=1}.
	\end{align}
\end{subequations}
\textbf{Case 3. The proof of $\bm{K{(s_2=2)}}\geq0$ at $\partial\overline{\Omega}_{s_2=2}$}
\begin{equation}
	K(s_2=2)=-s_4\underbrace{\big(1+(1-s_4)\frac{\chi_1+\chi_2}{2}\big)}_{\geq0}
	\underbrace{\big(w_0(\chi_1\chi_2-1)+2(1+\frac{\chi_1+\chi_2}{2})\big)}_{\geq0}K_{4}\:\:{\rm in}\:\:\overline{\Omega}_{s_2=2},
\end{equation}
where $K_4=s_4+\chi_1\chi_2s_4w_0-2+\frac{(\chi_1+\chi_2)^2}{4}(1-s_4)\big(2+s_4(w_0-1)\big)+s_4\frac{\chi_1+\chi_2}{2}(2-s_4)(1-w_0).$ Now, we focus on $K_4$ and obtain
\begin{subequations}\label{A45}
	\begin{align}
		&K_{4}(w_0=0)=\underbrace{(1-\frac{\chi_1+\chi_2}{2})(s_4-2)\big(\frac{\chi_1+\chi_2}{2}(1-s_4)+1\big)}_{\leq0}\:\:{\rm in}\:\:\overline{\Omega}_{s_2=2,w_0=0},\\
		&K_{4}(w_0=1)=\underbrace{s_4-2+\chi_1\chi_2s_4+\frac{(\chi_1+\chi_2)^2}{2}(1-s_4)}_{\leq0}\:\:{\rm in}\:\:\overline{\Omega}_{s_2=2,w_0=1},\\
		&K_4(\chi_1=-1)=s_4\big(1+\frac{\chi_2-1}{2}(2-s_4)(1-w_0)\big)-s_4w_0\chi_2-2\notag\\
		&\qquad\qquad\qquad-\frac{\big(s_4(w_0-1)+2\big)(s_4-1)(\chi_2-1)^2}{4}\:\:{\rm in}\:\:\overline{\Omega}_{s_2=2,\chi_1=-1},\label{K41}\\
		&K_4(\chi_1=1)=s_4\big(1+\frac{\chi_2+1}{2}(2-s_4)(1-w_0)\big)+s_4w_0\chi_2-2\notag\\
		&\qquad\qquad\quad\:\:-\frac{\big(s_4(w_0-1)+2\big)(s_4-1)(\chi_2+1)^2}{4}\:\:{\rm in}\:\: \overline{\Omega}_{s_2=2,\chi_1=1}.\label{K42}
	\end{align}
\end{subequations}
\textbf{Case 4. The proof of $\bm{K{(w_0=0)}}\geq0$ at $\partial\overline{\Omega}_{w_0=0}$}
\begin{subequations}
	\begin{align}
		&K(w_0=0,\chi_1=-1)=\underbrace{\frac{\big(1+(s_2-1)(s_4-1)\big)\big(s_2(1-\chi_2)+1+\chi_2\big)}{16}}_{\geq0} K_{51}K_{52}K_{53}\:\:{\rm in}\:\: \overline{\Omega}_{w_0=0,\chi_1=-1},\\
		&K(w_0=0,\chi_1=1)=\underbrace{\frac{\big(1-(1-s_2)(1-s_4)\big)\big(s_2(1+\chi_2)+1-\chi_2\big)}{16}}_{\geq0}K_{54}K_{55}K_{56}\:\:{\rm in}\:\: \overline{\Omega}_{w_0=0,\chi_1=1},
	\end{align}
\end{subequations}
where\begin{subequations}\label{A46}
	\begin{align}
		&K_{51}=2(s_2-1)(1+\chi_2)+s_2^2(\chi_2-3)+s_2^3(1-\chi_2)\:\:{\rm in}\:\: \overline{\Omega}_{w_0=0,\chi_1=-1},\\
		&K_{52}=1+\chi_2+(1-\chi_2)\big(s_2(2-s_2)(1-s_4)+s_4\big)\:\:{\rm in}\:\: \overline{\Omega}_{w_0=0,\chi_1=-1},\\
		&K_{53}=s_2+s_4(3-s_4)-4+s_2s_4(s_4-2)\notag\\
		&\qquad\:\:+\chi_2\big(4-3s_2+s_4(s_4-5)-s_2s_(s_4-4)\big)\:\:{\rm in}\:\: \overline{\Omega}_{w_0=0,\chi_1=-1},,\\
		&K_{54}=2(s_2-1)(1-2\chi_2)-s_2^2(3+\chi_2)+s_2^3(1+\chi_2)\:\:{\rm in}\:\: \overline{\Omega}_{w_0=0,\chi_1=1},\\
		&K_{55}=\chi_2\big(s_2-2+s_4+s_4^2(1-s_2)\big)+s_2s_4(2-s_4)-s_4(1-s_4)+2-s_2\:\:{\rm in}\:\: \overline{\Omega}_{w_0=0,\chi_1=1},\\
		&K_{56}=s_2(2-s_2)(1-s_4)+s_4-3+\chi_2(s_4-1)(s_2-1)^2\:\:{\rm in}\:\: \overline{\Omega}_{w_0=0,\chi_1=1}.
	\end{align}
\end{subequations}
\textbf{Case 5. The proof of $\bm{K{(w_0=1)}}\geq0$ at $\partial \Omega_{w_0=1}$}
\begin{subequations}
	\begin{align}
		&K(w_0=1,\chi_1=-1)=(2-s_2)\underbrace{\frac{\big(2+(\chi_2-1)(1-s_2)(1-s_4)\big)}{16}}_{\geq0}K_{61}K_{62}K_{63}\:\:{\rm in}\:\:\overline{\Omega}_{w_0=1,\chi_1=-1},\\
		&K(w_0=1,\chi_1=1)=-s_2\underbrace{\frac{\big(2-(1+\chi_2)(1-s_4)\big)\big(s_2s_4(1+\chi_2)+(1-\chi_2)(2-s_2)\big)}{16}}_{\geq0} K_{64}K_{65}\:\:{\rm in}\:\:\overline{\Omega}_{w_0=1,\chi_1=1},
	\end{align}
\end{subequations}
where
\begin{subequations}\label{A47}
	\begin{align}
		&K_{61}=\chi_2\big(4-3s_2+(s_2-2)s_4\big)+s_2(1-s_4)+2s(s_4-2)\:\:{\rm in}\:\:\overline{\Omega}_{w_0=1,\chi_1=-1},\\
		&K_{62}=\chi_2\big(1+s_2(s_2-2)-s_4(s_2-1)^2\big)+s_2(2-s_2)+s_4(s_2-1)^2+1\:\:{\rm in}\:\:\overline{\Omega}_{w_0=1,\chi_1=-1},\\
		&K_{63}=\chi_2(1-s_2)\big(s_2^2(1-s_4)+2(s_2s_4-1)\big)+s_2(1-s_4)(2-3s_2+s_2^2)-2\:\:{\rm in}\:\:\overline{\Omega}_{w_0=1,\chi_1=-1},\\
		&K_{64}=\chi_1\big(s_4(s_2-1)^2+(2-s_2)s_2-1\big)+s_2(2-s_2)+s_4(s_2-1)^2-3\:\:{\rm in}\:\:\overline{\Omega}_{w_0=1,\chi_1=1},\\
		&K_{65}=(1+\chi_2)(s_4-1)(2+s_2^2-3s_2)+s_2(s_2-1)\chi_2+(3s_2-2)(s_2-1)+2\:\:{\rm in}\:\:\overline{\Omega}_{w_0=1,\chi_1=1}.
	\end{align}
\end{subequations}
It is obvious that the functions in Eqs. (\ref{A41}), (\ref{A42}), (\ref{A45}), (\ref{A46}) and (\ref{A47}) have at most three variables, after some simple analysis, we can obtain $K_{2i}\geq0\:(i=1,3,4,5,6)$, $K_{22}\leq0$, $K_{3i}\geq0\:(i=1,2,3,4)$, $K_{4}(\chi_1=\pm1)\leq0$, $K_{5i}\geq0\:(i=2,5)$, $K_{5i}\leq0\:(i=1,3,4,6)$, $K_{6i}\leq0\:(i=1,3,4)$ and $K_{6i}\geq0\:(i=2,5)$, 
then $K(s_4=2)\geq0$ in the closed region $\overline{\Omega}_{s_4=2}$, $K(s_2=0)\geq0$ in the closed region $\overline{\Omega}_{s_2=0}$, $K(s_2=2)\geq0$ in the closed region $\overline{\Omega}_{s_2=2}$, $K(w_0=0)\geq0$ in the closed region $\Omega$, and $K(w_0=1)\geq0$ in closed region $\overline{\Omega}_{w_0=1}$.\\
\textbf{Case 6. The proof of $\bm{K{(\chi_1=-1)}}\geq0$ at $\partial\overline{\Omega}_{\chi_1=-1}$}
\begin{subequations}
	\begin{align}
		&K(\chi_1=-1,\chi_2=-1)=s_2(2-s_2)(2-s_4)\big(s_2(2-s_2)+s_4(s_2-1)^2\big)\notag\\
		&\quad\quad\quad\quad\quad\quad\quad\quad\quad\big(1-(s_2-1)^2+s_4(s_2-2w_0)(s_2-1)\big)^2\geq0\:\:{\rm in}\:\: \overline{\Omega}_{\chi_1=-1,\chi_2=-1},\\
		&K(\chi_1=-1,\chi_2=1)=-\big(1+(s_2-1)^2\big)\big(s_2+s_4(1-w_0)(1-s_2)\big)\notag\\
		&\quad\quad\quad\quad\quad\quad\quad\quad\quad\big(s_2-2+s_4(1-w_0)(1-s_2)\big)\geq0\:\:{\rm in}\:\: \overline{\Omega}_{\chi_1=-1,\chi_2=1},
	\end{align}
\end{subequations}
it is obvious that $K(\chi_1=-1)\geq0$ in closed region $\overline{\Omega}_{\chi_1=-1}$.\\
\textbf{Case 7. The proof of $\bm{K{(\chi_1=1)}}\geq0$ at $\partial\overline{\Omega}_{\chi_1=1}$}
\begin{subequations}
	\begin{align}
		&K(\chi_1=1,\chi_2=-1)=K(\chi_1=-1,\chi_2=1)\geq0\:\:{\rm in}\:\:\overline{\Omega}_{\chi_1=1,\chi_2=-1},\\
		&K(\chi_1=1,\chi_2=1)=s_2^3s_4(s_2 - 2)(s_2 + s_4 - s_2s_4)^2\big((s_4-2)(s_2-1)^2 + s_2(s_2-2)\big)
		\geq0\:\:{\rm in}\:\:\overline{\Omega}_{\chi_1=1,\chi_2=1},
	\end{align}
\end{subequations}
thus $K(\chi_1=1)\geq0$ in closed region $\overline{\Omega}_{\chi_1=1}$.


\begin{thebibliography}{00}
	\bibitem{Chen1998} S. Chen and G. Doolen, Lattice Boltzmann method for fluid flows, Annu. Rev. Fluid Mech. 30, 329 (1998).
	
	\bibitem{Succi2001} S. Succi, The Lattice Boltzmann Equation for Fluid Dynamics and Beyond (Oxford University Press, Oxford, 2001).
	
	\bibitem{Aidun2010} C. K. Aidun and J. R. Clausen, Lattice-Boltzmann method for complex flows, Annu. Rev. Fluid Mech. 42, 439 (2010).
	\bibitem{Guo2013} Z. Guo and C. Shu, Lattice Boltzmann Method and Its Applications in Engineering (World Scientific Publishing, Singapore, 2013).
	
	\bibitem{Kruger2017} T. Kr$\rm{\ddot{u}}$ger, H. Kusumaatmaja, A. Kuzmin, O. Shardt, G. Silva, and E. M. Viggen, The Lattice Boltzmann Method: Principles and Practice (Springer, Switzerland, 2017).
	\bibitem{Wang2019} H. Wang, X. Yuan, H. Liang, Z. Chai, and B. Shi, A brief review of the phase-field-based lattice Boltzmann method for multiphase flows, Capillary 2, 33 (2019).
	
	\bibitem{Huber2010}C. Huber, B. Chopard, and M. Manga, A lattice Boltzmann model for coupled diffusion, J. Comput. Phys. 229(20), 7956 (2010).
	\bibitem{Ancona1994}M. G. Ancona, Fully-lagrangian and lattice-boltzmann methods for solving systems of conservation equations, J. Comput. Phys. 115(1), 107 (1994).
	\bibitem{Suga2010}S. Suga, An accurate multi-level finite difference scheme for 1D diffusion equations derived from the lattice Boltzmann method, J. Stat. Phys. 140(3), 494 (2010).
	\bibitem{Lin2022} Y. Lin, N. Hong, B. Shi and Z. Chai, Multiple-relaxation-time lattice Boltzmann model-based four-level finite-difference scheme for one-dimensional diffusion equations, Phys. Rev. E 104(1), 015312 (2021).
	\bibitem{Silva2023} G. Silva, Discrete effects on the source term for the lattice Boltzmann modelling of one-dimensional reaction–diffusion equations, Comput. Fluids 251, 105735 (2023).
	
	
	\bibitem{Van2000}R. G. M. Van der Sman and M. H. Ernst, Convection-diffusion lattice Boltzmann scheme for irregular lattices, J. Comput. Phys. 160(2), 766 (2000).
	\bibitem{Gin2005}I. Ginzburg, Equilibrium-type and link-type lattice Boltzmann models for generic advection and anisotropic-dispersion equation, Adv. Wat. Resour. 28(11), 1171 (2005).
	\bibitem{Rasin2005}I. Rasin, S. Succi, and W. Miller, A multi-relaxation lattice kinetic method for passive scalar diffusion, J. Comput. Phys. 206(2), 453 (2005).
	\bibitem{Shi2009}B. Shi and Z. Guo, Lattice Boltzmann model for nonlinear convection-diffusion equations, Phys. Rev. E 79(1), 016701 (2009).
	\bibitem{Chopard2009}B. Chopard, J. L. Falcone, and J. Latt, The lattice Boltzmann advection-diffusion model revisited, Eur. Phys. J. Spec. Top. 171(1), 245 (2009).
	\bibitem{Yoshida2010}H. Yoshida and M. Nagaoka, Multiple-relaxation-time lattice Boltzmann model for the convection and anisotropic diffusion equation, J. Comput. Phys. 229(20), 7774 (2010).
	\bibitem{Gin2013}I. Ginzburg, Multiple anisotropic collisions for advection-diffusion lattice Boltzmann schemes, Adv. Wat. Resour. 51, 381 (2013).
	\bibitem{Chai2010}Z. Chai and T. S. Zhao, Lattice Boltzmann model for the convection-diffusion equation, Phys. Rev. E 87(6), 063309 (2013).
	\bibitem{Chai2014}Z. Chai, B. Shi, and Z. Guo, A multiple-relaxation-time lattice Boltzmann model for general nonlinear anisotropic convection-diffusion equations, J. Sci. Comput. 69(1), 355 (2016).
	\bibitem{Jettestuen2016}O. Aursj${\rm \varnothing}$, E. Jettestuen, J. L. Vinningland, and A. Hiorth, An improved lattice Boltzmann method for simulating advective-difusive processes in fluids, J. Comput. Phys. 332, 363 (2017).
	\bibitem{Li2017}L. Li, Multiple-time-scaling lattice Boltzmann method for the convection-diffusion equation, Phys. Rev. E 99(6), 063301 (2019).
	\bibitem{Michelert2022}J. Michelet, M. M. Tekitek and M. Berthier, Multiple relaxation time lattice Boltzmann schemes for advection-diffusion equations with application to radar image processing, J. Comput. Phys. 471, 111612 (2022).
	\bibitem{Dellacherie2014}S. Dellacherie, Construction and analysis of lattice Boltzmann methods applied to a 1D convection-diffusion equation, Acta. Appl. Math. 131(1), 69 (2014).
	\bibitem{Cui2016}S. Cui, N. Hong, B. Shi, and Z. Chai, Discrete effect on
	the halfway bounce-back boundary condition of multiple-relaxation-time lattice Boltzmann model for convection-diffusion equations, Phys. Rev. E 93(4), 043311 (2016).
	\bibitem{Chen2023} Y. Chen, Z. Chai and B. Shi. Fourth-order multiple-relaxation-time lattice Boltzmann model and equivalent finite-difference scheme for one-dimensional convection-diffusion equations, submitted.
		
	\bibitem{Hirabayashi2001}M. Hirabayashi, Y. Chen, and H. Ohashi, The lattice BGK model for the Poisson equation, JSME Int. J. Ser. B. 44(1), 45 (2001).
	\bibitem{Chai2008}Z. Chai and B. Shi, A novel lattice Boltzmann model for the Poisson equation, Appl. Math. Model. 32(10), 2050 (2008).
	\bibitem{Chai2019}Z. Chai, H. Liang, R. Du, and B. Shi, A lattice Boltzmann model for two-phase flow in porous media, SIAM J. Sci. Comput. 41(4), B746 (2019).

	\bibitem{Li2012}Q. Li, Z. Zheng, S. Wang, and J. Liu, A multi-level finite difference scheme for one-dimensional Burgers equation derived from the lattice Boltzmann Method, J. Appl. Math. 2012, 1 (2012).
	
	\bibitem{Qian1992} Y. H. Qian, D. d$'$Humi$\grave{{\rm e}}$re, and P. Lallemand, Lattice BGK	models for Navier-Stokes equation, Europhys. Lett. 17(6), 479 (1992).
	\bibitem{Chai2020} Z. Chai and B. Shi, Multiple-relaxation-time lattice Boltzmann method for the Navier-Stokes and nonlinear convection-diffusion equations: Modeling, analysis, and elements, Phys. Rev. E 102(2), 023306 (2020).
	
	
	\bibitem{Chai-1-2014}Z. Chai, T. S. Zhao, Nonequilibrium scheme for computing the flux of the convection-diffusion equation in the framework of the lattice Boltzmann method, Phys. Rev. E 90(1), 013305 (2014).
	\bibitem{Lallemand2000}P. Lallemand, L.S. Luo, Theory of the lattice Boltzmann method: dispersion, dissipation, isotropy, Galilean invariance, and stability, Phys. Rev. E 61, 6546 (2000).
	\bibitem{Pan2006}C. Pan, L.-S. Luo, C. T. Miller, An evaluation of lattice Boltzbmann schemes for porous medium flow simulation, Comput. Fluids 35(8), 898 (2006).
	\bibitem{Luo2011}L.-S. Luo, W. Liao, X. Chen, Y. Peng, and W. Zhang, Numerics
	of the lattice Boltzmann method: Effects of collision models
	on the lattice Boltzmann simulations, Phys. Rev. E 83, 056710
	(2011).	
	
	
	
	
	\bibitem{Chapman1990} S. Chapman and T. G. Cowling, The Mathematical Theory of
	Nonuniform Gases (Cambridge University Press, Cambridge, 1970).
	\bibitem{Ikenberry1956}E. Ikenberry and C. Truesdell, On the pressures and the flux
	of energy in a gas according to Maxwell's kinetic theory, J. Ration. Mech. Anal. 5(1), 1 (1956).
	\bibitem{Yong2016} W.-A. Yong, W. Zhao, and L.-S. Luo, Theory of the lattice
	Boltzmann method: Derivation of macroscopic equations via
	the Maxwell iteration, Phys. Rev. E 93(3), 033310 (2016).
	\bibitem{Holdych2004}D. J. Holdych, D. R. Noble, J. G. Georgiadis, and R. O. Buckius,
	Truncation error analysis of lattice Boltzmann methods,
	J. Comput. Phys. 193(2), 595 (2004).
	\bibitem{Wagner2006} A. Wagner, Thermodynamic consistency of liquid-gas lattice
	Boltzmann simulations, Phys. Rev. E 74(5), 056703 (2006).
	
	
	\bibitem{GinzburgD2009}D. d$'$Humi$\grave{{\rm e}}$re and I. Ginzburg, Viscosity independent numerical errors for lattice Boltzmann models: From recurrence
	equations to magic collision numbers, Comput. Math. Appl. 58(5), 823 (2009).
	\bibitem{Gin2012}I. Ginzburg, Truncation errors, exact and heuristic stability
	analysis of two-relaxation-times lattice Boltzmann schemes for
	anisotropic advection-diffusion equation, Commun. Comput.
	Phys. 11(5), 1439 (2012).
	
	\bibitem{Dubois2008}F. Dubois, Equivalent partial differential equations of a lattice Boltzmann scheme, Comput. Math. Appl. 55(7), 1441 (2008).
	\bibitem{Dubois2009}F. Dubois, Third order equivalent equation of lattice Boltzmann scheme, Discret. Contin. Dyn. Syst. 23(1-2), 221 (2009).
	\bibitem{Dubois2019}F. Dubois, Nonlinear fourth order Taylor expansion of lattice Boltzmann schemes, Asymptotic Anal. 127(4), 297 (2022).
	\bibitem{Junk}M. Junk, A finite difference interpretation of the lattice
	Boltzmann method, Numer. Meth. Part. Diff. Equ. 17(4), 383
	(2001).
	\bibitem{Inamuro}T. Inamuro, A lattice kinetic scheme for incompressible viscous
	flows with heat transfer, Philos. Trans. R. Soc. Lond. A 360(1792), 477
	(2002).
	
	\bibitem{Fort1953}E. C. Du Fort, and S. P. Frankel, Stability conditions in the
	numerical treatment of parabolic differential equations, Math.
	Comput. 7(43), 135 (1953).		
	
	\bibitem{Kwok1993} Y.-K. Kwok, K.-K. Tam, Stability analysis of three-level difference schemes for initial-boundary problems for multidimensional convective-diffusion equations, Commun. Numer. Methods Eng. 9(7), 595 (1993).
	
	\bibitem{Fucik2021} R. Fu\v{c}\'{i}k, R. Straka, Equivalent finite difference and partial differential equations for the lattice Boltzmann method. Comput. Math. Appl. 90, 96 (2021).
	
	
	\bibitem{Fucik2023} R. Fu\v{c}\'{i}k, P. Eichler, J. Klinkovsky, R. Straka and T.  Oberhuber, Lattice Boltzmann Method Analysis Tool (LBMAT), Numer. Algorithms, https://doi.org/10.1007/s11075-022-01476-8 (2022).
	
	\bibitem{Bellotti2022} T. Bellotti, B. Graille and M. Massot, Finite difference formulation of any lattice Boltzmann scheme, Numer. Math. 152(1), 1 (2022).
	\bibitem{Bellotti-1-2022} T. Bellotti, Rigorous derivation of the macroscopic equations for the lattice Boltzmann method via the corresponding Finite Difference scheme, arXiv.org, (2022).
	
	
	\bibitem{Chai-1-2016}Z. Chai, C. Huang, B. Shi, and Z. Guo, A comparative study on the lattice Boltzmann models for predicting effective diffusivity of porous media, Int. J. Heat Mass Transf. 98, 687 (2016).
	
	\bibitem{Miller} J. J. H. Miller, On the location of zeros of certain classes of polynomials with applications to numerical analysis, J. Inst. Math. Appl. 8(3), 397 (1971).
	\bibitem{Routh}E. J. Routh, A Treatise of the Stability of a Given State of Motion (Macmillan, London, 1877).
	\bibitem{Hurwitz} A. Hurwitz, Ueber die Bedingungen, unter welchen eine Gleichung nur Wurzeln mit negativen reellen Theilen besitzt, Math. Ann. 46(2), 273 (1895).
	\bibitem{Gantmacher} F. R. Gantmacher, Applications of the Theory of Matrices
	(Interscience Publishers, New York, 1959).

	 
	\bibitem{Hou} S.-H. Hou, Classroom Note: A Simple Proof of the Leverrier-Faddeev Characteristic Polynominal Algorithm. SIAM Rev. 40(3), 706-709 (1998). 
\end{thebibliography}
\end{document}